\pdfoutput=1
\RequirePackage{silence}
\documentclass[a4paper,11pt]{amsart}
\usepackage[hmarginratio={1:1},vmarginratio={1:1},lmargin=60.0pt,tmargin=60.0pt]{geometry}

\usepackage[numbers]{natbib}

\allowdisplaybreaks

\usepackage[utf8]{inputenc}

\usepackage{amssymb,amsmath,amsthm,amsfonts,mathrsfs,bbm}
\usepackage[table]{xcolor}
\usepackage{graphicx}
\usepackage{mathtools}

\usepackage{ytableau}

\definecolor{spinach}{RGB}{46,139,87}
\definecolor{tomato}{RGB}{255,99,71}
\definecolor{pumpkin}{RGB}{224,180,80}

\definecolor{orchid}{RGB}{143,40,194}
\definecolor{lava}{RGB}{207,16,32}
\definecolor{mydarkblue}{RGB}{10,10,150}

\usepackage[yyyymmdd,hhmmss]{datetime}
\usepackage{todonotes}

\usepackage{xparse}

\usepackage{enumitem}
\setlist[enumerate]{itemsep=0.15cm,label=\emph{\upshape(\alph*)}}
\setlist[enumerate,2]{itemsep=0.15cm,label=\emph{\upshape(\roman*)}}

\let\emph\relax
\DeclareTextFontCommand{\emph}{\em}

\usepackage{array}
\newcolumntype{C}{>{$}c<{$}}

\DeclarePairedDelimiterX{\set}[1]{\{}{\}}{\setargs{#1}}
\NewDocumentCommand{\setargs}{>{\SplitArgument{1}{|}}m}{\setargsaux#1}
\NewDocumentCommand{\setargsaux}{mm}
{\IfNoValueTF{#2}{#1} {#1\,\delimsize|\,\mathopen{}#2}}

\usepackage[all]{xy}
\usepackage{tikz}
\usetikzlibrary{cd}
\usetikzlibrary{decorations}
\usetikzlibrary{decorations.markings}
\usetikzlibrary{decorations.pathreplacing}
\usetikzlibrary{decorations.pathmorphing}
\usetikzlibrary{arrows.meta,shapes,positioning,matrix,calc}
\usetikzlibrary{shapes.callouts}
\tikzset{anchorbase/.style={baseline={([yshift=-0.5ex]current bounding box.center)}},
tinynodes/.style={font=\tiny,text height=0.25ex,text depth=0.05ex},
smallnodes/.style={font=\scriptsize,text height=0.75ex,text depth=0.15ex},
crossline/.style={preaction={draw=white,line width=10.0pt,-},preaction={draw=black,line width=1.8pt,-}},
usual/.style={line width=2.0,color=black},
affine/.style={line width=2.0,color=tomato,densely dotted},
levi/.style={decoration={markings,post length=0.25mm,pre length=0.25mm,mark=at position #1 with {\node[circle,radius=0.5cm,inner sep=-3.5pt,color=tomato,fill=tomato,opacity=0.5]{};}},postaction={decorate}},
levi/.default=0.5,
levia/.style={decoration={markings,post length=0.25mm,pre length=0.25mm,mark=at position #1 with {\node[circle,radius=0.5cm,inner sep=-3.5pt,color=tomato,fill=tomato,opacity=0.25]{};}},postaction={decorate}},
levia/.default=0.5,
levib/.style={decoration={markings,post length=0.25mm,pre length=0.25mm,mark=at position #1 with {\node[circle,radius=0.5cm,inner sep=-3.5pt,color=tomato,fill=tomato,opacity=0.5/3]{};}},postaction={decorate}},
levib/.default=0.5,
}
\tikzstyle directed=[postaction={decorate,decoration={markings,mark=at position #1 with {\arrow[line width=0.5mm, black]{>}}}}]
\tikzstyle rdirected=[postaction={decorate,decoration={markings,mark=at position #1 with {\arrow[line width=0.5mm, black]{<}}}}]

\newcommand{\ie}{\text{i.e.}}
\newcommand{\eg}{\text{e.g.}}

\newcommand{\etc}{\text{etc.}}

\newcommand{\ver}{\text{verbatim}}
\newcommand{\vive}{\text{vice versa}}

\newcommand{\loccit}{\text{loc. cit.}}

\newcommand{\acts}{\centerdot}
\renewcommand{\dots}{\text{...}}
\renewcommand{\vdots}{\rotatebox{90}{\text{...}}}
\renewcommand{\ddots}{\raisebox{0.2cm}{\rotatebox{-40}{\text{...}}}}
\newcommand{\placeholder}{{}_{-}}

\newcommand{\mystrut}{\rule[-0.2\baselineskip]{0pt}{0.9\baselineskip}}

\newcommand{\vcirc}{\circ_{v}}
\newcommand{\hcirc}{\circ_{h}}
\newcommand{\munit}{\mathbbm{1}}

\newcommand{\setstuff}[1]{\mathrm{#1}}
\newcommand{\catstuff}[1]{\mathbf{#1}}
\newcommand{\functorstuff}[1]{\mathtt{#1}}

\newcommand{\obstuff}[1]{\mathtt{#1}}
\newcommand{\morstuff}[1]{\mathrm{#1}}

\newcommand{\idmor}{\morstuff{id}}

\newcommand{\End}{\setstuff{End}}
\newcommand{\Hom}{\setstuff{Hom}}

\newcommand{\C}{\mathbb{C}}
\newcommand{\Z}{\mathbb{Z}}

\newcommand{\F}{\mathbb{F}}
\newcommand{\N}{\mathbb{Z}_{\geq 0}}

\newcommand{\qpar}{q}
\newcommand{\vpar}{v}
\newcommand{\Zv}{\Z_{v}}
\newcommand{\qnum}[1]{[#1]}
\newcommand{\qfac}[1]{[#1]!}
\newcommand{\qbin}[2]{{\textstyle\genfrac{[}{]}{0pt}{}{#1}{#2}}}
\newcommand{\qbinn}[2]{\genfrac{[}{]}{0pt}{}{#1}{#2}}

\newcommand{\K}{\mathbb{K}_{q}}
\newcommand{\Kone}{\mathbb{K}_{1}}
\newcommand{\Kstar}{\mathbb{K}_{\qpar}}

\newcommand{\U}{\mathbb{U}}
\newcommand{\upar}[1][n]{u_{#1}}

\newcommand{\epoly}[1][k]{e_{#1}}
\newcommand{\efunction}[1][k]{\tilde{e}_{#1}}
\newcommand{\epolyg}[2]{e_{#1}^{(#2)}}

\newcommand{\pivo}{\ast}
\newcommand{\danti}{\updownarrow}
\newcommand{\dinvo}{\leftrightarrow}

\newcommand{\completion}[1]{\functorstuff{AdId}(#1)}

\newcommand{\GLN}[1][n]{\mathrm{GL}_{#1}}
\newcommand{\Levi}{\mathrm{L}}
\newcommand{\sym}[1][m]{\mathrm{S}_{#1}}
\newcommand{\EXT}[1][\bullet]{{\textstyle\bigwedge^{#1}}}

\newcommand{\gln}[1][n]{\mathfrak{gl}_{#1}}
\newcommand{\sln}{\mathfrak{sl}_{n}}

\newcommand{\Uglnv}[1][\gln]{U_{\vpar}(#1)}
\newcommand{\Ugln}[1][\gln]{U_{\qpar}(#1)}
\newcommand{\Uglnone}[1][\gln]{U_{1}(#1)}
\newcommand{\Uglnstar}[1][\gln]{U_{\qpar}(#1)}

\newcommand{\vecrep}[1][\qpar]{\obstuff{V}_{#1}}
\newcommand{\ext}[2]{{\textstyle\bigwedge_{#1}^{#2}}}
\newcommand{\simple}[1][\lambda]{L(#1)}

\newcommand{\braid}[1][\qpar]{\hat{\morstuff{R}}_{#1}}

\newcommand{\levi}{\ell}
\newcommand{\cartan}{\mathfrak{h}}

\newcommand{\Ulv}[1][\levi]{U_{\vpar}(#1)}
\newcommand{\Ul}[1][\levi]{U_{\qpar}(#1)}

\newcommand{\Ulstar}[1][\levi]{U_{\qpar}(#1)}

\newcommand{\rep}[1][\levi]{\catstuff{Rep}_{\qpar}{#1}}
\newcommand{\repone}[1][\levi]{\catstuff{Rep}_{1}{#1}}
\newcommand{\repstar}[1][\levi]{\catstuff{Rep}_{\qpar}{#1}}

\newcommand{\fund}[1][\levi]{\catstuff{Fund}_{\qpar}{#1}}
\newcommand{\fundone}[1][\levi]{\catstuff{Fund}_{1}{#1}}
\newcommand{\fundstar}[1][\levi]{\catstuff{Fund}_{\qpar}{#1}}

\newcommand{\tilt}[1][\gln]{\catstuff{Tilt}_{\vpar}{#1}}

\newcommand{\lprod}{\circ_{h}^{\levi}}
\newcommand{\lbraid}[1][\qpar]{\hat{\morstuff{R}}_{#1}^{\levi}}

\newcommand{\mergemap}[2]{\raisebox{0.26cm}{\rotatebox{180}{$\morstuff{Y}$}}_{#1}^{#2}}
\newcommand{\splitmap}[2]{\morstuff{Y}_{#2}^{#1}}

\newcommand{\expl}[1][k]{\morstuff{Y}_{#1}^{1,\dots,1}}
\newcommand{\iexpl}[1][k]{\raisebox{0.26cm}{\rotatebox{180}{$\morstuff{Y}$}}^{#1}_{1,\dots,1}}

\newcommand{\cupl}[1][k]{\cup_{#1}^{\leftarrow}}
\newcommand{\cupr}[1][k]{\cup_{#1}^{\rightarrow}}
\newcommand{\capl}[1][k]{\cap_{#1}^{\leftarrow}}
\newcommand{\capr}[1][k]{\cap_{#1}^{\rightarrow}}

\newcommand{\wmapl}[1][{1,\obstuff{K}}]{\morstuff{r}_{#1}}
\newcommand{\wmapr}[1][{1,\obstuff{K}}]{\morstuff{r}^{#1}}

\newcommand{\lproj}[1][i]{\morstuff{pr}_{#1}^{\levi}}

\newcommand{\diacat}{\catstuff{Dia}}

\newcommand{\webv}[1][\gln]{\catstuff{Web}_{\vpar}{#1}}
\newcommand{\web}[1][\gln]{\catstuff{Web}_{\qpar}{#1}}

\newcommand{\webstar}[1][\gln]{\catstuff{Web}_{\qpar}{#1}}
\newcommand{\uob}[1][k]{\uparrow_{#1}}
\newcommand{\dob}[1][k]{\downarrow_{#1}}

\newcommand{\awebv}[1][\gln]{\catstuff{AWeb}_{\vpar}{#1}}
\newcommand{\aweb}[1][\gln]{\catstuff{AWeb}_{\qpar}{#1}}
\newcommand{\awebone}[1][\gln]{\catstuff{AWeb}_{1}{#1}}
\newcommand{\awebstar}[1][\gln]{\catstuff{AWeb}_{\qpar}{#1}}
\newcommand{\winding}[1][\obstuff{K}]{\rho_{#1}}
\newcommand{\iwinding}[1][\obstuff{K}]{\rho^{#1}}

\newcommand{\aff}[1][\webv]{\catstuff{Aff}({#1})}

\newcommand{\wprod}{\circ_{h}^{A}}

\newcommand{\ecirclel}[1][k]{\morstuff{c}^{\leftarrow}_{#1}}
\newcommand{\ecircler}[1][k]{\morstuff{c}^{\rightarrow}_{#1}}

\newcommand{\lweb}[1][\levi]{\catstuff{AWeb}_{\qpar}{#1}}
\newcommand{\lwebone}[1][\levi]{\catstuff{AWeb}_{1}{#1}}

\newcommand{\lideal}[1][\levi]{\mathcal{I}_{#1}}

\newcommand{\fulltwist}{\morstuff{ft}}

\newcommand{\wproj}[1][i]{\morstuff{pr}_{#1}^{w}}

\newcommand{\pfunctor}[1][\qpar]{\Gamma_{#1}}
\newcommand{\apfunctor}[1][\qpar]{\catstuff{A}\Gamma_{#1}}
\newcommand{\apfunctorstar}[1][\qpar]{\catstuff{A}\Gamma_{#1}}
\newcommand{\apfunctorl}[2]{\catstuff{A}\Gamma_{#1}^{#2}}

\newcommand{\ak}[1][{m,d}]{\mathcal{H}^{#1}_{\qpar}}
\newcommand{\akv}[1][{m,d}]{\mathcal{H}^{#1}_{\vpar}}
\newcommand{\akstar}[1][{m,d}]{\mathcal{H}^{#1}_{\qpar}}
\newcommand{\blob}[1][{m,d}]{\mathcal{B}^{#1}_{\qpar}}
\newcommand{\blobstar}[1][{m,d}]{\mathcal{B}^{#1}_{\qpar}}
\newcommand{\bloblevi}[1][\levi]{\mathcal{B}^{#1}_{\qpar}}
\newcommand{\bloblevistar}[1][\levi]{\mathcal{B}^{#1}_{\qpar}}

\newcommand{\blam}{\boldsymbol{\lambda}}
\newcommand{\lscalar}{s}
\newcommand{\parts}[1][{m,d}]{P(#1)}
\newcommand{\std}[1][\blam]{Std(#1)}

\newcommand{\akgens}[1][i]{\morstuff{T}_{#1}}
\newcommand{\jm}[1][i]{\morstuff{X}_{#1}}
\newcommand{\relement}[1][i]{\morstuff{R}_{#1}}

\newcommand{\idealak}[1][>1]{\mathcal{I}_{#1}}
\newcommand{\akideal}[1][>1]{\mathcal{J}_{#1}}

\newcommand{\akmap}[1][{m,d}]{\pi_{#1}}
\newcommand{\blmap}[1][{m,d}]{\overline{\pi}_{#1}}
\newcommand{\blmaplevi}[1][\levi]{\overline{\pi}_{#1}}

\newcommand{\basis}[1][{s,t,\sigma}]{\nu_{#1}^{\qpar}}
\newcommand{\basisset}[1][{\lweb[\levi]}]{B(#1)}
\newcommand{\basisone}[1][{s,t,\sigma}]{\nu_{#1}^{1}}
\newcommand{\basissetone}[1][{\lwebone[\levi]}]{B(#1)}

\newcommand{\lbasis}[1][{s,t,\sigma}]{\morstuff{n}_{#1}^{\qpar}}
\newcommand{\lbasisset}[1][{\fund[\levi]}]{B(#1)}
\newcommand{\lbasisone}[1][{s,t,\sigma}]{\morstuff{n}_{#1}^{1}}
\newcommand{\lbasissetone}[1][{\fundone[\levi]}]{B(#1)}

\usepackage{aliascnt}
\def\NewTheorem#1{%
\newaliascnt{#1}{equation}%
\newtheorem{#1}[#1]{#1}%
\aliascntresetthe{#1}%
\expandafter\def\csname #1autorefname\endcsname{#1}%
}
\def\equationautorefname~#1\null{(#1)\null}

\numberwithin{equation}{subsection}

\NewTheorem{Proposition}
\NewTheorem{Theorem}
\NewTheorem{Corollary}
\AtEndEnvironment{Corollary}{\null\hfill$\square$}%
\NewTheorem{Lemma}
\theoremstyle{definition}
\NewTheorem{Definition}
\NewTheorem{Notation}
\NewTheorem{Example}
\AtEndEnvironment{Example}{\null\hfill$\Diamond$}%
\theoremstyle{remark}
\NewTheorem{Remark}
\NewTheorem{Assumption}
\NewTheorem{Hypothesis}
\NewTheorem{Conjecture}
\NewTheorem{Question}
\NewTheorem{Observation}
\NewTheorem{Fact}
\NewTheorem{Conclusion}
\NewTheorem{lemma}

\usepackage[hypertexnames=false]{hyperref}
\usepackage{bookmark}
\hypersetup{
pdftoolbar=true,
pdfmenubar=true,
pdffitwindow=false,
pdfstartview={FitH},
pdftitle={Annular webs and Levi subalgebras},
pdfauthor={Abel Lacabanne, Daniel Tubbenhauer and Pedro Vaz},
pdfsubject={},
pdfcreator={Abel Lacabanne, Daniel Tubbenhauer and Pedro Vaz},
pdfproducer={Abel Lacabanne, Daniel Tubbenhauer and Pedro Vaz},
pdfkeywords={},
pdfnewwindow=true,
colorlinks=true,
linkcolor=mydarkblue,
citecolor=teal,
filecolor=magenta,
urlcolor=orchid,
linkbordercolor=lava,
citebordercolor=teal,
urlbordercolor=orchid,
linktocpage=true
}

\newcommand{\nnfootnote}[1]{%
\begin{NoHyper}
\renewcommand\thefootnote{}\footnote{#1}%
\addtocounter{footnote}{-1}%
\end{NoHyper}
}

\def\makeautorefname#1#2{\csdef{#1autorefname}{#2}}

\makeautorefname{section}{Section}%
\makeautorefname{subsection}{Section}%
\makeautorefname{subsubsection}{Section}%

\begin{document}
\title[Annular webs and Levi subalgebras]
{Annular webs and Levi subalgebras}
\author[A. Lacabanne, D. Tubbenhauer and P. Vaz]{Abel Lacabanne, Daniel Tubbenhauer and Pedro Vaz}

\address{A.L.: Laboratoire de Math{\'e}matiques Blaise Pascal (UMR 6620), Universit{\'e} Clermont Auvergne, Campus Universitaire des C{\'e}zeaux, 3 place Vasarely, 63178 Aubi{\`e}re Cedex, France,\newline \href{http://www.normalesup.org/~lacabanne}{www.normalesup.org/$\sim$lacabanne},
	\href{https://orcid.org/0000-0001-8691-3270}{ORCID 0000-0001-8691-3270}}
\email{abel.lacabanne@uca.fr}

\address{D.T.: The University of Sydney, School of Mathematics and Statistics F07, Office Carslaw 827, NSW 2006, Australia, \href{http://www.dtubbenhauer.com}{www.dtubbenhauer.com}, \href{https://orcid.org/0000-0001-7265-5047}{ORCID 0000-0001-7265-5047}}
\email{daniel.tubbenhauer@sydney.edu.au}

\address{P.V.: Institut de Recherche en Math{\'e}matique et Physique, 
	Universit{\'e} Catholique de Louvain, Chemin du Cyclotron 2,  
	1348 Louvain-la-Neuve, Belgium, \href{https://perso.uclouvain.be/pedro.vaz}{https://perso.uclouvain.be/pedro.vaz}, \href{https://orcid.org/0000-0001-9422-4707}{ORCID 0000-0001-9422-4707}}
\email{pedro.vaz@uclouvain.be}

\begin{abstract}
For any Levi subalgebra of the form $\mathfrak{l}=\mathfrak{gl}_{l_{1}}\oplus\dots\oplus\mathfrak{gl}_{l_{d}}\subseteq\mathfrak{gl}_{n}$ we construct a quotient of 
the category of annular quantum $\mathfrak{gl}_{n}$ webs that is 
equivalent to the category
of finite dimensional representations of quantum $\mathfrak{l}$ 
generated by exterior powers of the vector representation. This can be interpreted as an annular version of skew Howe duality, gives a description of the representation category of $\mathfrak{l}$ by additive idempotent completion, and a web version of the generalized blob algebra.
\end{abstract}

\nnfootnote{\textit{Mathematics Subject Classification 2020.} Primary: 17B37; Secondary: 18M05, 18M15, 20C08.}
\nnfootnote{\textit{Keywords.} Representations of quantum general linear groups and Levi subalgebras, annular webs, quantum skew Howe duality, diagrammatic presentation, generalized blob algebras.}

\maketitle

\tableofcontents

\arrayrulewidth=0.5mm
\setlength{\arrayrulewidth}{0.5mm}


\section{Introduction}\label{S:Introduction}


Throughout fix $n,l_{1},\dots,l_{d}\in\N$ with $\sum_{i=1}^{d}l_{i}=n$.


\subsection{Webs, and Schur--Weyl and Howe duality}\label{SS:IntroductionPartOne}


The so-called \emph{Schur--Weyl duality} has 
played a key role ever since the early days of representation theory. 
It relates representations of the symmetric group $\sym$
and the general linear group 
$\GLN=\GLN(\C)$, and has been 
generalized in many ways. The representation used to relate 
these two groups is $(\C^{n})^{\otimes m}$.

Let $\Levi=\GLN[l_{1}]\times\dots
\times\GLN[l_{d}]\subseteq\GLN$, 
and let us for simplicity stay over $\C$ 
for now. Two generalizations of Schur--Weyl duality 
are of crucial importance for this paper.
Firstly, the \emph{Schur--Weyl duality of $(\Z/d\Z)\wr\sym$}
(that is, type $G(d,1,m)$) 
and $\Levi$ from \cite{ArTeYa-schur-weyl-ariki-koike}, \cite{SaSh-schur-weyl-ariki-koike} and 
\cite{Hu-schur-weyl-ak-algebra} (see also \cite{MaSt-complex-reflection-groups} for a nice and self-contained discussion of this duality). 
Here the underlying representation is again $(\C^{n})^{\otimes m}$.
Secondly, \emph{skew (type A) Howe duality}, see 
\cite{Ho-remarks-invariant-theory} and \cite{Ho-perspectives-invariant-theory}, 
relating $\GLN[N]$ and $\GLN$ via their action on 
the exterior algebra $\EXT(\C^{N}\otimes\C^{n})$.

As explained in \cite{CaKaMo-webs-skew-howe}, a diagrammatic 
interpretation of skew Howe duality is given by 
\emph{(exterior $\GLN$) webs}. (The same diagrammatics goes 
under many names, including birdtracks \cite{Cv-bridtracks}
or spiders \cite{Ku-spiders-rank-2}.) In some sense, in \cite{CaKaMo-webs-skew-howe} skew Howe duality relating 
$\GLN[{N\in\N}]$ and 
$\GLN$ takes the form of an equivalence between 
the category of webs and the category of $\GLN$-representations generated by 
$\set[\big]{\EXT[k]\C^{n}|k\in\set{1,\dots,n}}$, with the web category 
being obtained by using all $\GLN[N]$ for $N\in\N$. 
After additive idempotent completion 
webs even give a diagrammatic interpretation of all finite dimensional $\GLN$-representations.

In this paper we show that an explicit quotient of the category of 
\emph{annular (exterior $\GLN$) webs} is equivalent 
to the category of $\Levi$-representations generated by the set 
$\set[\big]{\EXT[k]\C^{n}|k\in\set{1,\dots,n}}$. As before, additive idempotent completion gives a description of all finite dimensional $\Levi$-representations. This, in some sense, is a form 
of what could be called \emph{annular skew Howe duality} 
(we avoid the notion affine as its meaning is context depending)
or \emph{skew type $G(d,1,m\in\N)$ Howe duality}. 


\subsection{The main result and relations to other works}\label{SS:IntroductionPartTwo}


We now give a few details and change to the universal enveloping algebras.
We consider a \emph{Levi subalgebra} of the form $\levi=\gln[l_{1}]\oplus\dots\oplus\gln[l_{d}]\subseteq\gln$ (in this paper we write $\levi$ 
instead of the usual notation $\mathfrak{l}$ for readability).
Let $\K$ denote a field containing an element 
$\qpar\in\K$ that is not a root of unity and 
additional variables $\U=\set{\upar[1],\dots,\upar[d]}$ and their inverses, and 
let further $\Kone$ denote a field of characteristic zero 
containing $\U$ and their inverses. 
With these 
ground fields the category of finite dimensional 
$\Ugln[\levi]$-representations over $\K$ respectively of finite dimensional 
$\Uglnone[\levi]$-representations 
over $\Kone$ are semisimple. (We should warn the reader: as explained in the main body of the text there are some nontrivial quantization issues 
and we carefully need to distinguish the two cases over $\K$ and $\Kone$.)

In \autoref{S:AnnularWebs} we define a $\K$-linear 
category of annular webs 
$\aweb$ as well as a quotient $\aweb[\levi]$ by 
evaluating essential circles using the variables 
$\U$ and their inverses. Similarly over $\Kone$, where we 
write $\awebone$ and $\awebone[\levi]$. Let $\fund[\levi]$ 
respectively $\fundone[\levi]$ denote the categories 
of $\Ugln[\levi]$- and $\Uglnone[\levi]$-representations generated by 
the exterior powers of the vector representation.
Our main result is \autoref{T:EquivalenceMain} showing that 
$\aweb[\levi]$ is equivalent to $\fund[\levi]$ and that 
$\awebone[\levi]$ is pivotally equivalent to $\fundone[\levi]$.
The main ingredients in the proof of \autoref{T:EquivalenceMain} 
are the usual diagrammatic ideas, 
the Schur--Weyl type dualities from \cite{SaSh-schur-weyl-ariki-koike} 
as well as the \emph{explosion trick}, which utilizes the semisimplicity.

An almost direct consequence of \autoref{T:EquivalenceMain} is that 
the endomorphism algebras of annular webs corresponding to tensor products 
of the vector representation can be described explicitly. 
As we will see in \autoref{SS:EquivalenceAK} these are given 
by certain \emph{row quotients of Ariki--Koike algebras} 
(Ariki--Koike algebras are Hecke algebras of $(\Z/d\Z)\wr\sym$, see 
for example \cite{ArKo-hecke-algebra}, \cite{BrMa-hecke}, 
\cite{Ch-gelfandtzetlin}) as studied in \cite{LaVa-schur-weyl-ariki-koike}.
In the special case of two row quotients, which corresponds to $\levi=\gln[1]\oplus\dots\oplus\gln[1]$ 
being the \emph{Cartan subalgebra}, these algebras are Martin--Woodcock's 
\emph{generalized blob algebras} \cite{MaWo-generalized-blob}. 
We thus obtain a web description of generalized blob algebras, see 
\autoref{SS:EquivalenceBlob} for details.
This web description allows us to answer 
two conjectures of Cautis--Kamnitzer 
\cite[Conjectures 10.2 and 10.3]{CaKa-q-satake-sln} affirmatively, up to technicalities as detailed in \autoref{R:EquivalenceGenBlobThird}.

Let us also mention that our work is inspired by \cite{Qu-annular-skew-howe}
(which gives another, very honest, version of annular skew Howe duality), the aforementioned paper \cite{CaKa-q-satake-sln} as well as \cite{QuWe-extremal-projectors-2} 
(which proves \autoref{T:EquivalenceMain} for $\Kone=\C$ and the special case of the Cartan subalgebra). With respect to \cite{Qu-annular-skew-howe} 
and \cite{QuWe-extremal-projectors-2}, which are partially motivated 
by skein theory, we should also warn the reader that 
the monoidal structure on $\aweb[\levi]$ coming from skein theory and the one on $\fund[\levi]$ coming from the Hopf algebra 
structure of $\Ugln[\levi]$ do not seem 
to be compatible, see \autoref{SS:EquivalenceMonoidal} for a discussion.
\medskip

\noindent\textbf{Acknowledgments.}
We would like to thank Mikhail Khovanov and Alvaro Martinez Ruiz for helpful discussions, and the referee for their meticulous and careful reading (multiple times!) of the manuscript, and helpful suggestions.

AL was partially supported by a PEPS JCJC grant from INSMI (CNRS).
DT was supported, in part, by the Australian Research Council, and by 
their dedication akin to insanity.
PV was supported by the Fonds de la Recherche Scientifique - FNRS under Grant no. MIS-F.4536.19.
The last author visited the first author during this project, which was in part sponsored by the Universit{\'e} Clermont Auvergne. The hospitality of the Universit{\'e} Clermont Auvergne is gratefully acknowledged.


\section{Notations and conventions}\label{S:Prelim}


We start by specifying our notations.

\begin{Notation}\label{N:PrelimFixN}
Recall that we fixed $n,l_{1},\dots,l_{d}\in\N$ with $\sum_{i=1}^{d}l_{i}=n$.
These will be used via the general linear Lie algebra
$\gln$ and a Levi subalgebra $\levi$ given by $\levi=\gln[l_{1}]\oplus\dots\oplus\gln[l_{d}]\subseteq\gln$.
For $l_{1}=\dots=l_{d}=1$, so that 
$\levi=\gln[1]\oplus\dots\oplus\gln[1]$ is the Cartan 
subalgebra, we will write $\cartan$ instead of $\levi$.
\end{Notation}

\begin{Notation}\label{N:PrelimField}
We now specify our underlying ground rings.	
\begin{enumerate}

\item For essential circles, see \autoref{SS:AnnularWebsCircles} below,
we need extra polynomial elements (these can be ignored otherwise).
We denote by $\U=\set{\upar[1],\dots,\upar[d]}$ variables which will play this role.

\item Let $\Zv=\Z[\vpar,\vpar^{-1},\U,\U^{-1}]$ for some indeterminate $\vpar$.

\item We let $\K$ denote a field 
containing $\U$ and $\U^{-1}$ and an element $\qpar$ which is not a root of unity. Let also $\Kone$ be a field of characteristic zero containing $\U$ and $\U^{-1}$.
(Note that $\mathrm{char}(\K)$ is allowed to be a prime, but we assume that
$\mathrm{char}(\Kone)=0$.)

\item We also see $\K$ as the specialization and scalar extension $\placeholder\otimes_{\Zv}\K$ given by $\vpar\mapsto\qpar$, and $\Kone$ as the specialization and scalar extension $\placeholder\otimes_{\Zv}\Kone$ given by $\vpar\mapsto 1$. We will apply scalar extension to $\Zv$-linear categories, 
and since this will play an important role we will indicate the specialization 
accordingly. The two important specializations for $\K$ and $\Kone$ 
are distinguished by using $\qpar$ respectively $1$ as a subscript.

\item If not specified otherwise, then $\otimes$ denotes the tensor product over the ground ring, which is either $\Zv$, $\K$ or $\Kone$.

\end{enumerate}
\end{Notation}

Not everything in this paper is defined over $\Zv$, and 
statements for $\K$ are not strictly related to the ones for $\Kone$. 
If we use $\Zv$, then we can specialize without any problem. But if we do not work over $\Zv$, then, to not double parts of the text, we use the following crucial simplification:

\begin{Notation}\label{N:PrelimFields}
Throughout (on the diagrammatic and the representation theoretical sides)
we always specify our conventions for $\K$ and leave the 
analog conventions and lemmas for $\Kone$ implicit; the ones for 
$\Kone$ are always the $\qpar=1$ versions of the ones for $\K$.
When lemmas {\etc} are the same for $\K$ and $\Kone$, then 
we will use simply $\qpar$ ({\eg} we write $\K$ and not $\K$ and 
$\Kone$) to indicate this. 
We will stress whenever the statements for $\K$ and $\Kone$ are significantly different.
\end{Notation}

We will use quantum numbers, 
factorials and binomials
viewed as elements of $\Zv$.
That is, for $a\in\Z$ and $b\in\N$ we let $\qnum{0}=0$, $\qfac{0}=1=\qbin{a}{0}$, $\qnum{a}=-\qnum{-a}$ for $a<0$ and 
otherwise
\begin{gather*}
\qnum{a}=\vpar^{a-1}+\vpar^{a-3}+\dots+\vpar^{-a+3}+\vpar^{-a+1}
,\quad
\qfac{b}=\qnum{b}\qnum{b-1}\dots\qnum{1}
,\quad
\qbinn{a}{b}=
\frac{\qnum{a}\qnum{a-1}\dots\qnum{a-b+1}}{\qfac{b}}.
\end{gather*}
The following will be used silently throughout.

\begin{Lemma}\label{L:PrelimInvertible}
All quantum binomials are invertible in $\Kstar$.
\end{Lemma}

\begin{proof}
Easy and omitted.
\end{proof}

\begin{Notation}\label{N:PrelimNotationMonoidal}
We work with strict pivotal (thus, monoidal) categories, were we strictify 
categories if necessary (by the usual strictification theorems this 
restriction is for convenience only). We have two directions of composition, vertical $\vcirc$ and horizontal $\hcirc$, as well as a duality ${}^{\pivo}$ operation. The monoidal unit is 
denoted by $\munit$, and identity morphisms are denoted by $\idmor$. 
We will also distinguish objects and morphisms using different fonts, {\eg} $\obstuff{K}$ and $\morstuff{f}$.

As we will see, {\eg} in
\autoref{SS:EquivalenceMonoidal} below, it will turn out 
to be important to carefully distinguish the monoidal product on 
the various categories we consider. This will only play a role 
on the level of morphisms, and we use various symbols 
for monoidal products between morphisms if necessary.
\end{Notation}

\begin{Notation}\label{N:PrelimDiagramNotation}
We now summarize the diagrammatic conventions that we will use 
in this paper. 
\begin{enumerate}

\item The following illustration of the interchange law summarizes 
our reading conventions:
\begin{gather*}
(\idmor\hcirc\morstuff{g})\vcirc(\morstuff{f}\hcirc\idmor) 
=
\begin{tikzpicture}[anchorbase,scale=0.22,tinynodes]
\draw[spinach!35,fill=spinach!35] (-5.5,2) rectangle (-2,0.5);
\draw[spinach!35,fill=spinach!35] (-0.5,3) rectangle (3,4.5);
\draw[thick,densely dotted] (-5.5,2.5) node[left,yshift=-2pt]{$\vcirc$} 
to (3.5,2.5) node[right,yshift=-2pt]{$\vcirc$};
\draw[thick,densely dotted] (-1.25,5) 
node[above,yshift=-2pt]{$\hcirc$} to (-1.25,0) node[below]{$\hcirc$};
\draw[usual] (-5,2) to (-5,3.75) node[right]{$\dots$} to (-5,5);
\draw[usual] (-2.5,2) to (-2.5,5);
\draw[usual] (0,0) to (0,1.25) node[right]{$\dots$} to (0,3);
\draw[usual] (2.5,0) to (2.5,3);
\draw[usual] (-5,0) to (-5,0.25) node[right]{$\dots$} to (-5,0.5);
\draw[usual] (-2.5,0) to (-2.5,0.5);
\draw[usual] (2.5,4.5) to (2.5,5);
\draw[usual] (0,4.5) to (0,4.75) node[right]{$\dots$} to (0,5);
\node at (-3.75,1.0) {$\morstuff{f}$};
\node at (1.25,3.5) {$\morstuff{g}$};
\end{tikzpicture}
=
\begin{tikzpicture}[anchorbase,scale=0.22,tinynodes]
\draw[spinach!35,fill=spinach!35] (-5.5,1.75) rectangle (-2,3.25);
\draw[spinach!35,fill=spinach!35] (-0.5,1.75) rectangle (3,3.25);
\draw[usual] (-5,3.25) to (-5,4.125) node[right]{$\dots$} to (-5,5);
\draw[usual] (-2.5,3.25) to (-2.5,5);
\draw[usual] (0,0) to (0,0.625) node[right]{$\dots$} to (0,1.75);
\draw[usual] (2.5,0) to (2.5,1.75);
\draw[usual] (-5,0) to (-5,0.625) node[right]{$\dots$} to (-5,1.75);
\draw[usual] (-2.5,0) to (-2.5,1.75);
\draw[usual] (2.5,3.25) to (2.5,5);
\draw[usual] (0,3.25) to (0,4.125) node[right]{$\dots$} to (0,5);
\node at (-3.75,2.25) {$\morstuff{f}$};
\node at (1.25,2.25) {$\morstuff{g}$};
\end{tikzpicture}
=
\begin{tikzpicture}[anchorbase,scale=0.22,tinynodes]
\draw[spinach!35,fill=spinach!35] (5.5,2) rectangle (2,0.5);
\draw[spinach!35,fill=spinach!35] (0.5,3) rectangle (-3,4.5);
\draw[thick, densely dotted] (5.5,2.5) 
node[right,yshift=-2pt]{$\vcirc$} to (-3.5,2.5) node[left,yshift=-2pt]{$\vcirc$};
\draw[thick, densely dotted] (1.25,5) 
node[above,yshift=-2pt]{$\hcirc$} to (1.25,0) node[below]{$\hcirc$};
\draw[usual] (2.5,2) to (2.5,3.75) node[right]{$\dots$} to (2.5,5);
\draw[usual] (5,2) to (5,5);
\draw[usual] (0,0) to (0,3);
\draw[usual] (-2.5,0) to (-2.5,1.25) node[right]{$\dots$} to (-2.5,3);
\draw[usual] (2.5,0) to (2.5,0.25) node[right]{$\dots$} to (2.5,0.5);
\draw[usual] (5,0) to (5,0.5);
\draw[usual] (-2.5,4.5) to (-2.5,4.75) node[right]{$\dots$} to (-2.5,5);
\draw[usual] (0,4.5) to (0,5);
\node at (3.75,1.0) {$\morstuff{g}$};
\node at (-1.25,3.5) {$\morstuff{f}$};
\end{tikzpicture}
=(\morstuff{f}\hcirc\idmor)\vcirc(\idmor\hcirc\morstuff{g}).
\end{gather*}
That is, we read diagrams from bottom to top and left to right.

\item As we will recall below, webs are certain types of labeled (with numbers $a\in\N$) and oriented graphs. Some labels and 
orientations are determined by others, and we will often omit orientations and 
labels that can be recovered from the given data to avoid clutter.

\item If labels or orientations are omitted altogether, then 
the displayed webs are a shorthand for any web of the same shape and 
legit labels and orientations.

\item We use webs with edges labeled by $a\in\Z$, where 
we use the convention that edges of label $0$ are omitted from the illustrations, and edges with label not in $\N$ set the web to zero.
(We will use negative labels, but for objects and not for edges in webs.)

\item We also use strands labeled by objects, {\ie} for $\obstuff{K}=(k_{1},\dots,k_{m})$
\begin{gather*}
\begin{tikzpicture}[anchorbase,scale=1]
\draw[usual] (0,0)node[below]{$\obstuff{K}$} to (0,0.5)node[above]{$\obstuff{K}$};
\end{tikzpicture}
=
\begin{tikzpicture}[anchorbase,scale=1]
\draw[usual] (0,0)node[below]{$k_{1}$} to (0,0.5)node[above]{$k_{1}$};
\draw[usual] (1,0)node[below]{$k_{m}$} to (1,0.5)node[above]{$k_{m}$};
\node at (0.5,0.25) {$\dots$};
\end{tikzpicture}
\end{gather*}
to indicate an arbitrary (but finite) number of parallel strings.

\end{enumerate}
\end{Notation}

\section{Web categories in the plane}\label{S:Webs}

This section serves as a reminder on web categories 
and their basic properties.
Details can be found in many books and papers, {\eg} \cite{TuVi-monoidal-tqft} 
for general diagrammatics, and \cite{CaKaMo-webs-skew-howe} 
or, closer to our conventions,
\cite{LaTu-gln-webs} and the references therein for web categories. 
The proofs of the statements
below are easy or can be found in {\loccit}


\subsection{Preliminaries}\label{SS:WebsGeneral}


Let $R$ be a commutative ring with unit.

\begin{Definition}\label{D:WebsDia}
For the purpose of this section, a \emph{diagram category} 
$\diacat$ is a pivotal 
$R$-linear category with objects 
$\hcirc$-generated by $\uob$ with $\dob=(\uob)^{\pivo}$ for $k\in\N$ 
with $\munit$ being the empty word, and a braid group action on upwards objects, meaning morphisms $\braid[k,l]\colon\uob\hcirc\uob[l]\to\uob[l]\hcirc\uob$ for each simple braid group generator that 
satisfy the braid relations. 
\end{Definition}

We also write
$\obstuff{K}=(k_{1},\dots,k_{m})\in\Z^{m}$
for $m\in\N$ for the objects of $\diacat$, where we use the notations
$k\leftrightsquigarrow\uob$ and 
$-k\leftrightsquigarrow\dob$ for $k>0$.
We illustrate the (co)evaluation morphisms as
\begin{gather*}
\scalebox{0.98}{$\begin{tikzpicture}[anchorbase,scale=1]
\draw[usual,directed=0.99] (0,0)node[below]{$k$} to[out=90,in=180] (0.5,0.5)node[above]{\phantom{k}} to[out=0,in=90] (1,0)node[below]{${-}k$};
\end{tikzpicture}
\hspace*{-0.1cm}\colon
\uob\hcirc\dob\to\munit
,
\begin{tikzpicture}[anchorbase,scale=1]
\draw[usual,directed=0.99] (1,0)node[below]{$k$} to[out=90,in=0] (0.5,0.5)node[above]{\phantom{k}} to[out=180,in=90] (0,0)node[below]{${-}k$};
\end{tikzpicture}
\hspace*{-0.1cm}\colon
\dob\hcirc\uob\to\munit
,
\begin{tikzpicture}[anchorbase,scale=1]
\draw[usual,directed=0.99] (0,0)node[above]{${-}k$} to[out=270,in=180] (0.5,-0.5)node[below]{\phantom{k}} to[out=0,in=270] (1,0)node[above]{$k$};
\end{tikzpicture}
\hspace*{-0.1cm}\colon
\munit\to\dob\hcirc\uob
,
\begin{tikzpicture}[anchorbase,scale=1]
\draw[usual,directed=0.99] (1,0)node[above]{${-}k$} to[out=270,in=0] (0.5,-0.5)node[below]{\phantom{k}} to[out=180,in=270] (0,0)node[above]{$k$};
\end{tikzpicture}
\hspace*{-0.1cm}\colon
\munit\to\uob\hcirc\dob,$}
\end{gather*}
and the braid group action as \emph{$(k,l)$-crossings} (overcrossings and undercrossings):
\begin{gather*}
\text{over}:
\begin{tikzpicture}[anchorbase,scale=1]
\draw[usual,directed=0.99] (1,0)node[below]{$l$} to (0,1)node[above]{$l$};
\draw[usual,directed=0.99,crossline] (0,0)node[below]{$k$} to (1,1)node[above]{$k$};
\end{tikzpicture}
\colon\uob\hcirc\uob[l]\to\uob[l]\hcirc\uob
,\quad
\text{under}:
\begin{tikzpicture}[anchorbase,scale=1]
\draw[usual,directed=0.99] (0,0)node[below]{$k$} to (1,1)node[above]{$k$};
\draw[usual,directed=0.99,crossline] (1,0)node[below]{$l$} to (0,1)node[above]{$l$};
\end{tikzpicture}
\colon\uob\hcirc\uob[l]\to\uob[l]\hcirc\uob.
\end{gather*}

\begin{Definition}\label{D:WebsDiaAnti}
The \emph{diagrammatic antiinvolution} $(\placeholder)^{\danti}$ on 
a diagram category
is defined on objects by $(\obstuff{K}\hcirc\obstuff{L}\hcirc\dots)^{\danti}=
\obstuff{K}^{\pivo}\hcirc\obstuff{L}^{\pivo}\hcirc\dots$
and on morphisms as in \autoref{Eq:PrelimDiaInvo} below. 
The \emph{diagrammatic involution} $(\placeholder)^{\dinvo}$ on 
a monoidal diagram category is defined on objects by $(\obstuff{K}\hcirc\obstuff{L}\hcirc\dots)^{\dinvo}=\dots\hcirc\obstuff{L}\hcirc\obstuff{K}$
and on morphisms as in \autoref{Eq:PrelimDiaInvo} below.

\begin{gather}\label{Eq:PrelimDiaInvo}
\left(\begin{tikzpicture}[anchorbase,scale=0.22]
\draw[spinach!35,fill=spinach!35] (-5.5,2) rectangle (-2,0.5);
\draw[usual] (-5,2) to (-5,2.5);
\draw[usual] (-2.5,2) to (-2.5,2.5);
\draw[usual] (-5,0) to (-5,0.5);
\draw[usual] (-2.5,0) to (-2.5,0.5);
\node at (-3.75,1.3) {$\morstuff{f}$};
\end{tikzpicture}\right)^{\danti}
=
\begin{tikzpicture}[anchorbase,scale=0.22]
\draw[spinach!35,fill=spinach!35] (-5.5,2) rectangle (-2,0.5);
\draw[usual] (-5,2) to (-5,2.5);
\draw[usual] (-2.5,2) to (-2.5,2.5);
\draw[usual] (-5,0) to (-5,0.5);
\draw[usual] (-2.5,0) to (-2.5,0.5);
\node at (-3.75,1.35) {\reflectbox{\rotatebox{180}{$\morstuff{f}$}}};
\end{tikzpicture}
,\quad
\left(\begin{tikzpicture}[anchorbase,scale=0.22]
\draw[spinach!35,fill=spinach!35] (-5.5,2) rectangle (-2,0.5);
\draw[usual] (-5,2) to (-5,2.5);
\draw[usual] (-2.5,2) to (-2.5,2.5);
\draw[usual] (-5,0) to (-5,0.5);
\draw[usual] (-2.5,0) to (-2.5,0.5);
\node at (-3.75,1.3) {$\morstuff{f}$};
\end{tikzpicture}\right)^{\dinvo}
=
\begin{tikzpicture}[anchorbase,scale=0.22]
\draw[spinach!35,fill=spinach!35] (-5.5,2) rectangle (-2,0.5);
\draw[usual] (-5,2) to (-5,2.5);
\draw[usual] (-2.5,2) to (-2.5,2.5);
\draw[usual] (-5,0) to (-5,0.5);
\draw[usual] (-2.5,0) to (-2.5,0.5);
\node at (-3.75,1.3) {\reflectbox{$\morstuff{f}$}};
\end{tikzpicture}
.
\end{gather}
\end{Definition}

\begin{Lemma}\label{L:WebsDiaAnti}
The diagrammatic antiinvolution is an antiinvolution, 
and the diagrammatic involution is an involution.
\end{Lemma}

\begin{proof}
Easy and omitted.
\end{proof}

\begin{Definition}\label{D:WebsMate}
\emph{Mating} in $\diacat$ is the process of applying (co)evaluation morphisms. For example,
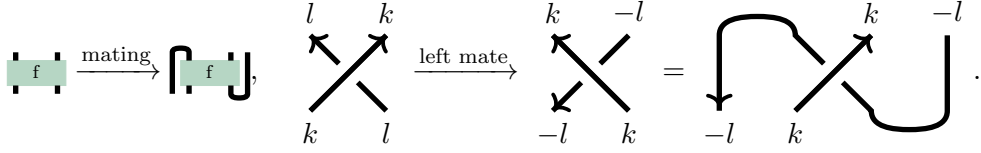
\begin{gather*}
\begin{tikzpicture}[anchorbase,scale=0.22,tinynodes]
\draw[spinach!35,fill=spinach!35] (-5.5,2) rectangle (-2,0.5);
\draw[usual] (-5,2) to (-5,2.5);
\draw[usual] (-2.5,2) to (-2.5,2.5);
\draw[usual] (-5,0) to (-5,0.5);
\draw[usual] (-2.5,0) to (-2.5,0.5);
\node at (-3.75,1.0) {$\morstuff{f}$};
\end{tikzpicture}
\xrightarrow{\text{mating}}
\begin{tikzpicture}[anchorbase,scale=0.22,tinynodes]
\draw[spinach!35,fill=spinach!35] (-5.5,2) rectangle (-2,0.5);
\draw[usual] (-5,2) to (-5,2.5) to[out=90,in=0] (-5.5,2.75) to[out=180,in=90] (-6,2.5) to (-6,0);
\draw[usual] (-2.5,2) to (-2.5,2.5);
\draw[usual] (-5,0) to (-5,0.5);
\draw[usual] (-1.5,2.5) to (-1.5,0) to[out=270,in=0] (-2,-0.25) to[out=180,in=270] (-2.5,0) to (-2.5,0.5);
\node at (-3.75,1.0) {$\morstuff{f}$};
\end{tikzpicture}
,\quad
\begin{tikzpicture}[anchorbase,scale=1]
\draw[usual,directed=0.99] (1,0)node[below]{$l$} to (0,1)node[above]{$l$};
\draw[usual,directed=0.99,crossline] (0,0)node[below]{$k$} to (1,1)node[above]{$k$};
\end{tikzpicture}
\xrightarrow{\text{left mate}}
\begin{tikzpicture}[anchorbase,scale=1]
\draw[usual,directed=0.99] (0,1)node[above]{${-}l$} to (-1,0)node[below]{${-}l$};
\draw[usual,directed=0.99,crossline] (0,0)node[below]{$k$} to (-1,1)node[above]{$k$};
\end{tikzpicture}
=
\begin{tikzpicture}[anchorbase,scale=1]
\draw[usual,directed=0.99] (2,1)node[above]{${-}l$} to (2,0) to[out=270,in=0] (1.5,-0.25) to[out=180,in=270] (1,0) to (0,1) to[out=90,in=0] (-0.5,1.25) to[out=180,in=90] (-1,1) to (-1,0)node[below]{${-}l$};
\draw[usual,directed=0.99,crossline] (0,0)node[below]{$k$} to (1,1)node[above]{$k$};
\end{tikzpicture}
.
\end{gather*}
To be precise, a mate of a morphism 
$\morstuff{f}$ is any morphism obtained from $\morstuff{f}$ 
by applying (co)evaluation morphisms. A left and right mate of $\morstuff{f}$ is a mate that only uses left or right (co)evaluation morphisms. Finally, \emph{the} left and right mate of a crossing is as indicated above.
\end{Definition}

Note that mating produces many morphisms from a given set of morphisms.

\begin{Lemma}\label{L:WebsTangle}
Suppose that $\End_{\diacat}(\munit)\cong\End_{\diacat}(\uob)\cong R$.
If the left mate of the $(k,l)$-overcrossing is invertible,
then all mates of the $(k,l)$-overcrossing and its inverse span 
a pivotal subcategory equivalent to $R$-linear (labeled) tangles. 
In general, if the left mate of the $(k,l)$-overcrossing is invertible,
then all mates of the $(k,l)$-overcrossing and its inverse span 
a pivotal subcategory equivalent to $R$-linear (labeled) framed tangles.
\end{Lemma}

\begin{proof}
Note that the assumption $\End_{\diacat}(\uob)\cong R$ 
implies the Reidemeister I relation up to scalars, which is 
enough to copy the argument in \cite[Theorem X11.2.2]{Ka-quantum-groups}.
The second claim follows similarly.
\end{proof}

We say $\Hom_{\diacat}(\obstuff{K},\obstuff{L})$ is \emph{determined by}
$\Hom_{\diacat}(\obstuff{K}^{\prime},\obstuff{L}^{\prime})$ if there exists
an isomorphism of $R$-modules between them.

\begin{Lemma}\label{L:WebsUpwards}
If the left mates of the $(k,l)$-overcrossings are invertible,
then all hom-spaces in $\diacat$ are determined by upwards hom-spaces.
\end{Lemma}

\begin{proof}
By \autoref{L:WebsTangle}, the assumptions imply that 
\begin{gather*}
\begin{tikzpicture}[anchorbase,scale=0.22,tinynodes]
\draw[spinach!35,fill=spinach!35] (-5.5,2) rectangle (-2,0.5);
\draw[usual,directed=0.99] (-5,3) to (-5,2);
\draw[usual,directed=0.99] (-2.5,2) to (-2.5,3);
\draw[usual,directed=0.99] (-5,0.5) to (-5,-0.5);
\draw[usual,directed=0.99] (-2.5,0.5) to (-2.5,-0.5);
\node at (-3.75,1.0) {$\morstuff{f}$};
\end{tikzpicture}
\mapsto
\begin{tikzpicture}[anchorbase,scale=0.22,tinynodes]
\draw[usual,directed=0.99] (-0.5,-1.5) to (-0.5,3) to[out=90,in=0] (-2.75,3.75) to[out=180,in=90] (-5,3) to (-5,2);
\draw[usual,directed=0.99] (-5,0.5) to (-5,-0.5) to[out=270,in=0] (-6,-1) to[out=180,in=270] (-7,-0.5) to (-7,5);
\draw[usual,directed=0.99] (-2.5,0.5) to (-2.5,-0.5) to[out=270,in=0] (-6,-2) to[out=180,in=270] (-9.5,-0.5) to (-9.5,5);
\draw[usual,crossline,directed=0.99] (-2.5,2) to (-2.5,5);
\draw[spinach!35,fill=spinach!35] (-5.5,2) rectangle (-2,0.5);
\node at (-3.75,1.0) {$\morstuff{f}$};
\end{tikzpicture}
\end{gather*}
is an isomorphism. This easily generalizes.
\end{proof}

\begin{Definition}\label{D:WebsExplosion}
A left invertible morphism $\uob[k{+}l]\to\uob[k]\hcirc\uob[l]$ is called 
\emph{explosion}.
\end{Definition}

Up to rescaling, explosion morphisms and their left inverses can be illustrated by
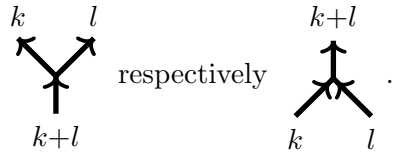
\begin{gather*}
\begin{tikzpicture}[anchorbase,scale=1,yscale=-1]
\draw[usual,directed=0.99] (0.5,0.5) to (0,0)node[above]{$k$};
\draw[usual,directed=0.99] (0.5,0.5) to (1,0)node[above]{$l$};
\draw[usual,directed=0.99] (0.5,1)node[below]{$k{+}l$} to (0.5,0.5);
\end{tikzpicture}
\text{ respectively }
\begin{tikzpicture}[anchorbase,scale=1]
\draw[usual,directed=0.9] (0,0)node[below]{$k$} to (0.5,0.5);
\draw[usual,directed=0.9] (1,0)node[below]{$l$} to (0.5,0.5);
\draw[usual,directed=0.99] (0.5,0.5) to (0.5,1)node[above]{$k{+}l$};
\end{tikzpicture}
.
\end{gather*}

\begin{Lemma}\label{L:WebsExplosion}
If the left mates of the $(k,l)$-overcrossings 
are invertible and all explosion morphisms exist,
then all hom-spaces in $\diacat$ are determined by
end-spaces between objects of the form $1^{\hcirc k}=\uob[1]\hcirc\dots\hcirc\uob[1]$ ($k$ factors).
\end{Lemma}

\begin{proof}
Turn all strands upwards using \autoref{L:WebsUpwards}, and then explode 
the strands inductively.
\end{proof}


\subsection{Exterior \texorpdfstring{$\mathfrak{gl}_{n}$}{gln}-webs}\label{SS:WebsGln}


We now recall the category of exterior $\gln$ webs.

\begin{Definition}\label{D:PrelimWebsGln}
The \emph{(exterior $\gln$) web category $\webv$} is the
diagram category for $R=\Zv$ with $\hcirc$-generating objects 
of categorical dimension $\qbin{n}{k}$
and $\vcirc$-$\hcirc$-generating morphisms
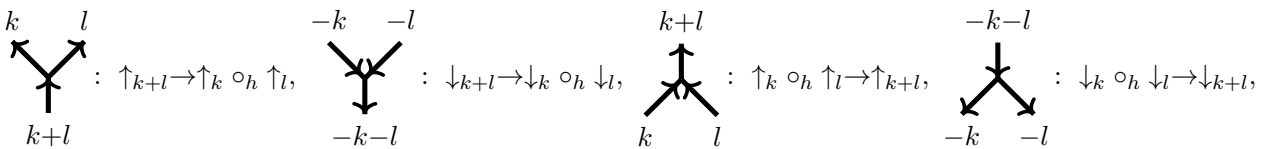
\begin{gather*}
\scalebox{0.95}{$\begin{tikzpicture}[anchorbase,scale=1,yscale=-1]
\draw[usual,directed=0.99] (0.5,0.5) to (0,0)node[above]{$k$};
\draw[usual,directed=0.99] (0.5,0.5) to (1,0)node[above]{$l$};
\draw[usual,directed=0.99] (0.5,1)node[below]{$k{+}l$} to (0.5,0.5);
\end{tikzpicture}
\hspace*{-0.1cm}\colon
\uob[k{+}l]\to\uob\hcirc\uob[l]
,
\begin{tikzpicture}[anchorbase,scale=1,yscale=-1]
\draw[usual,directed=0.9] (0,0)node[above]{${-}k$} to (0.5,0.5);
\draw[usual,directed=0.9] (1,0)node[above]{${-}l$} to (0.5,0.5);
\draw[usual,directed=0.99] (0.5,0.5) to (0.5,1)node[below]{${-}k{-}l$};
\end{tikzpicture}
\hspace*{-0.1cm}\colon
\dob[k{+}l]\to\dob\hcirc\dob[l]
,
\begin{tikzpicture}[anchorbase,scale=1]
\draw[usual,directed=0.9] (0,0)node[below]{$k$} to (0.5,0.5);
\draw[usual,directed=0.9] (1,0)node[below]{$l$} to (0.5,0.5);
\draw[usual,directed=0.99] (0.5,0.5) to (0.5,1)node[above]{$k{+}l$};
\end{tikzpicture}
\hspace*{-0.1cm}\colon
\uob\hcirc\uob[l]\to\uob[k{+}l]
,
\begin{tikzpicture}[anchorbase,scale=1]
\draw[usual,directed=0.99] (0.5,0.5) to (0,0)node[below]{${-}k$};
\draw[usual,directed=0.99] (0.5,0.5) to (1,0)node[below]{${-}l$};
\draw[usual,directed=0.99] (0.5,1)node[above]{${-}k{-}l$} to (0.5,0.5);
\end{tikzpicture}
\hspace*{-0.1cm}\colon
\dob\hcirc\dob[l]\to\dob[k{+}l]
,$}
\end{gather*}
such that the left mates of the $(k,l)$-overcrossings are invertible.
The relations imposed on $\webv$ are \emph{isotopies} 
(not displayed zigzag and trivalent-slide relations, see {\eg} 
\cite[Section 2]{LaTu-gln-webs} for details),
the \emph{exterior relation}, \emph{associativity}, \emph{coassociativity}, \emph{digon removal}, and 
\emph{dumbbell-crossing relation}.
That is, we take the quotient by the $\vcirc$-$\hcirc$-ideal generated by isotopies and
\begin{gather*}
\begin{tikzpicture}[anchorbase,scale=1]
\draw[usual] (0,0)node[below]{${>}n$} to (0,1.5)node[above]{${>}n$};
\end{tikzpicture}
=0
,\hspace*{-0.15cm}
\begin{tikzpicture}[anchorbase,scale=1]
\draw[usual] (0,0)node[below]{$k$} to (1,1);
\draw[usual] (1,0)node[below]{$l$} to (0.5,0.5);
\draw[usual] (2,0)node[below]{$m$} to (1,1);
\draw[usual,directed=0.99] (1,1) to (1,1.5)node[above]{$k{+}l{+}m$};
\end{tikzpicture}
\hspace*{-0.15cm}=\hspace*{-0.15cm}
\begin{tikzpicture}[anchorbase,scale=1]
\draw[usual] (0,0)node[below]{$k$} to (1,1);
\draw[usual] (1,0)node[below]{$l$} to (1.5,0.5);
\draw[usual] (2,0)node[below]{$m$} to (1,1);
\draw[usual,directed=0.99] (1,1) to (1,1.5)node[above]{$k{+}l{+}m$};
\end{tikzpicture}
\hspace*{-0.15cm},\hspace*{-0.15cm}
\begin{tikzpicture}[anchorbase,scale=1]
\draw[usual,rdirected=0.05] (0,0)node[above]{$k$} to (1,-1);
\draw[usual,rdirected=0.1] (1,0)node[above]{$l$} to (0.5,-0.5);
\draw[usual,rdirected=0.05] (2,0)node[above]{$m$} to (1,-1);
\draw[usual] (1,-1) to (1,-1.5)node[below]{$k{+}l{+}m$};
\end{tikzpicture}
\hspace*{-0.15cm}=\hspace*{-0.15cm}
\begin{tikzpicture}[anchorbase,scale=1]
\draw[usual,rdirected=0.05] (0,0)node[above]{$k$} to (1,-1);
\draw[usual,rdirected=0.1] (1,0)node[above]{$l$} to (1.5,-0.5);
\draw[usual,rdirected=0.05] (2,0)node[above]{$m$} to (1,-1);
\draw[usual] (1,-1) to (1,-1.5)node[below]{$k{+}l{+}m$};
\end{tikzpicture}
\hspace*{-0.15cm},
\begin{tikzpicture}[anchorbase,scale=1,rounded corners]
\draw[usual] (0.5,0.35) to (0,0.75)node[left]{$k$} to (0.5,1.15);
\draw[usual] (0.5,0.35) to (1,0.75)node[right]{$l$} to (0.5,1.15);
\draw[usual,directed=0.99] (0.5,1.15) to (0.5,1.5)node[above]{$k{+}l$};
\draw[usual] (0.5,0.35) to (0.5,0)node[below]{$k{+}l$};
\end{tikzpicture}
=\qbinn{k+l}{k}
\cdot\hspace*{-0.15cm}
\begin{tikzpicture}[anchorbase,scale=1]
\draw[usual,directed=0.99] (0.5,0)node[below]{$k{+}l$} to (0.5,1.5)node[above]{$k{+}l$};
\end{tikzpicture}
,\\
\begin{tikzpicture}[anchorbase,scale=1]
\draw[usual,directed=0.99] (0,0)node[below]{$k$} to (0.5,0.5) to (0.5,1) to (0,1.5)node[above]{$r$};
\draw[usual,directed=0.99] (1,0)node[below]{$l$} to (0.5,0.5) to (0.5,1) to (1,1.5)node[above]{$s$};
\end{tikzpicture}
=
(-1)^{kl}\sum_{k-r=a-b}
(-\vpar)^{(k-a)(l-b)}
\begin{tikzpicture}[anchorbase,scale=1]
\draw[usual] (1,0.4)node[left,yshift=-0.05cm]{$b$} to (0,1.1);
\draw[usual,crossline] (0,0.4)node[right,yshift=-0.1cm]{$a$} to (1,1.1);
\draw[usual,directed=0.99] (0,0)node[below]{$k$} to (0,1.5)node[above]{$r$};
\draw[usual,directed=0.99] (1,0)node[below]{$l$} to (1,1.5)node[above]{$s$};
\end{tikzpicture}
=
(-1)^{kl}\sum_{k-r=a-b}
(-\vpar)^{-(k-a)(l-b)}
\begin{tikzpicture}[anchorbase,scale=1]
\draw[usual] (0,0.4)node[right,yshift=-0.1cm]{$a$} to (1,1.1);
\draw[usual,crossline] (1,0.4)node[left,yshift=-0.05cm]{$b$} to (0,1.1);
\draw[usual,directed=0.99] (0,0)node[below]{$k$} to (0,1.5)node[above]{$r$};
\draw[usual,directed=0.99] (1,0)node[below]{$l$} to (1,1.5)node[above]{$s$};
\end{tikzpicture}
,
\end{gather*}
together with their $(\placeholder)^{\danti}$-duals.
\end{Definition}

We call morphisms in $\webv$ \emph{(exterior $\gln$) webs}.

\begin{Remark}\label{R:PrelimWebsGln}
The dumbbell-crossing relation is not new
and can be deduced from Green's book on the Schur algebra \cite{Gr-poly-reps} via 
an interpretation of webs as elements in the Schur algebra. 
Consequently, this relation is also called 
the \emph{Schur relation}.
\end{Remark}

It follows for example from \cite[Sections 2 and 5]{LaTu-gln-webs} that there 
is a well-defined functor from the pivotal category of 
topological webs (where webs are defined as plane labeled oriented trivalent graphs up to planar isotopy) to $\webv$ which we will use to draw webs in a topological fashion. 

\begin{Lemma}\label{L:WebsCrossings}
In $\webv$ we have
\begin{gather}\label{Eq:PrelimCrossing}
\begin{tikzpicture}[anchorbase,scale=1]
\draw[usual,directed=0.99] (1,0)node[below]{$l$} to (0,1)node[above]{$l$};
\draw[usual,directed=0.99,crossline] (0,0)node[below]{$k$} to (1,1)node[above]{$k$};
\end{tikzpicture}
=
(-1)^{kl}
\sum_{b-a=k-l}
(-\vpar)^{k-b}
\begin{tikzpicture}[anchorbase,scale=1]
\draw[usual,directed=0.99] (0,0)node[below]{$k$} to (0,1.5)node[above]{$l$};
\draw[usual,directed=0.99] (1,0)node[below]{$l$} to (1,1.5)node[above]{$k$};
\draw[usual] (0,0.4) to node[below]{$b$} (1,0.6);
\draw[usual] (1,0.9) to node[above]{$a$} (0,1.1);
\end{tikzpicture}
,\quad
\begin{tikzpicture}[anchorbase,scale=1]
\draw[usual,directed=0.99] (0,0)node[below]{$k$} to (1,1)node[above]{$k$};
\draw[usual,directed=0.99,crossline] (1,0)node[below]{$l$} to (0,1)node[above]{$l$};
\end{tikzpicture}
=
(-1)^{kl}
\sum_{b-a=k-l}
(-\vpar)^{-k+b}
\begin{tikzpicture}[anchorbase,scale=1]
\draw[usual,directed=0.99] (0,0)node[below]{$k$} to (0,1.5)node[above]{$l$};
\draw[usual,directed=0.99] (1,0)node[below]{$l$} to (1,1.5)node[above]{$k$};
\draw[usual] (0,0.4) to node[below]{$b$} (1,0.6);
\draw[usual] (1,0.9) to node[above]{$a$} (0,1.1);
\end{tikzpicture}
.
\end{gather}
\end{Lemma}	

\begin{proof}
This is explained in \cite[Section 5]{LaTu-gln-webs}.
\end{proof}

For completeness, the $k=l=1$ case of \autoref{Eq:PrelimCrossing} is
\begin{gather}\label{Eq:PrelimThinCrossing}
\begin{tikzpicture}[anchorbase,scale=1]
\draw[usual,directed=0.99] (1,0)node[below]{$1$} to (0,1)node[above]{$1$};
\draw[usual,directed=0.99,crossline] (0,0)node[below]{$1$} to (1,1)node[above]{$1$};
\end{tikzpicture}
=
\vpar\cdot
\begin{tikzpicture}[anchorbase,scale=1]
\draw[usual,directed=0.99] (0,0)node[below]{$1$} to (0,1)node[above]{$1$};
\draw[usual,directed=0.99] (1,0)node[below]{$1$} to (1,1)node[above]{$1$};
\end{tikzpicture}
-
\begin{tikzpicture}[anchorbase,scale=1]
\draw[usual] (0,0)node[below]{$1$} to (0.5,0.25);
\draw[usual] (1,0)node[below]{$1$} to (0.5,0.25);
\draw[usual] (0.5,0.25) to (0.5,0.75);
\draw[usual,directed=0.99] (0.5,0.75) to (0,1)node[above]{$1$};
\draw[usual,directed=0.99] (0.5,0.75) to (1,1)node[above]{$1$};
\end{tikzpicture}
,\quad
\begin{tikzpicture}[anchorbase,scale=1]
\draw[usual,directed=0.99] (0,0)node[below]{$1$} to (1,1)node[above]{$1$};
\draw[usual,directed=0.99,crossline] (1,0)node[below]{$1$} to (0,1)node[above]{$1$};
\end{tikzpicture}
=
\vpar^{-1}\cdot
\begin{tikzpicture}[anchorbase,scale=1]
\draw[usual,directed=0.99] (0,0)node[below]{$1$} to (0,1)node[above]{$1$};
\draw[usual,directed=0.99] (1,0)node[below]{$1$} to (1,1)node[above]{$1$};
\end{tikzpicture}
-
\begin{tikzpicture}[anchorbase,scale=1]
\draw[usual] (0,0)node[below]{$1$} to (0.5,0.25);
\draw[usual] (1,0)node[below]{$1$} to (0.5,0.25);
\draw[usual] (0.5,0.25) to (0.5,0.75);
\draw[usual,directed=0.99] (0.5,0.75) to (0,1)node[above]{$1$};
\draw[usual,directed=0.99] (0.5,0.75) to (1,1)node[above]{$1$};
\end{tikzpicture}
.
\end{gather}

\begin{Lemma}\label{L:WebsGlnProperties}
The following hold in $\webv$.
\begin{enumerate}

\item The crossings satisfy the \emph{Reidemeister II and III relations}, and the \emph{Reidemeister I relation} holds up to scalars, that is,
\begin{gather}\label{Eq:WebsRM1}
\begin{tikzpicture}[anchorbase,scale=1]
\draw[usual,directed=0.99] (1,0) to[out=270,in=0] (0.5,-0.35) to[out=180,in=270] (0,1)node[above]{$k$};
\draw[usual,crossline] (0,-1)node[below]{$k$} to[out=90,in=180] (0.5,0.35) to[out=0,in=90] (1,0);
\end{tikzpicture}
=
\vpar^{k(-k+n+1)}\cdot
\begin{tikzpicture}[anchorbase,scale=1]
\draw[usual,crossline,directed=0.99] (0,-1)node[below]{$k$} to (0,1)node[above]{$k$};
\end{tikzpicture}
=
\begin{tikzpicture}[anchorbase,scale=1,xscale=-1]
\draw[usual] (0,-1)node[below]{$k$} to[out=90,in=180] (0.5,0.35) to[out=0,in=90] (1,0);
\draw[usual,crossline,directed=0.99] (1,0) to[out=270,in=0] (0.5,-0.35) to[out=180,in=270] (0,1)node[above]{$k$};
\end{tikzpicture}
,\quad
\begin{tikzpicture}[anchorbase,scale=1]
\draw[usual] (0,-1)node[below]{$k$} to[out=90,in=180] (0.5,0.35) to[out=0,in=90] (1,0);
\draw[usual,crossline,directed=0.99] (1,0) to[out=270,in=0] (0.5,-0.35) to[out=180,in=270] (0,1)node[above]{$k$};
\end{tikzpicture}
=
\vpar^{k(k-n-1)}\cdot
\begin{tikzpicture}[anchorbase,scale=1]
\draw[usual,crossline,directed=0.99] (0,-1)node[below]{$k$} to (0,1)node[above]{$k$};
\end{tikzpicture}
=
\begin{tikzpicture}[anchorbase,scale=1,xscale=-1]
\draw[usual,directed=0.99] (1,0) to[out=270,in=0] (0.5,-0.35) to[out=180,in=270] (0,1)node[above]{$k$};
\draw[usual,crossline] (0,-1)node[below]{$k$} to[out=90,in=180] (0.5,0.35) to[out=0,in=90] (1,0);
\end{tikzpicture}
,
\end{gather}
together with their $(\placeholder)^{\danti}$-duals.
Various \emph{naturality relations} hold, see {\eg} \cite[Section 2]{LaTu-gln-webs}.

\item \emph{Square switches}, see {\eg} \cite[Lemma 5.6]{LaTu-gln-webs}. Various other relations that we do not explicitly use, see {\eg} \cite[Section 2]{LaTu-gln-webs}, also hold.

\item \emph{Explosion} holds for $\webstar$ (but not in $\webv$), that is
\begin{gather*}
\begin{tikzpicture}[anchorbase,scale=1]
\draw[usual,directed=0.99] (0.5,0)node[below]{$k{+}l$} to (0.5,1.5)node[above]{$k{+}l$};
\end{tikzpicture}
=\qbinn{k+l}{k}^{-1}
\cdot
\begin{tikzpicture}[anchorbase,scale=1,rounded corners]
\draw[usual] (0.5,0.35) to (0,0.75)node[left]{$k$} to (0.5,1.15);
\draw[usual] (0.5,0.35) to (1,0.75)node[right]{$l$} to (0.5,1.15);
\draw[usual,directed=0.99] (0.5,1.15) to (0.5,1.5)node[above]{$k{+}l$};
\draw[usual] (0.5,0.35) to (0.5,0)node[below]{$k{+}l$};
\end{tikzpicture}
,
\end{gather*}
together with its $(\placeholder)^{\danti}$-dual. Moreover, the $(k,l)$-overcrossings satisfy explosion as well, {\ie}
\begin{gather}\label{Eq:WebsExplosion}
\begin{tikzpicture}[anchorbase,scale=1]
\draw[usual] (1,0)node[below]{$l$} to (0,1)node[above]{$l$};
\draw[usual,crossline] (0,0)node[below]{$k$} to (1,1)node[above]{$k$};
\end{tikzpicture}
=\qfac{k}^{-1}\qfac{l}^{-1}\cdot
\scalebox{0.75}{$\begin{tikzpicture}[anchorbase,scale=1]
\draw[spinach!50,fill=spinach!50] (-1.25,1) rectangle (2.25,1.5);
\draw[usual] (-0.5,0)node[below]{$k$} to (-0.5,0.5) to (-1,1);
\draw[usual] (-0.5,.5) to (0,1);
\draw[usual] (-1,1.5) to (-0.5,2) to (-0.5,2.5)node[above]{$l$};
\draw[usual] (0,1.5) to (-0.5,2);
\draw[usual] (1.5,0)node[below]{$l$} to (1.5,0.5) to (1,1);
\draw[usual] (1.5,.5) to (2,1);
\draw[usual] (1,1.5) to (1.5,2) to (1.5,2.5)node[above]{$k$};
\draw[usual] (2,1.5) to (1.5,2);
\node at (-0.5,0.8) {$\dots$};
\node at (0.5,1.25) {$\morstuff{x}$};
\node at (-0.5,1.7) {$\dots$};
\node at (1.5,0.8) {$\dots$};
\node at (1.5,1.7) {$\dots$};
\end{tikzpicture}$}
,\quad
\morstuff{x}
=
\begin{tikzpicture}[anchorbase,scale=1]
\draw[usual] (2,0) to (0,2);
\draw[usual] (3,0) to (1,2)node[above,xshift=-0.5cm]{\text{$l$ strands}};
\draw[usual,crossline] (0,0) to (2,2);
\draw[usual,crossline] (1,0)node[below,xshift=-0.5cm]{\text{$k$ strands}} to (3,2);
\node at (0.55,0) {$\dots$};
\node at (2.45,0) {$\dots$};
\node at (0.55,2) {$\dots$};
\node at (2.45,2) {$\dots$};
\end{tikzpicture}
,
\end{gather}
as well as a similar formula for the $(k,l)$-undercrossings.\qed

\end{enumerate}
\end{Lemma}

Note that thus \autoref{L:WebsTangle} applies in $\webstar$.
Let us also note the following, partially explaining why explosion 
works well in practice:

\begin{Lemma}\label{L:WebsSemisimple}
The additive idempotent completion of $\webstar$ is semisimple.
\end{Lemma}

\begin{proof}
By \autoref{L:WebsTangle} and the existence of certain projectors, 
see {\eg} \cite{RoTu-symmetric-howe}
or \cite[Section 2.3]{TuVaWe-super-howe}. To be precise (and to fix notation for the rest of the paper), the projectors are
\begin{gather*}
\qfac{k}^{-1}
\cdot
\begin{tikzpicture}[anchorbase,scale=1]
\draw[usual] (0,0)node[below]{$1$} to (0.5,0.25);
\draw[usual] (1,0)node[below]{$1$} to (0.5,0.25);
\draw[usual] (0.5,0.25) to (0.5,0.5)node[right]{$k$} to (0.5,0.75);
\draw[usual,directed=0.99] (0.5,0.75) to (0,1)node[above]{$1$};
\draw[usual,directed=0.99] (0.5,0.75) to (1,1)node[above]{$1$};
\node at (0.5,0) {$\dots$};
\node at (0.5,1) {$\dots$};
\end{tikzpicture}
,\quad
\begin{tikzpicture}[anchorbase,scale=1]
\draw[usual] (0.5,0.5) to (0,0)node[below]{$1$};
\draw[usual] (0.5,0.5) to (1,0)node[below]{$1$};
\draw[usual] (0.5,1)node[above]{$k$} to (0.5,0.5);
\node at (0.5,0) {$\dots$};
\end{tikzpicture}
=
\begin{tikzpicture}[anchorbase,scale=1]
\draw[usual] (0,0)node[below]{$1$} to (0.75,0.75);
\draw[usual] (0.5,0)node[below]{$1$} to (0.25,0.25);
\draw[usual] (1,0)node[below]{$1$} to (0.5,0.5);
\draw[usual] (1.5,0)node[below]{$1$} to (0.75,0.75);
\draw[usual] (0.75,0.75) to (0.75,1)node[above]{$k$};
\node at (0.75,0) {$\dots$};
\end{tikzpicture}
,
\end{gather*}
where the dots indicate an inductive construction as illustrated on the right 
(the order of how these are constructed is irrelevant due 
to associativity and coassociativity).
\end{proof}


\section{Annular webs}\label{S:AnnularWebs}


This section discusses our main diagram categories of this paper. 
Similar constructions have appeared in many texts, 
{\eg} \cite{CaKa-q-satake-sln} or \cite{QuWe-extremal-projectors-2}.


\subsection{The annular web category}\label{SS:AnnularWebs}


The following definition does not use any $\hcirc$ structure.

\begin{Definition}\label{D:AnnularWebsDefinition}
The \emph{(exterior $\gln$) annular web category $\awebv$} is
the category obtained from $\webv$ by adding extra 
$\vcirc$-generators
\begin{gather*}
\winding=
\begin{tikzpicture}[anchorbase,scale=1]
\draw[usual,rounded corners] (0,0)node[below]{$k_{1}$} to (0,0.5) to (-0.5,0.5);
\draw[usual,rounded corners] (2.5,0.5) to (2,0.5) to (2,1)node[above]{$k_{1}$};
\draw[usual] (1,0)node[below]{$k_{2}$} to (0,1)node[above]{$k_{2}$};
\draw[usual] (2,0)node[below]{$k_{m}$} to (1,1)node[above]{$k_{m}$};
\node at (1,0.5) {$\dots$};
\draw[affine] (-0.5,0) to (-0.5,1);
\draw[affine] (2.5,0) to (2.5,1);
\end{tikzpicture}
\,,\quad
\iwinding=
\begin{tikzpicture}[anchorbase,scale=1]
\draw[usual,rounded corners] (0,0)node[above]{$k_{1}$} to (0,-0.5) to (-0.5,-0.5);
\draw[usual,rounded corners] (2.5,-0.5) to (2,-0.5) to (2,-1)node[below]{$k_{1}$};
\draw[usual] (1,0)node[above]{$k_{2}$} to (0,-1)node[below]{$k_{2}$};
\draw[usual] (2,0)node[above]{$k_{m}$} to (1,-1)node[below]{$k_{m}$};
\node at (1,-0.5) {$\dots$};
\draw[affine] (-0.5,0) to (-0.5,-1);
\draw[affine] (2.5,0) to (2.5,-1);
\end{tikzpicture}
\,,
\end{gather*}
for each 
$\obstuff{K}=(k_{1},\dots,k_{m})\in\Z^{m}$ to $\webv$ modulo the $\vcirc$-ideal 
generated by the relations
\begin{gather}\label{Eq:AnnularWebsRel1}
\begin{tikzpicture}[anchorbase,scale=1,rounded corners]
\draw[usual,rounded corners] (0,-0.01) to (0,0) to (0,0.5) to (-0.5,0.5);
\draw[usual,rounded corners] (2.5,0.5) to (2,0.5) to (2,1);
\draw[usual,rounded corners] (0,-1) to (1,0) to (0,1);
\draw[usual,rounded corners] (1,-1) to (2,0) to (1,1);
\draw[usual,rounded corners] (0,0) to (0,-0.5) to (-0.5,-0.5);
\draw[usual,rounded corners] (2.5,-0.5) to (2,-0.5) to (2,-1);
\node at (1,0.5) {$\dots$};
\node at (1,-0.5) {$\dots$};
\draw[affine] (-0.5,-1) to (-0.5,1);
\draw[affine] (2.5,-1) to (2.5,1);
\end{tikzpicture}
=
\begin{tikzpicture}[anchorbase,scale=1]
\draw[usual] (0,-1) to (0,1);
\draw[usual] (1,-1) to (1,1);
\draw[usual] (2,-1) to (2,1);
\node at (0.5,0) {$\dots$};
\draw[affine] (-0.5,-1) to (-0.5,1);
\draw[affine] (2.5,-1) to (2.5,1);
\end{tikzpicture}	
,\quad
\begin{tikzpicture}[anchorbase,scale=1,xscale=-1,rounded corners]
\draw[usual,rounded corners] (0,-0.01) to (0,0) to (0,0.5) to (-0.5,0.5);
\draw[usual,rounded corners] (2.5,0.5) to (2,0.5) to (2,1);
\draw[usual,rounded corners] (0,-1) to (1,0) to (0,1);
\draw[usual,rounded corners] (1,-1) to (2,0) to (1,1);
\draw[usual,rounded corners] (0,0) to (0,-0.5) to (-0.5,-0.5);
\draw[usual,rounded corners] (2.5,-0.5) to (2,-0.5) to (2,-1);
\node at (1,0.5) {$\dots$};
\node at (1,-0.5) {$\dots$};
\draw[affine] (-0.5,-1) to (-0.5,1);
\draw[affine] (2.5,-1) to (2.5,1);
\end{tikzpicture}
=
\begin{tikzpicture}[anchorbase,scale=1]
\draw[usual] (2,-1) to (2,1);
\draw[usual] (0,-1) to (0,1);
\draw[usual] (1,-1) to (1,1);
\node at (1.5,0) {$\dots$};
\draw[affine] (-0.5,-1) to (-0.5,1);
\draw[affine] (2.5,-1) to (2.5,1);
\end{tikzpicture}
\,,
\end{gather}
\begin{gather}\label{Eq:AnnularWebsRel2}
\begin{tikzpicture}[anchorbase,scale=1]
\draw[usual] (0,-1) to (0.5,-0.5) to (1,-1);
\draw[usual,rounded corners] (0.5,-0.5) to (0.5,0) to (-0.5,0);
\draw[usual,rounded corners] (2.5,1) to (2.5,0) to (3.5,0);
\draw[usual] (1.5,-1)node[below]{$\obstuff{K}$} to (1.5,1)node[above]{$\obstuff{K}$};
\draw[affine] (-0.5,-1) to (-0.5,1);
\draw[affine] (3.5,-1) to (3.5,1);
\end{tikzpicture}
=
\begin{tikzpicture}[anchorbase,scale=1]
\draw[usual,rounded corners] (0,-1) to (0,-0.25) to (-0.5,-0.25);
\draw[usual,rounded corners] (1,-1) to (1,0) to (-0.5,0);
\draw[usual] (2.5,1) to (2.5,0.5) to (3,0) to (3.5,0);
\draw[usual] (2.5,0.5) to (2,0) to (2,-0.02);
\draw[usual,rounded corners] (2,0) to (2,-0.25) to (3.5,-0.25);
\draw[usual] (1.5,-1)node[below]{$\obstuff{K}$} to (1.5,1)node[above]{$\obstuff{K}$};
\draw[affine] (-0.5,-1) to (-0.5,1);
\draw[affine] (3.5,-1) to (3.5,1);
\end{tikzpicture}
\,,
\end{gather}
\begin{gather}\label{Eq:AnnularWebsRel3}
\begin{tikzpicture}[anchorbase,scale=1]
\draw[usual] (0,-1) to[out=90,in=180] (0.5,-0.5) to[out=0,in=90] (1,-1);
\draw[usual] (1.5,-1)node[below]{$\obstuff{K}$} to (1.5,1)node[above]{$\obstuff{K}$};
\draw[affine] (-0.5,-1) to (-0.5,1);
\draw[affine] (3.5,-1) to (3.5,1);
\end{tikzpicture}
=
\begin{tikzpicture}[anchorbase,scale=1]
\draw[usual,rounded corners] (0,-1) to (0,-0.25) to (-0.5,-0.25);
\draw[usual,rounded corners] (1,-1) to (1,0) to (-0.5,0);
\draw[usual] (2,-0.02) to(2,0) to[out=90,in=180] (2.5,0.5) to[out=0,in=90] (3,0) to (3.5,0);
\draw[usual,rounded corners] (2,0) to (2,-0.25) to (3.5,-0.25);
\draw[usual] (1.5,-1)node[below]{$\obstuff{K}$} to (1.5,1)node[above]{$\obstuff{K}$};
\draw[affine] (-0.5,-1) to (-0.5,1);
\draw[affine] (3.5,-1) to (3.5,1);
\end{tikzpicture}
\,,
\end{gather}
together with the $(\placeholder)^{\danti}$- and 
$(\placeholder)^{\dinvo}$-duals of the bottom two relations.
\end{Definition}

We call morphisms in $\awebv$ \emph{annular (exterior $\gln$) webs}, 
and $\winding$ and $\iwinding$ are called \emph{coils}.

\begin{Remark}\label{R:AnnularWebsCrossings}
We think of the coils as 
crossings in front of the annulus, {\eg}
\begin{gather*}
\winding[(1,1)]
\leftrightsquigarrow
\begin{tikzpicture}[anchorbase,scale=1]
\draw[usual,rounded corners] (0,0)node[below]{$1$} to (0,0.25) to (-0.5,0.25);
\draw[usual,rounded corners,directed=0.99] (2.5,0.75) to (2,0.75) to (2,1)node[above]{$1$};
\draw[usual,directed=0.99] (1,0)node[below]{$1$} to (1,1)node[above]{$1$};
\draw[very thick,crossline] (-0.5,0.25) to (-0.6,0.25) to[out=180,in=180] (-0.6,0.5) to (2.6,0.5) to[out=0,in=0] (2.6,0.75) to (2.5,0.75);
\draw[affine] (-0.5,0) to (-0.5,1);
\draw[affine] (2.5,0) to (2.5,1);
\end{tikzpicture}
.
\end{gather*}
This convention comes because we follow 
\cite{SaSh-schur-weyl-ariki-koike} later on for computations. 
Using the inverse braiding compared to the definitions in 
\cite[Theorem 3.2]{SaSh-schur-weyl-ariki-koike} translate to coils passing behind the annulus.
\end{Remark}

\begin{Remark}\label{R:AnnularWebsDefinition2}
The name annular webs comes from the interpretation of 
the pictures in \autoref{D:AnnularWebsDefinition} as embedded 
in an annulus. For example,
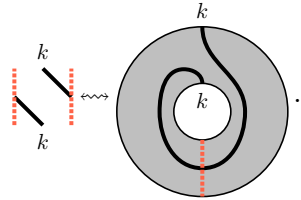
\begin{gather*}
\scalebox{0.75}{$\begin{tikzpicture}[anchorbase,scale=1]
\draw[usual] (0,0)node[below]{$k$} to (-0.5,0.5);
\draw[usual] (0.5,0.5) to (0,1)node[above]{$k$};
\draw[affine] (-0.5,0) to (-0.5,1);
\draw[affine] (0.5,0) to (0.5,1);
\end{tikzpicture}
\leftrightsquigarrow
\begin{tikzpicture}[anchorbase,scale=1]
\draw[thick,fill=gray!50] (0,0) circle (1.5cm);
\draw[thick,fill=white] (0,0) circle (0.5cm);
\draw[usual] (0,0.5)node[below]{$k$} to[out=90,in=0] (-0.325,0.75) to[out=180,in=90] (-0.75,0) to[out=270,in=180] (0,-1) to[out=0,in=270] (0.75,0)to[out=90,in=270] (0,1.5)node[above]{$k$};
\draw[affine] (0,-0.5) to (0,-1.5);
\end{tikzpicture}$}
\,.
\end{gather*}
\end{Remark}

The following elements defined via Lagrange interpolation play 
a crucial role later on:

\begin{Definition}\label{D:AnnularWebsProjectors}
For $i\in\set{1,\dots,d}$ define 
\begin{gather*}
\wproj=
\prod_{j\neq i}
\frac{\winding[1]-\upar[j]}{\upar[i]-\upar[j]}
=\prod_{j\neq i}
\frac{
\begin{tikzpicture}[anchorbase,scale=0.5]
\draw[usual,rounded corners] (0,0) to (0,0.5) to (-0.5,0.5);
\draw[usual,rounded corners] (0.5,0.5) to (0,0.5) to (0,1);
\draw[affine] (-0.5,0) to (-0.5,1);
\draw[affine] (0.5,0) to (0.5,1);
\end{tikzpicture}
-\upar[j]}{\upar[i]-\upar[j]}
\in\End_{\awebstar}(1),
\end{gather*}
which we call \emph{web block projectors}.
\end{Definition}


\subsection{Properties of annular webs}\label{SS:AnnularWebsProps}


We can endow $\awebv$ with a monoidal structure $\wprod$ as follows.
On objects 
$\obstuff{K}\hcirc\obstuff{K}^{\prime}=
\obstuff{K}\wprod\obstuff{K}^{\prime}$ is just the concatenation, 
{\ie} if $\obstuff{K}=(k_{1},\dots,k_{r})$ 
and $\obstuff{K}^{\prime}=(k_{1}^{\prime},\dots,k_{s}^{\prime})$, then 
$\obstuff{K}\hcirc\obstuff{K}^{\prime}=(k_{1},\dots,k_{r},k_{1}^{\prime},\dots,k_{s}^{\prime})$. 
On morphisms $\wprod$ we use
\begin{gather}\label{Eq:AnnularWebsPropsMonoidal}
\begin{tikzpicture}[anchorbase,scale=1,yscale=1]
\draw[spinach!35,fill=spinach!35] (-0.5,-0.25) to (-0.5,0.25) to (1.5,0.25) to (1.5,-0.25) to (-0.5,-0.25);
\draw[usual] (0,-1) to (0,-0.25);
\draw[usual] (0,0.25) to (0,1);
\draw[usual] (1,-1) to (1,-0.25);
\draw[usual] (1,0.25) to (1,1);
\node at (0.5,-0.625) {$\dots$};
\node at (0.5,0.625) {$\dots$};
\draw[affine] (-0.5,-1) to (-0.5,1);
\draw[affine] (1.5,-1) to (1.5,1);
\node at (0.5,0) {$\morstuff{f}$};
\end{tikzpicture}	
\wprod
\begin{tikzpicture}[anchorbase,scale=1,yscale=1]
\draw[spinach!35,fill=spinach!35] (-0.5,-0.25) to (-0.5,0.25) to (1.5,0.25) to (1.5,-0.25) to (-0.5,-0.25);
\draw[usual] (0,-1) to (0,-0.25);
\draw[usual] (0,0.25) to (0,1);
\draw[usual] (1,-1) to (1,-0.25);
\draw[usual] (1,0.25) to (1,1);
\node at (0.5,-0.625) {$\dots$};
\node at (0.5,0.625) {$\dots$};
\draw[affine] (-0.5,-1) to (-0.5,1);
\draw[affine] (1.5,-1) to (1.5,1);
\node at (0.5,0) {$\morstuff{g}$};
\end{tikzpicture}
=
\begin{tikzpicture}[anchorbase,scale=1,yscale=1]
\draw[spinach!35,fill=spinach!35] (-0.5,0.25) to (3.5,0.25) to (3.5,0.75) to (-0.5,0.75) to (-0.5,0.25);
\draw[spinach!35,fill=spinach!35] (-0.5,-0.25) to (3.5,-0.25) to (3.5,-0.75) to (-0.5,-0.75) to (-0.5,-0.25);
\draw[usual] (0,-1) to (0,-0.75);
\draw[usual,crossline] (0,-0.25) to (0,1);
\draw[usual] (1,-1) to (1,-0.75);
\draw[usual,crossline] (1,-0.25) to (1,1);
\draw[usual] (2,-1) to (2,-0.85);
\draw[usual] (2,-0.15) to (2,0.25);
\draw[usual] (2,0.75) to (2,1);
\draw[usual] (3,-1) to (3,-0.85);
\draw[usual] (3,-0.15) to (3,0.25);
\draw[usual] (3,0.75) to (3,1);
\node at (0.5,-1) {$\dots$};
\node at (0.5,1) {$\dots$};
\node at (2.5,-1) {$\dots$};
\node at (2.5,1) {$\dots$};
\draw[affine] (-0.5,-1) to (-0.5,1);
\draw[affine] (3.5,-1) to (3.5,1);
\node at (0.5,-0.5) {$\morstuff{f}$};
\node at (2.5,0.5) {$\morstuff{g}$};
\end{tikzpicture}
\,,
\end{gather}
using the $(k,l)$-crossings from \autoref{Eq:PrelimCrossing} 
and their mates so that $\morstuff{f}$ 
is in the front and $\morstuff{g}$ is in the back.

\begin{Remark}\label{R:AnnularWebsMonoidal}
\autoref{Eq:AnnularWebsPropsMonoidal} is a standard construction in skein theory, see {\eg} \cite{PrSi-skein-algebra}.
\end{Remark}

\begin{Lemma}\label{L:AnnularWebsMonoidal}
The monoidal structure
$\wprod$ 
and $\munit=\emptyset$ endows $\awebv$
with the structure of a pivotal category with duality given by cups and caps.
\end{Lemma}

\begin{proof}
Easy and omitted.
\end{proof}

\begin{Lemma}\label{L:AnnularWebsProperties}
The following holds in $\awebv$.
\begin{enumerate}

\item We have an \emph{annular digon removal}, that is
\begin{gather*}
\begin{tikzpicture}[anchorbase,scale=1]
\draw[usual] (0.5,-1)node[below]{$k{+}l$} to (0.5,-0.5) to (0,0) to (-0.5,0)node[left,yshift=-0.05cm]{$k$};
\draw[usual,rounded corners] (0.5,-0.5) to (1,0) to (1,0.25) to (-0.5,0.25)node[left,yshift=0.05cm]{$l$};
\draw[usual] (2.5,1)node[above]{$k{+}l$} to (2.5,0.5) to (3,0) to (3.5,0)node[right,yshift=0.05cm]{$l$};
\draw[usual,rounded corners] (2.5,0.5) to (2,0) to (2,-0.25) to (3.5,-0.25)node[right,yshift=-0.05cm]{$k$};
\draw[usual] (1.5,-1)node[below]{$\obstuff{K}$} to (1.5,1)node[above]{$\obstuff{K}$};
\draw[affine] (-0.5,-1) to (-0.5,1);
\draw[affine] (3.5,-1) to (3.5,1);
\end{tikzpicture}
=
\qbinn{k+l}{k}
\cdot
\begin{tikzpicture}[anchorbase,scale=1]
\draw[usual,rounded corners] (0.5,-1)node[below]{$k{+}l$} to (0.5,0) to (-0.5,0);
\draw[usual,rounded corners] (2.5,1)node[above]{$k{+}l$} to (2.5,0) to (3.5,0);
\draw[usual] (1.5,-1)node[below]{$\obstuff{K}$} to (1.5,1)node[above]{$\obstuff{K}$};
\draw[affine] (-0.5,-1) to (-0.5,1);
\draw[affine] (3.5,-1) to (3.5,1);
\end{tikzpicture}
\,,
\end{gather*}
together with its $(\placeholder)^{\danti}$-dual.
Various other annular versions of the relations in \autoref{L:WebsGlnProperties} 
hold as well (but are not stated since we do not use them).

\item All \emph{half-slides} of merges, splits, cups and caps, {\eg}
\begin{gather*}
\begin{tikzpicture}[anchorbase,scale=1]
\draw[usual] (0.5,-1) to (0.5,-0.5) to (0,0) to (-0.5,0);
\draw[usual,rounded corners] (0.5,-0.5) to (1,0) to (1,1);
\draw[usual,rounded corners] (2,1) to (2,0) to (3.5,0);
\draw[usual] (1.5,-1)node[below]{$\obstuff{K}$} to (1.5,1)node[above]{$\obstuff{K}$};
\draw[affine] (-0.5,-1) to (-0.5,1);
\draw[affine] (3.5,-1) to (3.5,1);
\end{tikzpicture}
=
\begin{tikzpicture}[anchorbase,scale=1]
\draw[usual,rounded corners] (0.5,-1) to (0.5,-0.5) to (-0.5,-0.5);
\draw[usual,rounded corners] (1,1) to (1,0.5) to (-0.5,0.5);
\draw[usual] (2.5,0) to (3,0.5) to (3.5,0.5);
\draw[usual] (2,1) to (2,0.5) to (2.5,0);
\draw[usual,rounded corners] (2.5,0) to (2.5,-0.5) to (3.5,-0.5);
\draw[usual] (1.5,-1)node[below]{$\obstuff{K}$} to (1.5,1)node[above]{$\obstuff{K}$};
\draw[affine] (-0.5,-1) to (-0.5,1);
\draw[affine] (3.5,-1) to (3.5,1);
\end{tikzpicture}
\,.
\end{gather*}

\item \emph{Annular explosion} 
holds for $\awebstar$ (but not in $\awebv$), that is
\begin{gather*}
\begin{tikzpicture}[anchorbase,scale=1]
\draw[usual,rounded corners] (0.5,-1)node[below]{$k{+}l$} to (0.5,0) to (-0.5,0);
\draw[usual,rounded corners] (2.5,1)node[above]{$k{+}l$} to (2.5,0) to (3.5,0);
\draw[usual] (1.5,-1)node[below]{$\obstuff{K}$} to (1.5,1)node[above]{$\obstuff{K}$};
\draw[affine] (-0.5,-1) to (-0.5,1);
\draw[affine] (3.5,-1) to (3.5,1);
\end{tikzpicture}
=\qbinn{k+l}{k}^{-1}
\cdot
\begin{tikzpicture}[anchorbase,scale=1]
\draw[usual] (0.5,-1)node[below]{$k{+}l$} to (0.5,-0.5) to (0,0) to (-0.5,0)node[left,yshift=-0.05cm]{$k$};
\draw[usual,rounded corners] (0.5,-0.5) to (1,0) to (1,0.25) to (-0.5,0.25)node[left,yshift=0.05cm]{$l$};
\draw[usual] (2.5,1)node[above]{$k{+}l$} to (2.5,0.5) to (3,0) to (3.5,0)node[right,yshift=0.05cm]{$l$};
\draw[usual,rounded corners] (2.5,0.5) to (2,0) to (2,-0.25) to (3.5,-0.25)node[right,yshift=-0.05cm]{$k$};
\draw[usual] (1.5,-1)node[below]{$\obstuff{K}$} to (1.5,1)node[above]{$\obstuff{K}$};
\draw[affine] (-0.5,-1) to (-0.5,1);
\draw[affine] (3.5,-1) to (3.5,1);
\end{tikzpicture}
\,
\end{gather*}
together with its $(\placeholder)^{\danti}$-dual.

\end{enumerate}
\end{Lemma}

\begin{proof}
We get
\begin{gather}\label{Eq:AnnularWebsProperties}
\begin{tikzpicture}[anchorbase]
\draw[usual] (0,0)node[below]{$k$} to[out=90,in=0] (-0.5,0.5);
\draw[usual,directed=1] (0.5,0.5) to[out=180,in=270] (0,1)node[above]{$k$};
\draw[affine] (-0.5,0) to (-0.5,1);
\draw[affine] (0.5,0) to (0.5,1);
\end{tikzpicture}
=\frac{1}{\qfac{k}}\cdot
\begin{tikzpicture}[anchorbase]
\draw[usual] (2.5,-0.5) to (2,0) to[out=90,in=0] (1.5,0.5)node[below,xshift=0.12cm]{$1$}node[above,xshift=0.12cm,yshift=-0.05cm]{$\vdots$};
\draw[usual] (2.5,-0.5) to (3,0) to[out=90,in=0] (1.5,1)node[above,xshift=0.12cm]{$1$};	
\draw[usual] (3,2) to  (2.5,1.5) to[out=270,in=180] (4,0.5)node[below,xshift=-0.12cm]{$1$}node[above,xshift=-0.12cm,yshift=-0.05cm]{$\vdots$};
\draw[usual] (3,2) to (3.5,1.5) to[out=270,in=180] (4,1)node[above,xshift=-0.12cm]{$1$};	
\draw[usual,directed=1] (3,2) to (3,2.5)node[above]{$k$};
\draw[usual] (2.5,-1)node[below]{$k$} to (2.5,-0.5);
\draw[affine] (1.5,-1) to (1.5,2.5);
\draw[affine] (4,-1) to (4,2.5);
\end{tikzpicture}
\,,
\end{gather}
by (plain) explosion and \autoref{Eq:AnnularWebsRel2}.
The claimed relations can be proven using this.
\end{proof}

The thin coils suffice as $\vcirc$-generators:

\begin{Lemma}\label{L:AnnularExplosion}
The morphisms $\winding[(k,\obstuff{K})]$ and $\iwinding[(k,\obstuff{K})]$ in 
$\awebstar$ (but not in $\awebv$) can be defined inductively
from $\winding[(\pm 1,\obstuff{K})]$ and $\iwinding[(\pm 1,\obstuff{K})]$. Moreover, the morphisms $\winding[(-1,\obstuff{K})]$ and $\iwinding[(-1,\obstuff{K})]$ in $\awebv$ can be defined from $\iwinding[(1,-1,\obstuff{K},-1)]$ and $\winding[(1,-1,\obstuff{K},-1)]$.
\end{Lemma}

\begin{proof}
The pictures
\begin{gather*}
\begin{tikzpicture}[anchorbase,scale=1]
\draw[usual,rounded corners] (0.5,-1)node[below]{$k$} to (0.5,0) to (-0.5,0);
\draw[usual,rounded corners] (2.5,1)node[above]{$k$} to (2.5,0) to (3.5,0);
\draw[usual] (1.5,-1)node[below]{$\obstuff{K}$} to (1.5,1)node[above]{$\obstuff{K}$};
\draw[affine] (-0.5,-1) to (-0.5,1);
\draw[affine] (3.5,-1) to (3.5,1);
\end{tikzpicture}
=\qfac{k}^{-1}
\cdot
\begin{tikzpicture}[anchorbase,scale=1]
\draw[usual] (0.5,-1)node[below]{$k$} to (0.5,-0.5) to (0,0) to (-0.5,0)node[left,yshift=-0.05cm]{$1$};
\draw[usual,rounded corners] (0.5,-0.5) to (1,0) to (1,0.25) to (-0.5,0.25)node[left,yshift=0.05cm]{$1$};
\draw[usual] (2.5,1)node[above]{$k$} to (2.5,0.5) to (3,0) to (3.5,0)node[right,yshift=0.05cm]{$1$};
\draw[usual,rounded corners] (2.5,0.5) to (2,0) to (2,-0.25) to (3.5,-0.25)node[right,yshift=-0.05cm]{$1$};
\draw[usual] (1.5,-1)node[below]{$\obstuff{K}$} to (1.5,1)node[above]{$\obstuff{K}$};
\draw[affine] (-0.5,-1) to (-0.5,1);
\draw[affine] (3.5,-1) to (3.5,1);
\node at (0.5,0) {$\dots$};
\node at (2.5,0) {$\dots$};
\end{tikzpicture}
\,,
\\
\begin{tikzpicture}[anchorbase,scale=1]
\draw[usual,rounded corners,rdirected=0.06] (0.5,-1)node[below]{$-1$} to (0.5,0) to (-0.5,0);
\draw[usual,rounded corners,directed=0.99] (2.5,1)node[above]{$-1$} to (2.5,0) to (3.5,0);
\draw[usual] (1.5,-1)node[below]{$\obstuff{K}$} to (1.5,1)node[above]{$\obstuff{K}$};
\draw[affine] (-0.5,-1) to (-0.5,1);
\draw[affine] (3.5,-1) to (3.5,1);
\end{tikzpicture}
=
\begin{tikzpicture}[anchorbase,scale=1]
\draw[usual,rdirected=0.06] (1.5,-1)node[below]{$-1$} to (1.5,0.25) to[in=0,out=90] (1,0.75) to[in=90,out=180] (0.5,0.25);
\draw[usual,rounded corners] (0.5,0.25) to (0.5,0) to (-0.25,0);
\draw[usual,directed=0.99] (3.5,1)node[above]{$-1$} to (3.5,-.25) to[in=180,out=270] (4,-0.75) to[in=270,out=0] (4.5,-.25);
\draw[usual,rounded corners] (4.5,-0.25) to (4.5,0) to (5.25,0);
\draw[usual] (2.5,-1)node[below]{$\obstuff{K}$} to (2.5,1)node[above]{$\obstuff{K}$};
\draw[affine] (-0.25,-1) to (-0.25,1);
\draw[affine] (5.25,-1) to (5.25,1);
\end{tikzpicture}
\,,
\end{gather*}
define the morphisms as claimed. 
Their inverses are the $(\placeholder)^{\dinvo}$-duals of these pictures.
\end{proof}

The following lemma compares 
$\awebv$ to the construction in \cite{CaKa-q-satake-sln}:

\begin{Lemma}\label{L:AnnularWebsAffinization}
The category $\awebv$ is equivalent as a diagram category to the affinization $\aff$
(in the sense of {\eg} \cite[Definition 2.1]{MoSa-affinization}) of $\webv$.
\end{Lemma}

\begin{proof}
We only sketch the proof: as often in diagrammatic algebra matching a generator-relation presentation with a ``all diagrams'' definition is lengthy and we omit some details.

First, there is an essentially surjective functor $\Gamma$ from $\awebv$ to $\aff$ that puts a plane web into the annulus. Next, $\aff$
is defined by adjoining more morphisms and relations to $\awebv$, namely 
one coil and its inverse for each $\obstuff{K}$ and relations \cite[equation (2.5)]{MoSa-affinization}.
But using the coils in \autoref{D:AnnularWebsDefinition} one can define 
these more general coils following \cite[equation (2.5), left]{MoSa-affinization} which satisfy \cite[equation (2.5), right]{MoSa-affinization}, showing that $\Gamma$ is full. 
Faithfulness of $\Gamma$ can then be deduced from \autoref{T:EquivalenceMain} below, by showing that the functor therein factors through $\aff$ via $\Gamma$.

Alternatively, one can match the generator-relation presentation of $\awebv$ with the generator-relation presentation from \cite{HaOl-actions-tensor-categories},
and then the topological presentation of $\aff$ with the topological presentation of \cite{HaOl-actions-tensor-categories} and the result follows via \cite[Proposition 10]{HaOl-actions-tensor-categories}. (With the caveat that 
\cite{HaOl-actions-tensor-categories} only discusses tangles and the description therein needs to be extended to webs. That is straightforward but lengthy.)
\end{proof}


\subsection{Quotient by essential circles}\label{SS:AnnularWebsCircles}


We now define quotients of $\awebv$.

\begin{Definition}\label{D:AnnularWebsEssential}
The left and right \emph{essential $k$-circles} are defined to be
\begin{gather*}
\ecirclel=
\begin{tikzpicture}[anchorbase,scale=1]
\draw[usual,directed=0.99] (0.5,0.5)to (0,0.5)node[below]{$k$} to (-0.5,0.5);
\draw[affine] (-0.5,0) to (-0.5,1);
\draw[affine] (0.5,0) to (0.5,1);
\end{tikzpicture}=
\begin{tikzpicture}[anchorbase,scale=1]
\draw[usual,directed=0.55] (0.5,0.5) to[out=90,in=0] (0.25,0.75) to[out=180,in=90] (0,0.5)node[right]{$k$} to[out=270,in=0] (-0.25,0.25) to[out=180,in=270] (-0.5,0.5);
\draw[affine] (-0.5,0) to (-0.5,1);
\draw[affine] (0.5,0) to (0.5,1);
\end{tikzpicture}
\,,\quad
\ecircler=
\begin{tikzpicture}[anchorbase,scale=1]
\draw[usual,directed=0.99] (-0.5,0.5)to (0,0.5)node[below]{$k$} to (0.5,0.5);
\draw[affine] (-0.5,0) to (-0.5,1);
\draw[affine] (0.5,0) to (0.5,1);
\end{tikzpicture}
=
\begin{tikzpicture}[anchorbase,scale=1]
\draw[usual,directed=0.55] (-0.5,0.5) to[out=90,in=180] (-0.25,0.75) to[out=0,in=90] (0,0.5)node[right]{$k$} to[out=270,in=180] (0.25,0.25) to[out=0,in=270] (0.5,0.5);
\draw[affine] (-0.5,0) to (-0.5,1);
\draw[affine] (0.5,0) to (0.5,1);
\end{tikzpicture}
\,.
\end{gather*}
We also say \emph{essential circles} for short.
\end{Definition}

Note that essential circles are nontrivial endomorphism 
of $\munit$. We want to evaluate them. To this end, 
let $\epoly$ denote the $k$th elementary symmetric 
polynomial in $n$ variables, {\ie} $\epoly=\epoly(Z_{1},\dots,Z_{n})$.

\begin{Definition}\label{D:AnnularWebsLeviEval}
The \emph{Levi evaluation} for $\levi$ of the essential circles 
is defined to be 
\begin{gather}\label{Eq:AnnularWebsLeviEval}
\begin{gathered}
\begin{tikzpicture}[anchorbase,scale=1]
\draw[usual,directed=0.99] (0.5,0.5)to (0,0.5)node[below]{$k$} to (-0.5,0.5);
\draw[affine] (-0.5,0) to (-0.5,1);
\draw[affine] (0.5,0) to (0.5,1);
\end{tikzpicture}
-
\epoly[k](\vpar^{-1}\upar[1],\vpar^{-3}\upar[1],\dots,\vpar^{-2l_{1}+1}\upar[1],\dots,\vpar^{-1}\upar[d],\vpar^{-3}\upar[d],\dots,\vpar^{-2l_{d}+1}\upar[d])
\cdot
\begin{tikzpicture}[anchorbase,scale=1]
\draw[affine] (-0.5,0) to (-0.5,1);
\draw[affine] (0.5,0) to (0.5,1);
\end{tikzpicture}
\,,
\\
\begin{tikzpicture}[anchorbase,scale=1]
\draw[usual,directed=0.99] (-0.5,0.5)to (0,0.5)node[below]{$k$} to (0.5,0.5);
\draw[affine] (-0.5,0) to (-0.5,1);
\draw[affine] (0.5,0) to (0.5,1);
\end{tikzpicture}
-
\epoly[k](\vpar\upar[1]^{-1},\vpar^{3}\upar[1]^{-1},\dots,\vpar^{2l_{1}-1}\upar[1]^{-1},\dots,\vpar\upar[d]^{-1},\vpar^{3}\upar[d]^{-1},\dots,\vpar^{2l_{d}-1}\upar[d]^{-1})
\cdot
\begin{tikzpicture}[anchorbase,scale=1]
\draw[affine] (-0.5,0) to (-0.5,1);
\draw[affine] (0.5,0) to (0.5,1);
\end{tikzpicture}
\,,
\end{gathered}
\end{gather}
which are elements in $\End_{\awebv}(\munit)$.
\end{Definition}

\begin{Example}\label{E:AnnularWebsLeviEval}
As an extreme case take $\levi=\gln$. The 
formulas in \autoref{Eq:AnnularWebsLeviEval} then become
\begin{gather*}
\begin{tikzpicture}[anchorbase,scale=1]
	\draw[usual,directed=0.99] (0.5,0.5)to (0,0.5)node[below]{$k$} to (-0.5,0.5);
	\draw[affine] (-0.5,0) to (-0.5,1);
	\draw[affine] (0.5,0) to (0.5,1);
\end{tikzpicture}
-\vpar^{-kn}\qbin{n}{k}\upar[1]^{k}
\cdot
\begin{tikzpicture}[anchorbase,scale=1]
	\draw[affine] (-0.5,0) to (-0.5,1);
	\draw[affine] (0.5,0) to (0.5,1);
\end{tikzpicture}
\,,\quad
\begin{tikzpicture}[anchorbase,scale=1]
	\draw[usual,directed=0.99] (-0.5,0.5)to (0,0.5)node[below]{$k$} to (0.5,0.5);
	\draw[affine] (-0.5,0) to (-0.5,1);
	\draw[affine] (0.5,0) to (0.5,1);
\end{tikzpicture}
-\vpar^{kn}\qbin{n}{k}\upar[1]^{-k}
\cdot
\begin{tikzpicture}[anchorbase,scale=1]
	\draw[affine] (-0.5,0) to (-0.5,1);
	\draw[affine] (0.5,0) to (0.5,1);
\end{tikzpicture}
\,.
\end{gather*}
The appearing scalars are multiples of the categorical dimension of $\uob[k]$, which is the value of the usual circle in the web calculus. Since we want to 
eventually evaluate essential circles to these scalars, this might be 
a hint of a connection to annular webs obtained from evaluation representations 
for the affine Lie algebra, see {\eg} \cite[Section 3]{Qu-annular-skew-howe}.
\end{Example}

The quotient $\lweb$ of $\awebstar$ by an ideal 
$\lideal$ defined later in \autoref{SS:EquivalenceStatement} gives a diagrammatic description of quantum $\levi$-representations.

For $\Kone$ it will turn out that $\lideal$ is the 
two-sided $\vcirc$-$\wprod$-ideal generated by 
the Levi evaluations \autoref{Eq:AnnularWebsLeviEval}.
Being careful with the underlying monoidal structure (details are 
given in \autoref{T:EquivalenceMain}), 
the same holds for $\K$, hence the name.


\section{Representation theory of Levi subalgebras}\label{S:Levi}


This section discusses the representation categories of this paper. 
The below is (partially) well-known and we will be brief whenever 
appropriate. A lot of details and also background 
can be found in texts such as \cite{Ja-lectures-qgroups}.


\subsection{The general linear representation category}\label{SS:LeviGln}


We start with notations regarding the general linear quantum algebra.
Let $\Uglnv$ be the \emph{divided power quantum enveloping algebra for $\gln$}, where we use the conventions, excluding the Hopf algebra structure, from \cite{AnPoWe-representation-qalgebras} 
in the special case of $\gln$ 
(using $K_{i}^{\pm 1}=L_{i}^{\pm 1}L_{i+1}^{\mp 1}$).
The algebra $\Uglnv$ specializes to either the $\K$-algebra $\Ugln$,
for which we now recall the relevant formulas, 
or the $\Kone$-algebra $\Uglnone$.

The algebra $\Ugln$ is generated by $L_{i}^{\pm 1}$ for $i\in\set{1,\dots,n}$ 
(these are inverses) 
and $E_{i}$, $F_{i}$ for $i\in\set{1,\dots,n-1}$ and the Hopf algebra 
structure used in this paper is
\begin{alignat*}{3}
\Delta(E_{i})&=E_{i}\otimes L_{i}L_{i+1}^{-1}+1\otimes E_{i}
,\quad
\hspace{0.15cm}\epsilon(E_{i})&&=0
,\quad
S(E_{i})&&=-E_{i}L_{i}^{-1}L_{i+1},
\\
\Delta(F_{i})&=F_{i}\otimes 1+L_{i}^{-1}L_{i+1}\otimes F_{i}
,\quad
\epsilon(F_{i})&&=0
,\quad
S(F_{i})&&=-L_{i}L_{i+1}^{-1}F_{i},
\end{alignat*}
with $L_{i}^{\pm 1}$ being group like.

The \emph{vector representation} $\vecrep=\vecrep(\gln)=\K\set{v_{1},\dots,v_{n}}$ of $\Ugln$ is the given $\K$-vector space with action
\begin{gather*}
L_{i}^{\pm 1}\acts v_{j}=\qpar^{\pm\delta_{i,j}}v_{j}
,\quad
E_{i}\acts v_{j}=\delta_{i,j-1}v_{j-1}
,\quad
F_{i}\acts v_{j}=\delta_{i,j}v_{j+1}
.
\end{gather*}
Let $T\vecrep$ be the tensor algebra. The \emph{$k$th quantum exterior power} 
$\ext{\qpar}{k}\vecrep$ is defined as the degree $k$ part (in the usual sense) of the \emph{quantum exterior algebra} given by 
\begin{gather}\label{Eq:LeviExtAlgebra}
\ext{\qpar}{\bullet}\vecrep
=
\bigoplus_{k\in\N}\ext{\qpar}{k}\vecrep
=T\vecrep/\langle v_{h}\otimes v_{h},v_{j}\otimes v_{i}+\qpar^{-1}v_{i}\otimes v_{j}|i<j\rangle_{\text{two-sided $\otimes$-ideal}}.
\end{gather}
The exterior powers are $\Ugln$-representations by using the Hopf algebra 
structure, and so is $\munit=\K$ itself and all the duals of the above, denoted 
by using negative powers: $\ext{\qpar}{k}\vecrep=(\ext{\qpar}{-k}\vecrep)^{\pivo}$ for $k\in\Z_{<0}$.

\begin{Lemma}\label{L:LeviExtPower}
For $k\in\N$,
the $\Kstar$-vector space $\ext{\qpar}{k}\vecrep[\qpar]$ has a basis given by $\set[\big]{v_{S}=v_{i_{1}}\otimes\dots\otimes v_{i_{k}}|S=(i_{1}<\dots<i_{k})\text{ for }i_{j}\in\set{1,\dots,n}}$. If $-k\in\N$, then the $\Kstar$-vector space $\ext{\qpar}{k}\vecrep[\qpar]$ has a basis given by $\set[\big]{v_{S}^{\pivo}=v_{i_{-k}}\otimes\dots\otimes v_{i_{1}}|S=(i_{1}<\dots<i_{-k})\text{ for }i_{j}\in\set{1,\dots,n}}$.
\end{Lemma}

\begin{proof}
Easy and omitted.
\end{proof}

\begin{Notation}\label{N:LeviExtPower}
We also use the notation $v_{S}$ from \autoref{L:LeviExtPower} more generally 
for any $S=(i_{1},\dots,i_{k})$ for $i_{j}\in\set{1,\dots,n}$, and use the 
usual set operations on them. 
Recall that such expressions need to potentially be reordered 
using \autoref{Eq:LeviExtAlgebra} to 
match the basis of \autoref{L:LeviExtPower}.
\end{Notation}

We now consider so-called 
$\Ugln$-representations of type $1$ which, as usual, is not a serious 
restriction, see {\eg} \cite[Section 5.2]{Ja-lectures-qgroups} 
for details.

\begin{Definition}\label{D:LeviGlnRep}
Let $\rep[\gln]$ denote the category of finite dimensional $\Ugln$-representations
of type $1$. We view $\rep[\gln]$ as pivotal using the above Hopf algebra
structure on $\Ugln$. Let further $\fund[\gln]$ 
denote the full pivotal subcategory 
with objects of the form $\ext{\qpar}{\obstuff{K}}\vecrep=\ext{\qpar}{k_{1}}\vecrep\otimes\dots\otimes\ext{\qpar}{k_{m}}\vecrep$ 
for $\obstuff{K}=(k_{1},\dots,k_{m})\in\Z^{m}$ and $m\in\N$.
\end{Definition}

We call $\rep[\gln]$ the \emph{representation category} of $\Ugln$ and 
$\fund[\gln]$ its \emph{fundamental category}. (We use the same terminology for $\levi$ defined below.) The following is crucial, 
but well-known, and will be used throughout. 
To state it let $X^{+}_{\gln}\subset\Z^{n}$ denote the set of dominant 
integral $\gln$-weight, {\ie} tuples $\lambda=(\lambda_{1},\dots,\lambda_{n})$ 
such that $\lambda_{1}\geq\dots\geq\lambda_{n}$.

\begin{Lemma}\label{L:LeviGlnRepSemi}
We have the following.
\begin{enumerate}

\item The category $\repstar[\gln]$ is semisimple, its simple 
objects can be indexed by $\lambda\in X^{+}_{\gln}$
and their characters are given by Weyl's character formula.

\item $\fundstar[\gln]$ is 
pivotally equivalent to $\repstar[\gln]$
upon additive idempotent completion.

\end{enumerate}
\end{Lemma}

We will denote the simple objects in \autoref{L:LeviGlnRepSemi}.(a) by $\simple[\lambda]$ for $\lambda\in X^{+}_{\gln}$. We do not need them explicitly, but their construction is well-known (the $\simple[\lambda]$ are often called Weyl modules).

\begin{proof}
Claim (a). See {\eg} \cite[Theorems 5.15 and 5.17]{Ja-lectures-qgroups} or \cite[Section 6]{AnPoWe-representation-qalgebras}.

Claim (b). By (a) classical theory applies, see {\eg} \cite[Theorems 5.15 and 5.17]{Ja-lectures-qgroups} or \cite[Section 6]{AnPoWe-representation-qalgebras} for details.
\end{proof}

Let us now list generating $\Ugln$-equivariant morphisms that will 
be the images of the generators of $\web$. 
The notation is hopefully suggestive.

For tuples $S,T$ as in \autoref{N:LeviExtPower} let $|S<T|=\big|\set{(s,t)\in S\times T|s<t}\big|$ and $|S,N|=|S<N|-|N<S|$ for $N=(1,2,\dots,n)$. Using such tuples and this notation we define (here $k,l\in\N$):
\begin{gather}\label{Eq:LeviGens1}
\begin{aligned}
\mergemap{k,l}{k{+}l}&\colon\ext{\qpar}{k}\vecrep\otimes\ext{\qpar}{l}\vecrep\to
\ext{\qpar}{k{+}l}\vecrep
,\quad
v_{S}\otimes v_{T}
\mapsto\delta_{S\cap T,\emptyset}(-\qpar)^{-|T<S|}v_{S\cup T},
\\
\splitmap{k,l}{k{+}l}&\colon\ext{\qpar}{k{+}l}\vecrep\to
\ext{\qpar}{k}\vecrep\otimes\ext{\qpar}{l}\vecrep
,\quad
v_{U}\mapsto(-1)^{kl}\sum_{S\sqcup T=U,|S|=k}(-\qpar)^{|S<T|}
v_{S}\otimes v_{T},
\\
\capl&\colon\ext{\qpar}{{-}k}\vecrep\otimes\ext{\qpar}{k}\vecrep\to\K
,\quad
v_{S}^{\pivo}\otimes v_{T}\mapsto\delta_{S,T}
,\\
\capr&\colon\ext{\qpar}{k}\vecrep\otimes\ext{\qpar}{{-}k}\vecrep\to\K
,\quad
v_{S}\otimes v_{T}^{\pivo}\mapsto\qpar^{|S,N|}\delta_{S,T}
,
\\
\cupl&\colon\K\to\ext{\qpar}{k}\vecrep\otimes\ext{\qpar}{{-}k}\vecrep
,\quad
1\mapsto\sum_{|S|=k}v_{S}\otimes v_{S}^{\pivo}
,
\\
\cupr&\colon\K\to\ext{\qpar}{{-}k}\vecrep\otimes\ext{\qpar}{k}\vecrep
,\quad
1\mapsto\sum_{|S|=k}\qpar^{-|S,N|}v_{S}^{\pivo}\otimes v_{S}
.
\end{aligned}
\end{gather}

\begin{Lemma}\label{L:LeviGlnMaps}
The morphisms in \autoref{Eq:LeviGens1} are $\vcirc$-$\otimes$-generators of $\fundstar[\gln]$.
\end{Lemma}

\begin{proof}
A careful check of the relations shows that these 
maps are $\Uglnstar$-equivariant. That they generate follows from 
\autoref{L:LeviGlnRepSemi} and classical theory.
\end{proof}

We have an \emph{algebraic version of explosion}:

\begin{Lemma}\label{L:LeviExplosion}
\emph{Explosion} holds for $\fundstar[\gln]$, that is
\begin{gather*}
\idmor_{k{+}l}
=\qbin{k+l}{k}^{-1}\mergemap{k,l}{k{+}l}\splitmap{k,l}{k{+}l}.
\end{gather*}
\end{Lemma}

\begin{proof}
A direct computation.
\end{proof}

We denote by $\expl$ and $\iexpl$ the successive explosion of $k$ strands.

We also have the following \emph{$(1,1)$-overcrossings} and \emph{$(1,1)$-undercrossings}:
\begin{align*}
\braid[{1,1}]&=\qpar\idmor_{1,1}-\splitmap{{1,1}}{2}\mergemap{{1,1}}{2}
\colon
\vecrep\otimes\vecrep\to\vecrep\otimes\vecrep
,\quad
v_{i}\otimes v_{j}\mapsto
\begin{cases}
\qpar v_{i}\otimes v_{i}&\text{if }i=j,
\\
v_{j}\otimes v_{i}&\text{if }i<j,
\\
v_{j}\otimes v_{i}+(\qpar-\qpar^{-1})v_{i}\otimes v_{j}&\text{if }i>j,
\end{cases}
\\
\braid[{1,1}]^{-1}&=\qpar^{-1}\idmor_{1,1}-\splitmap{{1,1}}{2}\mergemap{{1,1}}{2}
\colon
\vecrep\otimes\vecrep\to\vecrep\otimes\vecrep
,\quad
v_{i}\otimes v_{j}\mapsto
\begin{cases}
\qpar^{-1}v_{i}\otimes v_{i}&\text{if }i=j,
\\
v_{j}\otimes v_{i}&\text{if }i>j,
\\
v_{j}\otimes v_{i}+(\qpar^{-1}-\qpar)v_{i}\otimes v_{j}&\text{if }i<j.
\end{cases}
\end{align*}
We also get crossings $\braid[{k,l}]^{\pm 1}\colon\ext{\qpar}{k}\vecrep\otimes\ext{\qpar}{l}\vecrep\to\ext{\qpar}{l}\vecrep\otimes\ext{\qpar}{k}\vecrep$ for all $k,l\in\Z$ by mimicking \autoref{Eq:PrelimCrossing} and mating.

\begin{Lemma}\label{L:LeviBraid}
The crossings satisfy the \emph{Reidemeister II and III relations} and 
various \emph{naturality relations}, and can be alternatively defined by explosion.
\end{Lemma}

\begin{proof}
Well-known and omitted (for the statement about the alternative 
definition using explosion see \autoref{L:WebsGlnProperties} 
imported via \autoref{T:EquivalenceGln}).
\end{proof}


\subsection{The Levi representation category}\label{SS:LeviRepCat}


Recall that we have fixed 
$\levi=\gln[l_{1}]\oplus\dots\oplus\gln[l_{d}]$ of which we think as being
\begin{gather*}
\levi
=
\left(
\begin{array}{ccc}
\gln[l_{1}] & 0 & 0 \\
0 & \ddots & 0 \\
0 & 0 & \gln[l_{d}]
\end{array}
\right)\subset\gln
,\quad
\text{generators in $\gln[l_{i}]$}\colon
\left\{
\begin{array}{cccc}
L_{i,1} & E_{i,1} & \phantom{0} & \phantom{0} \\
F_{i,1} & \ddots & \ddots & \phantom{0} \\
\phantom{0} & \ddots & \ddots & E_{i,l_{i}-1} \\
\phantom{0} & \phantom{0} & F_{i,l_{i}-1} & L_{i,l_{i}}
\end{array}
\right.
\;,
\end{gather*}
where we reindex the elements $L_{i}^{\pm}$, $E_{i}$ and $F_{i}$ as indicated. 
(Note that all $L_{i}^{\pm}$ appear in this reindexing, but not all $E_{i}$ and $F_{i}$.)

\begin{Definition}\label{D:LeviAlgebra}
Let $\Ulv$ be the $\Zv$-subalgebra of $\Uglnv$ generated 
by these $L_{i,k}^{\pm}$, $E_{i,k}$ and $F_{i,k}$. We endow 
$\Ulv$ with the structure of a Hopf algebra by restricting 
the one for $\Uglnv$.
\end{Definition}

The following lemma gives us a 
\emph{block decomposition} and will be used without further reference:

\begin{Lemma}\label{L:LeviDecomposition}
We have $\Ulv\cong\Uglnv[{\gln[l_{1}]}]\otimes\dots\otimes\Uglnv[{\gln[l_{d}]}]$ as $\Zv$-algebras.
\end{Lemma}

\begin{proof}
By definition.
\end{proof}

The representation theory of 
$\Ul$ is easy (knowing the representation theory for 
$\Ugln[{\gln[m]}]$),
but nevertheless we state a few lemmas that we will use.
For starters, note that all $\Ugln$-representations restrict to $\Ul$-representations. Note also that 
the vector representation $\vecrep(\gln[l_{i}])$ of $\Ugln[{\gln[l_{i}]}]$ 
is a $\Ul$-representation with action inflated to $\Ul$. The same 
holds for the exterior powers.

\begin{Lemma}\label{L:LeviModules}
As $\Ulstar$-representations we have $\ext{\qpar}{k}\vecrep{\cong}\ext{\qpar}{k}\bigoplus_{i=1}^{d}\vecrep[\qpar](\gln[l_{i}]){\cong}\bigoplus_{k_{1}+\dots+k_{d}=k}\ext{\qpar}{k_{1}}\vecrep[\qpar](\gln[l_{1}])$ $\otimes\dots\otimes\ext{\qpar}{k_{d}}\vecrep[\qpar](\gln[l_{d}])$. We also have 
$\vecrep[\qpar]^{\otimes k}\cong\bigoplus_{k_{1}+\dots+k_{d}=k}\vecrep[\qpar](\gln[l_{1}])^{\otimes k_{1}}\otimes\dots\otimes\vecrep[\qpar](\gln[l_{d}])^{\otimes k_{d}}$.
\end{Lemma}

Note that \autoref{L:LeviModules} implies that $\vecrep[\qpar]$ is not simple as a $\Ugln[\levi]$-representation.

\begin{proof}
The case $k=1$ is clear by, for example, using the usual 
diagrammatic description of $\vecrep[\qpar]$ (the actions of the $L_{i}^{\pm 1}$ are omitted in the following illustration):
\begin{gather*}
\begin{tikzcd}[ampersand replacement=\&,column sep=1cm,row sep=1cm]
\colorbox{orchid!35}{\mystrut$v_{8}$}
\ar[r,"E_{3,1}",yshift=0.075cm,orchid]
\&
\colorbox{orchid!35}{\mystrut$v_{7}$}
\ar[r,"E_{6}",yshift=0.075cm,gray,densely dotted]
\ar[l,"F_{3,1}",yshift=-0.075cm,orchid]
\&
\colorbox{tomato!35}{\mystrut$v_{6}$}
\ar[r,"E_{2,2}",yshift=0.075cm,tomato]
\ar[l,"F_{6}",yshift=-0.075cm,gray,densely dotted]
\&
\colorbox{tomato!35}{\mystrut$v_{5}$}
\ar[r,"E_{2,1}",yshift=0.075cm,tomato]
\ar[l,"F_{2,2}",yshift=-0.075cm,tomato]
\&
\colorbox{tomato!35}{\mystrut$v_{4}$}
\ar[r,"E_{3}",yshift=0.075cm,gray,densely dotted]
\ar[l,"F_{2,1}",yshift=-0.075cm,tomato]
\&
\colorbox{spinach!35}{\mystrut$v_{3}$}
\ar[r,"E_{1,2}",yshift=0.075cm,spinach]
\ar[l,"F_{3}",yshift=-0.075cm,gray,densely dotted]
\&
\colorbox{spinach!35}{\mystrut$v_{2}$}
\ar[r,"E_{1,1}",yshift=0.075cm,spinach]
\ar[l,"F_{1,2}",yshift=-0.075cm,spinach]
\&
\colorbox{spinach!35}{\mystrut$v_{1}$}
\ar[l,"F_{1,1}",yshift=-0.075cm,spinach]
\end{tikzcd}
.
\end{gather*}
This illustrates the case $l_{1}=3$, $l_{2}=3$, $l_{3}=2$ 
and $n=8$. The case of general $k$ and the second isomorphism are similar and omitted.
\end{proof}

We now reindex the basis of $\vecrep$ to 
$\set{v_{1,1},\dots,v_{1,l_{1}},\dots,v_{d,1},\dots,v_{d,l_{d}}}$, 
which then induces a reindexing of the basis of $\ext{\qpar}{k}\vecrep$ 
that we will use below.

\begin{Lemma}\label{L:LeviModules2}
As $\Ulstar$-representations we have $\vecrep[\qpar]^{\otimes k}\cong\ext{\qpar}{k}\vecrep\oplus W$ and no simple constituent of $\ext{\qpar}{k}\vecrep$ appears in $W$.
\end{Lemma}

\begin{proof}
By looking at highest weight vectors and classical theory,
this follows directly from the first and the second decomposition in \autoref{L:LeviModules}.
\end{proof}

\begin{Definition}\label{D:LeviRep}
Let $\rep[\levi]$ denote the category of finite dimensional $\Ul$-representations
of type $1$. We view $\rep[\levi]$ as pivotal using the above Hopf algebra
structure on $\Ugln$. Let further $\fund[\levi]$ 
denote the full pivotal subcategory 
with objects of the form $\ext{\qpar}{\obstuff{K}}\vecrep$ 
for $\obstuff{K}=(k_{1},\dots,k_{m})\in\Z^{m}$ and $m\in\N$.
\end{Definition}

We write $\lprod$ for the monoidal structure on $\rep[\levi]$ 
and $\fund[\levi]$ on the morphism level, and $\hcirc$ on the object level.
With contrast to \autoref{Eq:AnnularWebsPropsMonoidal}, 
a picture for the monoidal structure on $\fund[\levi]$ is
\begin{gather}\label{Eq:LeviPropsMonoidal}
\begin{tikzpicture}[anchorbase,scale=1,yscale=1]
\draw[spinach!35,fill=spinach!35] (-0.5,-0.25) to (-0.5,0.25) to (1.5,0.25) to (1.5,-0.25) to (-0.5,-0.25);
\draw[usual] (0,-1) to (0,-0.25);
\draw[usual] (0,0.25) to (0,1);
\draw[usual] (1,-1) to (1,-0.25);
\draw[usual] (1,0.25) to (1,1);
\node at (0.5,-0.625) {$\dots$};
\node at (0.5,0.625) {$\dots$};
\node at (0.5,0) {$\morstuff{f}$};
\end{tikzpicture}	
\lprod
\begin{tikzpicture}[anchorbase,scale=1,yscale=1]
\draw[spinach!35,fill=spinach!35] (-0.5,-0.25) to (-0.5,0.25) to (1.5,0.25) to (1.5,-0.25) to (-0.5,-0.25);
\draw[usual] (0,-1) to (0,-0.25);
\draw[usual] (0,0.25) to (0,1);
\draw[usual] (1,-1) to (1,-0.25);
\draw[usual] (1,0.25) to (1,1);
\node at (0.5,-0.625) {$\dots$};
\node at (0.5,0.625) {$\dots$};
\node at (0.5,0) {$\morstuff{g}$};
\end{tikzpicture}
=
\begin{tikzpicture}[anchorbase,scale=1,yscale=1]
\draw[spinach!35,fill=spinach!35] (1.5,0.25) to (3.5,0.25) to (3.5,0.75) to (1.5,0.75) to (1.5,0.25);
\draw[spinach!35,fill=spinach!35] (-0.5,-0.25) to (1.5,-0.25) to (1.5,-0.75) to (-0.5,-0.75) to (-0.5,-0.25);
\draw[usual] (0,-1) to (0,-0.75);
\draw[usual,crossline] (0,-0.25) to (0,1);
\draw[usual] (1,-1) to (1,-0.75);
\draw[usual,crossline] (1,-0.25) to (1,1);
\draw[usual] (2,-1) to (2,0.25);
\draw[usual] (2,0.75) to (2,1);
\draw[usual] (3,-1) to (3,0.25);
\draw[usual] (3,0.75) to (3,1);
\node at (0.5,-1) {$\dots$};
\node at (0.5,1) {$\dots$};
\node at (2.5,-1) {$\dots$};
\node at (2.5,1) {$\dots$};
\node at (0.5,-0.5) {$\morstuff{f}$};
\node at (2.5,0.5) {$\morstuff{g}$};
\end{tikzpicture}
.
\end{gather}

\begin{Remark}\label{R:LeviMonoidal} 
\autoref{Eq:LeviPropsMonoidal} is a standard construction in 
the theory of Hopf algebras, see {\eg} \cite{TuVi-monoidal-tqft}.
\end{Remark}

\begin{Lemma}\label{L:LeviRepSemi}
We have the following.
\begin{enumerate}

\item The category $\repstar[\levi]$ is semisimple, and its simple 
objects are of the form $\simple[\lambda_{1}]\otimes\dots\otimes\simple[\lambda_{d}]$ 
with the factors being simple objects of $\repstar[{\gln[l_{i}]}]$.

\item The additive idempotent completion of $\fundstar[\levi]$ is 
pivotally equivalent to $\repstar[\levi]$.

\end{enumerate}
\end{Lemma}

\begin{proof}
\autoref{L:LeviGlnRepSemi} applies componentwise.
\end{proof}

We define the \emph{Levi $(1,1)$-overcrossings} and \emph{Levi $(1,1)$-undercrossing} to be
\begin{gather*}
\lbraid[{1,1}]\colon\vecrep\otimes\vecrep\to\vecrep\otimes\vecrep
,\quad
v_{i,j}\otimes v_{k,l}\mapsto
\begin{cases}
\braid[{1,1}](v_{i,j}\otimes v_{i,l})&\text{if }i=k,
\\
v_{k,l}\otimes v_{i,j}&\text{else},
\end{cases}
\\
(\lbraid[{1,1}])^{-1}\colon\vecrep\otimes\vecrep\to\vecrep\otimes\vecrep
,\quad
v_{i,j}\otimes v_{k,l}\mapsto
\begin{cases}
(\braid[{1,1}])^{-1}(v_{i,j}\otimes v_{i,l})&\text{if }i=k,
\\
v_{k,l}\otimes v_{i,j}&\text{else}.
\end{cases}
\end{gather*}
In other words, $\lbraid[{1,1}]$ is the respective $\gln[l_{i}]$ braiding 
within one block, and the swap map otherwise, and similarly 
for its inverse. Note that these maps are in general not $\Ugln$-equivariant.

\begin{Example}\label{E:LeviBraiding}
In the extreme case that $\levi=\cartan$ the 
Levi $(1,1)$-overcrossings equals the Levi $(1,1)$-undercrossings
equals the swap map.
\end{Example}

\begin{Definition}\label{D:LeviBraiding}
For $k,l\in\N$ define $\lbraid[{k,l}]$ as the composition
\begin{gather*}
\ext{\qpar}{k}\vecrep\otimes\ext{\qpar}{l}\vecrep
\xrightarrow{\morstuff{i}}
\vecrep^{\otimes (k+l)}
\xrightarrow{\morstuff{x}}
\vecrep^{\otimes (k+l)}
\xrightarrow{\morstuff{p}}
\ext{\qpar}{l}\vecrep\otimes\ext{\qpar}{k}\vecrep
,
\end{gather*}
where are $\morstuff{i}$ and $\morstuff{p}$ are inclusion and projection, respectively, 
and $\morstuff{x}$ is defined as in \autoref{Eq:WebsExplosion} but with Levi crossings.
Define $(\lbraid[{k,l}])^{-1}$ similarly.
\end{Definition}

We also get various mates of which we think as rotated versions of 
the ones in \autoref{D:LeviBraiding}.

\begin{Lemma}\label{L:LeviBraid2}
The Levi crossings are $\Ulstar$-equivariant, 
satisfy the \emph{Reidemeister II and III relations} and 
various \emph{naturality relations}.
\end{Lemma}

\begin{proof}
This follows from a calculation, \autoref{L:LeviDecomposition} and \autoref{L:LeviBraid}.
\end{proof}

Define now the following \emph{coil maps}:
\begin{gather}\label{Eq:LeviGens2}
\begin{gathered}
\wmapl\colon
\vecrep\otimes\ext{\qpar}{\obstuff{K}}\vecrep\to\ext{\qpar}{\obstuff{K}}\vecrep\otimes\vecrep
,\quad
v_{i,j}\otimes w\mapsto\upar[i]\lbraid[{1,\obstuff{K}}](v_{i,j}\otimes w)
,
\\
\wmapr\colon
\ext{\qpar}{\obstuff{K}}\vecrep\otimes\vecrep\to\vecrep\otimes\ext{\qpar}{\obstuff{K}}\vecrep
,\quad
v\otimes w\mapsto\upar[i]^{-1}(\lbraid[{1,\obstuff{K}}])^{-1}(w\otimes v_{i,j})
.
\end{gathered}
\end{gather}
Note that $\wmapl$ and $\wmapr$ are inverses. The coil maps are not $\Ugln$-equivariant in general.

\begin{Lemma}\label{L:LeviMaps}
The morphisms in \autoref{Eq:LeviGens1} and 
\autoref{Eq:LeviGens2} are $\vcirc$-$\lprod$-generators of $\fundstar[\levi]$.
\end{Lemma}

\begin{proof}
The morphisms in \autoref{Eq:LeviGens1} are $\Uglnstar$-equivariant, 
so they are also $\Ulstar$-equivariant, and one easily checks 
that the morphisms in \autoref{Eq:LeviGens2} are $\Ulstar$-equivariant.
That these generate follows from \autoref{L:LeviGlnMaps} and
\autoref{L:LeviModules}.
\end{proof}

For $\obstuff{K}=\emptyset$ we will write $\wmapl[1]=\wmapl[1,\emptyset]$.

\begin{Definition}\label{D:LeviProjectors}
For $i\in\set{1,\dots,d}$ define 
\begin{gather*}
\lproj=
\prod_{j\neq i}
\frac{\wmapl[1]-\upar[j]}{\upar[i]-\upar[j]}
\in\End_{\fund[\levi]}(\vecrep),
\end{gather*}
which we call \emph{Levi block projectors}.
\end{Definition}

\begin{Lemma}\label{L:LeviProjectors}
We have $\lproj\lproj[j]=\delta_{i,j}\lproj$ and 
$\idmor_{\vecrep[\qpar]}=\sum_{i=1}^{d}\lproj$. 
These projectors realize the decomposition
$\vecrep[\qpar]\cong\bigoplus_{i=1}^{d}\vecrep[\qpar](\gln[l_{i}])$.
\end{Lemma}

\begin{proof}
Note that $\wmapl[1]$ is given by multiplication 
by $\upar[i]$ on $v_{i,j}$. Thus, the formula 
for $\lproj$ is the usual Lagrange-type 
interpolation and the claims follow.
\end{proof}


\subsection{Levi crossings}\label{SS:LeviCrossings}


The Levi crossings are not $\Ugln$-equivariant
in general, and there is no planar web picture for it. However, 
it will be helpful to have the 
following diagrammatic notation.
For the Levi $(k,l)$-overcrossings and the Levi $(k,l)$-undercrossings we use:
\begin{gather*}
\text{over}:
\lbraid[{k,l}]
\leftrightsquigarrow
\begin{tikzpicture}[anchorbase,scale=1]
\draw[usual,directed=0.99] (1,0)node[below]{$l$} to (0,1)node[above]{$l$};
\draw[usual,directed=0.99,levia=0.5,crossline] (0,0)node[below]{$k$} to (1,1)node[above]{$k$};
\end{tikzpicture}
,\quad
\text{under}:
(\lbraid[{k,l}])^{-1}
\leftrightsquigarrow
\begin{tikzpicture}[anchorbase,scale=1]
\draw[usual,directed=0.99] (0,0)node[below]{$k$} to (1,1)node[above]{$k$};
\draw[usual,directed=0.99,levia=0.5,crossline] (1,0)node[below]{$l$} to (0,1)node[above]{$l$};
\end{tikzpicture}
.
\end{gather*}
We also use rotated pictures for their mates.

By \autoref{L:LeviBraid2}, we have 
the \emph{Reidemeister II and III relations}, {\eg}
\begin{gather*}
\begin{tikzpicture}[anchorbase,scale=1]
\draw[usual,rounded corners] (1,0) to (0,1) to (1,2);
\draw[usual,levia=0.25,levia=0.75,crossline,rounded corners] (0,0) to (1,1) to (0,2);
\end{tikzpicture}
=
\begin{tikzpicture}[anchorbase,scale=1]
\draw[usual] (0,0) to (0,2);
\draw[usual] (1,0) to (1,2);
\end{tikzpicture}
\,,\quad
\begin{tikzpicture}[anchorbase,scale=1]
\draw[usual] (2,0) to (0,2);
\draw[usual,levi=0.75,crossline,rounded corners] (1,0) to (0,1) to (1,2);
\draw[usual,levia,levia=0.25,crossline] (0,0) to (2,2);
\end{tikzpicture}
=
\begin{tikzpicture}[anchorbase,scale=1]
\draw[usual] (2,0) to (0,2);
\draw[usual,crossline,levi=0.25,rounded corners] (1,0) to (2,1) to (1,2);
\draw[usual,crossline,levia,levia=0.75] (0,0) to (2,2);
\end{tikzpicture}
.
\end{gather*}

We can use this to define \emph{Levi braids} associated to any braid word. 
Of particular importance will be the (positive) \emph{Levi full twist} 
on $k$ strands (denoted by a box notation). By definition, this map is the square of the positive 
lift, using Levi overcrossings, of the longest word in the symmetric group on $\set{1,\dots,k}$. For example, for $k=4$ this full twist is
\begin{gather*}
\begin{tikzpicture}[anchorbase,scale=1]
\draw[tomato!50,fill=tomato!50] (-0.25,0) rectangle (1.75,0.5);
\draw[usual] (0,-0.5) to (0,0);
\draw[usual] (0,0.5) to (0,1);
\draw[usual] (0.5,-0.5) to (0.5,0);
\draw[usual] (0.5,0.5) to (0.5,1);
\draw[usual] (1,-0.5) to (1,0);
\draw[usual] (1,0.5) to (1,1);
\draw[usual] (1.5,-0.5) to (1.5,0);
\draw[usual] (1.5,0.5) to (1.5,1);
\node at (0.75,0.25) {$\fulltwist$};
\end{tikzpicture}
=
\left(
\scalebox{0.75}{$\begin{tikzpicture}[anchorbase,scale=1]
\draw[usual] (3,0) to (0,3);
\draw[usual,crossline,levi=0.83,rounded corners] (2,0) to (0,2) to (1,3);
\draw[usual,crossline,levia=0.49,levia=0.67,rounded corners] (1,0) to (0,1) to (2,3);
\draw[usual,levib=0.16,levib=0.33,levib=0.5,crossline] (0,0) to (3,3);
\end{tikzpicture}$}
\right)^{2}
.
\end{gather*}
We also have the usual naturality relations such as
\begin{gather*}
\begin{tikzpicture}[anchorbase,scale=1]
\draw[usual] (1,-1) to (-1,0) to (-1,1);
\draw[usual,crossline,levi=0.32] (0,0) to (-1,-1);
\draw[usual,crossline,levi=0.64] (1,0) to (0,-1);
\draw[usual] (0,-0.02) to(0,0) to[out=90,in=180] (0.5,0.5) to[out=0,in=90] (1,0) to (1,-0.02);
\end{tikzpicture}
=
\begin{tikzpicture}[anchorbase,scale=1]
\draw[usual] (0,0) to[out=90,in=180] (0.5,0.5) to[out=0,in=90] (1,0);
\draw[usual] (2,0) to (0,2);
\end{tikzpicture}
\quad\text{or}\quad
\begin{tikzpicture}[anchorbase,scale=1]
\draw[usual] (1,-1) to (-1,0) to (-1,1);
\draw[usual,crossline,levi=0.32] (0,0) to (-1,-1);
\draw[usual,crossline,levi=0.64] (1,0) to (0,-1);
\draw[usual] (0.5,0.5) to (0,0);
\draw[usual] (0.5,0.5) to (1,0) to (1,-0.03);
\draw[usual] (0.5,0.5) to (0.5,1);
\end{tikzpicture}
=
\begin{tikzpicture}[anchorbase,scale=1]
\draw[usual] (2,0) to (0,2);
\draw[usual,crossline,levi=0.67] (2,2) to (0.5,0.5);
\draw[usual] (0.5,0.5) to (0,0);
\draw[usual] (0.5,0.5) to (1,0);
\end{tikzpicture}
,
\end{gather*}
including the various $(\placeholder)^{\danti}$ and $(\placeholder)^{\dinvo}$-duals.

However, we need to be careful with the \emph{Reidemeister I 
relation} as $\ext{\qpar}{k}\vecrep$ needs not to be 
simple as a $\Ulstar$-representation, see \autoref{L:LeviModules}. Nevertheless, we still have the following lemma.

\begin{Lemma}\label{L:LeviReidemeisterOne}
The Levi crossings are diagonal matrices in the basis 
given by the decomposition of 
$\ext{\qpar}{k}\vecrep[\qpar]$ from \autoref{L:LeviModules}.
\end{Lemma}

\begin{proof}
Note that the decomposition of 
$\ext{\qpar}{k}\vecrep[\qpar]$ 
from \autoref{L:LeviModules} is multiplicity free, 
and Schur's lemma applies. (We stress that Schur's lemma in this setting does not need the underlying field to be algebraically closed, see {\eg} 
\cite[Corollary 7.4]{AnPoWe-representation-qalgebras} or \cite[Remark 2.29]{AnTu-tilting}.)
\end{proof}

One can check that the diagonal entries mentioned in 
\autoref{L:LeviReidemeisterOne} are given by products of 
the Reidemeister I scalars in \autoref{Eq:WebsRM1}.


\section{The equivalence}\label{S:Equivalence}


We now state and prove our main result.


\subsection{A reminder on the \texorpdfstring{$\gln$}{gln} story}\label{SS:EquivalenceReminder}


We first recall the relationship between $\gln$ webs 
and the representation theory of $\Ugln$.
Define a functor 
\begin{gather*}
\pfunctor=\pfunctor^{ext}(\gln)\colon\web\to\fund[\gln]
\end{gather*}
sending the object $\obstuff{K}$ to 
$\ext{\qpar}{\obstuff{K}}\vecrep$ and the 
generating morphisms of $\web$ to the following maps:
\begin{gather*}
\scalebox{0.95}{$\begin{tikzpicture}[anchorbase,scale=1,yscale=-1]
\draw[usual,directed=0.99] (0.5,0.5) to (0,0)node[above]{$k$};
\draw[usual,directed=0.99] (0.5,0.5) to (1,0)node[above]{$l$};
\draw[usual] (0.5,1)node[below]{$k{+}l$} to (0.5,0.5);
\end{tikzpicture}
\hspace*{-0.1cm}\mapsto
\splitmap{k,l}{k{+}l}
,
\begin{tikzpicture}[anchorbase,scale=1]
\draw[usual] (0,0)node[below]{$k$} to (0.5,0.5);
\draw[usual] (1,0)node[below]{$l$} to (0.5,0.5);
\draw[usual,directed=0.99] (0.5,0.5) to (0.5,1)node[above]{$k{+}l$};
\end{tikzpicture}
\hspace*{-0.1cm}\mapsto
\mergemap{k,l}{k{+}l}
,
\begin{tikzpicture}[anchorbase,scale=1]
\draw[usual,directed=0.99] (0,0)node[below]{$k$} to[out=90,in=180] (0.5,0.5)node[above]{\phantom{k}} to[out=0,in=90] (1,0)node[below]{${-}k$};
\end{tikzpicture}
\hspace*{-0.1cm}\mapsto
\capr
,
\begin{tikzpicture}[anchorbase,scale=1]
\draw[usual,directed=0.99] (1,0)node[below]{$k$} to[out=90,in=0] (0.5,0.5)node[above]{\phantom{k}} to[out=180,in=90] (0,0)node[below]{${-}k$};
\end{tikzpicture}
\hspace*{-0.1cm}\mapsto
\capl
,
\begin{tikzpicture}[anchorbase,scale=1]
\draw[usual,directed=0.99] (0,0)node[above]{${-}k$} to[out=270,in=180] (0.5,-0.5)node[below]{\phantom{k}} to[out=0,in=270] (1,0)node[above]{$k$};
\end{tikzpicture}
\hspace*{-0.1cm}\mapsto
\cupr
,
\begin{tikzpicture}[anchorbase,scale=1]
\draw[usual,directed=0.99] (1,0)node[above]{${-}k$} to[out=270,in=0] (0.5,-0.5)node[below]{\phantom{k}} to[out=180,in=270] (0,0)node[above]{$k$};
\end{tikzpicture}
\hspace*{-0.1cm}\mapsto
\cupl
,$}
\end{gather*}
and the downwards merges and splits to the respective mates.

\begin{Theorem}\label{T:EquivalenceGln}
The functor $\pfunctor$ is an equivalence of pivotal categories, and
it induces an equivalence of pivotal categories between the additive idempotent completion of $\web$ and $\rep[\gln]$.
\end{Theorem}

\begin{proof}
That $\pfunctor$ is fully faithful is \cite[Theorem 3.3.1]{CaKaMo-webs-skew-howe} and the fact 
that the relevant hom-spaces stay of the same dimension 
when restricting from $\gln$ to $\sln$. The second claim follows from \autoref{L:LeviGlnRepSemi}.(b), or 
\cite[Theorem 3.3.1]{CaKaMo-webs-skew-howe} and flatness of restriction 
from $\gln$ to $\sln$.
\end{proof}


\subsection{The statement}\label{SS:EquivalenceStatement}


We now extend the functor $\pfunctor$ into a functor
\begin{gather*}
\apfunctor=\apfunctor^{ext}(\levi)\colon\aweb\to\fund[\levi].
\end{gather*}
On objects and the generators of $\web$ the functor $\apfunctor$ is defined to be $\pfunctor$. We define $\apfunctor$ on the two coils $\winding$ 
and $\iwinding$ for $k_{1}=1$ by
\begin{gather*}
\begin{tikzpicture}[anchorbase,scale=1]
\draw[usual,rounded corners] (0,0)node[below]{$1$} to (0,0.5) to (-0.5,0.5);
\draw[usual,rounded corners] (2.5,0.5) to (2,0.5) to (2,1)node[above]{$1$};
\draw[usual] (1,0)node[below]{$k_{2}$} to (0,1)node[above]{$k_{2}$};
\draw[usual] (2,0)node[below]{$k_{m}$} to (1,1)node[above]{$k_{m}$};
\node at (1,0.5) {$\dots$};
\draw[affine] (-0.5,0) to (-0.5,1);
\draw[affine] (2.5,0) to (2.5,1);
\end{tikzpicture}
\mapsto
\wmapl
,\quad
\begin{tikzpicture}[anchorbase,scale=1]
\draw[usual,rounded corners] (0,0)node[above]{$1$} to (0,-0.5) to (-0.5,-0.5);
\draw[usual,rounded corners] (2.5,-0.5) to (2,-0.5) to (2,-1)node[below]{$1$};
\draw[usual] (1,0)node[above]{$k_{2}$} to (0,-1)node[below]{$k_{2}$};
\draw[usual] (2,0)node[above]{$k_{m}$} to (1,-1)node[below]{$k_{m}$};
\node at (1,-0.5) {$\dots$};
\draw[affine] (-0.5,0) to (-0.5,-1);
\draw[affine] (2.5,0) to (2.5,-1);
\end{tikzpicture}
\mapsto
\wmapr
.
\end{gather*}

We call the kernel of the functor $\apfunctorstar$ 
the \emph{Levi ideal} 
and denote it by $\lideal$.

\begin{Lemma}\label{L:AnnularWebsLeviEvalCat}
The Levi ideal $\lideal$ is a two-sided $\vcirc$-ideal in $\awebstar$.
\end{Lemma}

\begin{proof}
The kernel is a 
two-sided $\vcirc$-ideal.
\end{proof}

Thus, we get a well-defined category by the following definition.

\begin{Definition}\label{D:AnnularWebsLeviEvalCat}
Let $\lweb$ denote the quotient of $\awebstar$ 
by the Levi ideal $\lideal$.
\end{Definition}

\begin{Theorem}\label{T:EquivalenceMain}
Let $\completion{\placeholder}$ denote additive idempotent completion.
\begin{enumerate}

\item We have the commuting diagram of categories
\begin{gather*}
\begin{tikzcd}[ampersand replacement=\&,column sep=2cm,row sep=1cm]
\aweb\ar[rd,"\apfunctor"']\ar[r,"\text{\autoref{Eq:AnnularWebsLeviEval}}"] \& \aweb[\levi]\ar[r,"\completion{\placeholder}"]\ar[d,"\cong","\apfunctorl{\qpar}{\levi}"'] \& \completion{\aweb[\levi]}\ar[d,"\completion{\apfunctorl{\qpar}{\levi}}","\cong"']
\\
\phantom{.} \& \fund[\levi]\ar[r,"\completion{\placeholder}"'] \& \rep[\levi].
\end{tikzcd}
\end{gather*}
The Levi ideal $\lideal$
is the two-sided $\vcirc$-$\lprod$-ideal generated by the Levi evaluations.
Here $\lprod$ is the pullback of the monoidal structure from $\fund[\levi]$ to 
$\aweb$.

\item We have the commuting diagram of pivotal categories
\begin{gather*}
\begin{tikzcd}[ampersand replacement=\&,column sep=2cm,row sep=1cm]
\awebone\ar[rd,"{\apfunctor[1]}"']\ar[r,"\text{\autoref{Eq:AnnularWebsLeviEval}}"] \& \awebone[\levi]\ar[r,"\completion{\placeholder}"]\ar[d,"\cong","\apfunctorl{1}{\levi}"'] \& \completion{\awebone[\levi]}\ar[d,"\completion{\apfunctorl{1}{\levi}}","\cong"']
\\
\phantom{.} \& \fundone[\levi]\ar[r,"\completion{\placeholder}"'] \& \repone[\levi].
\end{tikzcd}
\end{gather*}
The Levi ideal $\lideal$
is the two-sided $\vcirc$-$\wprod$-ideal generated by the Levi evaluations.

\end{enumerate}
\end{Theorem}

The proof of \autoref{T:EquivalenceMain} is postponed to 
\autoref{S:Proof}, since 
we want to focus on applications of this theorem first.
For the rest of the section 
we assume that \autoref{T:EquivalenceMain} holds with the 
exception of the next subsection where we only assume that $\apfunctor$ is well-defined.


\subsection{Monoidal behavior of the main functor}\label{SS:EquivalenceMonoidal}


The categories $\aweb$ and $\rep[\levi]$ are endowed with monoidal structures that are natural from two different perspectives, as 
explained in \autoref{R:AnnularWebsMonoidal}. However, 
as we will elaborate now, these need not to be the same under the 
equivalence \autoref{T:EquivalenceMain}.

\begin{Lemma}\label{L:EquivalenceMonoidal}
Assume that the functor $\apfunctorstar$ is well-defined.
\begin{enumerate}

\item The functor $\apfunctor[1]$ is pivotal.

\item The functor $\apfunctor$ is not monoidal 
(and thus not pivotal).

\end{enumerate}
\end{Lemma}

\begin{proof}
(a). Let us consider $\qpar=1$. 
In this case the braiding on $\repone[\levi]$ is given by
permutation. Comparing \autoref{Eq:AnnularWebsPropsMonoidal} 
and \autoref{Eq:LeviPropsMonoidal}, and observing that 
coils in $\repone[\levi]$ are permutations, up to diagonal entries of the form $\upar[i]$ 
on blocks, shows that 
the functor $\apfunctor[1]$ is monoidal.
Pivotality is then clear.

(b). For $\qpar\neq 1$, one can check that the images under $\apfunctor$ of  \autoref{Eq:AnnularWebsPropsMonoidal} is not the same as
\autoref{Eq:LeviPropsMonoidal}. Explicitly, taking
$\obstuff{K}=\obstuff{L}=(1)$, 
$\morstuff{f}=\winding$ and $\morstuff{g}=\idmor_{\obstuff{L}}$ 
verifies that $\apfunctor$ is not monoidal. 
\end{proof}

\begin{Remark}\label{R:EquivalenceMonoidalChoice}
We note that \autoref{L:EquivalenceMonoidal} 
shows that the choice which 
side goes over or under in \autoref{Eq:AnnularWebsPropsMonoidal} matters and gives different results on the representation theoretical side.
This indicates that one might need to use the notion 
of module categories rather than monoidal categories 
to describe the representation theory associated to $\aweb$.
This is similar as in, for example, \cite{SaTu-bcd-webs} 
or \cite{EhSt-nw-algebras-howe} 
(via \cite[Remark 12]{MaSt-complex-reflection-groups}), so 
coideal subalgebras might play a role. 

The quantization
issue that we are facing in \autoref{L:EquivalenceMonoidal} is also potentially 
related to the classification of K-matrices as, for example, in \cite{Mu-characters-reflection-equation-algebra} 
where all nondegenerate solutions to the reflection 
equation in the quantum case have satisfy a 
minimal polynomial of order $\leq 2$. Diagrammatically K-matrices 
correspond to coils, but for $\levi=\gln[l_{1}]\oplus\dots\oplus\gln[l_{d}]\subseteq\gln$ these coils need to have a 
minimal polynomial of order $d$.
\end{Remark}

\begin{Remark}\label{R:EquivalenceMonoidal}
One could use the equivalence in 
\autoref{T:EquivalenceMain} to pullback 
the monoidal structure of $\rep[\levi]$, 
resulting in a monoidal structure on $\aweb[\levi]$ 
that is distinct from the one 
we give in \autoref{L:AnnularWebsMonoidal} above. 
This pullback monoidal structure would not 
satisfy the conditions in 
\cite[(2.11)]{MoSa-affinization}, 
which are necessary for it to be 
unique, see \cite[Proposition 2.5]{MoSa-affinization}. 
Conversely, one could push the monoidal structure 
from $\aweb[\levi]$ to $\rep[\levi]$, resulting in a 
monoidal structure on $\rep[\levi]$ satisfying 
\cite[(2.11)]{MoSa-affinization}.
\end{Remark}


\subsection{Ariki--Koike algebras and annular webs}\label{SS:EquivalenceAK}


The \emph{Ariki--Koike algebra} $\ak$ 
from \cite{ArKo-hecke-algebra}, \cite{BrMa-hecke}, 
\cite{Ch-gelfandtzetlin}, using different conventions, is defined as follows.

\begin{Definition}\label{D:EquivalenceGenAK}
Fix $m\in\N$, the number of strands, and
let $\akv$ denote the $\Zv$-algebra with algebra 
generators $\akgens[0],\akgens[1],\dots,\akgens[{m-1}]$ modulo 
the two-sided ideal generated by
\begin{gather*}
\prod_{k=1}^{d}(\akgens[0]-\upar[k])
,\quad
\akgens[0]\akgens[1]\akgens[0]\akgens[1]-\akgens[1]\akgens[0]\akgens[1]\akgens[0],
\\
(\akgens[i]-\qpar)(\akgens[i]+\qpar^{-1})\text{ if }i>0
,\quad
\akgens[i]\akgens[j]\akgens[i]-\akgens[j]\akgens[i]\akgens[j]\text{ if }|i-j|=1
,\quad
\akgens[i]\akgens[j]-\akgens[j]\akgens[i]\text{ if }|i-j|>1,
\end{gather*}
where $i,j\in\set{1,\dots,m-1}$.
\end{Definition}

The Ariki--Koike algebra acts on $1^{\hcirc m}$:

\begin{Proposition}\label{P:EquivalenceGenAK}
We have a surjective $\Kstar$-algebra homomorphism
\begin{gather*}
\akmap\colon
\akstar
\twoheadrightarrow
\End_{\awebstar[\levi]}(1^{\hcirc m})
,\quad
\akgens[0]\mapsto
\begin{tikzpicture}[anchorbase,scale=1]
\draw[usual,rounded corners] (0,0)node[below]{$1$} to (0,0.4) to (-0.5,0.4);
\draw[usual,directed=0.99] (1,0)node[below]{$1^{\hcirc(m-1)}$} to (1,1)node[above]{$1^{\hcirc(m-1)}$};
\draw[usual,crossline] (1.5,0.6) to (0.5,0.6);
\draw[usual,rounded corners,directed=0.99] (0.5,0.6) to (0,0.6) to (0,1)node[above]{$1$};
\draw[affine] (-0.5,0) to (-0.5,1);
\draw[affine] (1.5,0) to (1.5,1);
\end{tikzpicture}
,\quad
\akgens[i]\mapsto
\begin{tikzpicture}[anchorbase,scale=1]
\draw[usual,directed=0.99] (1,0)node[below]{$1$} to (0,1)node[above]{$1$};
\draw[usual,directed=0.99,crossline] (0,0)node[below]{$1$} to (1,1)node[above]{$1$};
\end{tikzpicture}
.
\end{gather*}
Here the bottom left strand of the image of $\akgens[i]$ 
is the $i$th strand from the left.
\end{Proposition}

\begin{proof}
By \autoref{T:EquivalenceMain}, this follows from 
\cite[Theorem 4.2]{SaSh-schur-weyl-ariki-koike} after adjustment of conventions.
\end{proof}

We define the usual \emph{Jucys--Murphy elements} as follows.

\begin{Definition}\label{D:EquivalenceJM}
Define elements of $\akv$ by $\jm[1]=\akgens[0]$
and for $i\in\Z_{\geq 1}$
recursively $\jm[i]=\akgens[i]\jm[{i-1}]\akgens[i]$.
\end{Definition}

Let $\parts$ denote the set of \emph{$d$-partitions of $m$} 
(which we identify with $d$-tuples of Young diagrams in 
the English convention), and 
for $\blam\in\parts$ let 
$\std$ denote the set of all \emph{standard $d$-tableaux of shape $\blam$}.

\begin{Lemma}\label{L:EquivalenceGenBlob}
For all $\blam\in\parts$ there exists a simple $\akstar$-representation 
$V_{\blam}$, 
and these form a complete and nonredundant 
set of simple $\akstar$-representations.
Moreover, we have an $\K$-algebra isomorphism $\akstar\cong\bigoplus_{\blam\in\parts}\End_{\Kstar}(V_{\blam})$. 
Finally, $V_{\blam}$ can be given a 
$\Kstar$-basis $\set{v_{T}|T\in\std}$ such that $\jm$ acts by
\begin{gather}\label{Eq:EquivalenceJM}
\jm\acts v_{T}=\upar[a]\qpar^{2b-2c}v_{T},
\end{gather}
where $a$ is the component of the entry $i$ in $T$, $b$ 
is the column number where $i$ appears and $c$ is the row number of $i$.
\end{Lemma}

\begin{proof}
Our assumptions on the involved parameters 
imply that $\akstar$ is semisimple and the lemma 
follows from \cite[Theorem 3.7]{ArKo-hecke-algebra} 
and results in the same section, {\eg} 
\cite[Proposition 3.16]{ArKo-hecke-algebra}.
\end{proof}

Let $\idealak[{>l_{1},\dots,>l_{d}}]\subset\ak$ denote the two-sided ideal
generated by the idempotents, realizing the Artin--Wedderburn 
decomposition in \autoref{L:EquivalenceGenBlob}, for 
$d$-partitions of $m$ with strictly more than $l_{i}$ rows 
in the $i$th entry. 

\begin{Example}\label{E:EquivalenceBoxNotation}
An important special case is 
$\idealak[{>1,\dots,>1}]\subset\ak$. In this case 
the $d$-partitions indexing the idempotents not in 
$\idealak[{>1,\dots,>1}]$ are of the form
\begin{gather*}
\ytableausetup{centertableaux,textmode,boxsize=0.4cm}
\Big(
\underbrace{\begin{ytableau}
\phantom{a} & \none[\dots] & \phantom{a}
\end{ytableau}}_{k_{1}}
,\dots,
\underbrace{\begin{ytableau}
\phantom{a} & \none[\dots] & \phantom{a}
\end{ytableau}}_{k_{d}}
\Big)
\end{gather*}
for $(k_{1},\dots,k_{d})\in\N^{d}$.
These are so-called one row $d$-partitions.
\end{Example}

\begin{Remark}\label{R:EquivalenceBoxNotation}
The image of 
$\idealak[{>l_{1},\dots,>l_{d}}]\subset\ak$
under $\akmap$ from \autoref{P:EquivalenceGenAK} is 
forcing a condition on the minimal polynomial of coils.
\end{Remark}

The following definition appears in \cite[Section 2C]{LaVa-schur-weyl-ariki-koike}:

\begin{Definition}\label{D:EquivalenceGenBlobk}
Let $\bloblevi$ be the algebra quotient of $\ak$ 
by $\idealak[{>l_{1},\dots,>l_{d}}]$.
\end{Definition}

\begin{Proposition}\label{P:EquivalenceGenBlobk}
The map $\akmap$ from \autoref{P:EquivalenceGenAK} 
induces a $\Kstar$-algebra isomorphism
\begin{gather*}
\blmaplevi\colon\bloblevistar
\xrightarrow{\cong}
\End_{\awebstar[\levi]}(1^{\hcirc m})
.
\end{gather*}
Thus, the kernel of $\akmap$ is $\ker(\akmap)=\idealak[{>l_{1},\dots,>l_{d}}]$.
\end{Proposition}

\begin{proof}
As in the proof of \autoref{P:EquivalenceGenAK}.
\end{proof}


\subsection{Cartan subalgebras and generalized blob algebras}\label{SS:EquivalenceBlob}


We now consider the case of the Cartan subalgebra in detail.

\begin{Proposition}\label{P:EquivalenceMainCartan}
The functor $\apfunctor$ descends to an equivalence of categories $\apfunctorl{\qpar}{\cartan}\colon\aweb[\cartan]\to\fund[\cartan]$, and 
it induces an equivalence of categories between the additive idempotent completion of $\aweb[\cartan]$ and $\rep[\cartan]$.
\end{Proposition}

\begin{proof}
Directly from \autoref{T:EquivalenceMain} 
and \autoref{L:EquivalenceMonoidal}.
\end{proof}

\begin{Remark}\label{R:EquivalenceMainCartan}
\autoref{P:EquivalenceMainCartan} should be compared with 
\cite[Corollary 43]{QuWe-extremal-projectors-2}, and can be seen 
as a quantum version of that corollary.
\end{Remark}

Now recall the so-called \emph{generalized blob algebra} $\blob$
in the sense of \cite{MaWo-generalized-blob} (which is a special case 
of \autoref{D:EquivalenceGenBlobk}):

\begin{Definition}\label{D:EquivalenceGenBlob}
Let $\blob$ be the algebra quotient of $\ak$ 
by $\idealak[{>1,\dots,>1}]$.
\end{Definition}

\begin{Proposition}\label{P:EquivalenceGenBlob}
The map $\akmap$ from \autoref{P:EquivalenceGenAK} 
induces a $\Kstar$-algebra isomorphism
\begin{gather*}
\blmap\colon\blobstar
\xrightarrow{\cong}
\End_{\awebstar[\cartan]}(1^{\hcirc m})
.
\end{gather*}
Thus, the kernel of $\akmap$ is $\ker(\akmap)=\idealak[{>1,\dots,>1}]$.
\end{Proposition}

\begin{proof}
By \autoref{P:EquivalenceMainCartan}, this follows from \cite[Theorem 3.1]{ArTeYa-schur-weyl-ariki-koike}.
\end{proof}


\subsection{On two conjectures about end-spaces in annular webs}\label{SS:EquivalenceConjectures}


Recall from \autoref{SS:AnnularWebsCircles} that $\epoly$ 
denotes the $k$th elementary symmetric polynomial in $d$ variables. 
In the special case of the Cartan subalgebra we have $d=n$, and we 
let $\epoly=\epoly(\upar[1],\dots,\upar[d])$.

\begin{Definition}\label{D:EquivalenceGenBlobRelations}
Define $\epolyg{k}{1}=\epoly$ and for $i\in\Z_{\geq 1}$
recursively $\epolyg{k}{i}=\epolyg{k}{i-1}+(\qpar^{2}-1)
(\jm[i-1]\epolyg{k-1}{i-1}-\jm[i-1]^{2}\epolyg{k-2}{i-1}+\dots\pm\jm[i-1]^{k})$, and let $\akideal[>1]\subset\ak$ denote the two-sided ideal generated by
$\relement[i]=\jm[i]^{d}-\epolyg{1}{i}\jm[i]^{d-1}+\dots\pm\epolyg{d}{i}$. Let 
$\akideal[2]\subset\ak$ denote the two-sided ideal generated by $\relement[2]$.
\end{Definition}

\begin{Proposition}\label{P:EquivalenceGenBlobSecond}
The kernel of $\akmap$ 
is alternatively given by $\ker(\akmap)=\akideal[2]=\akideal[>1]$. 
The same holds for $\qpar=1$.
\end{Proposition}

\begin{proof}
We start by making two claims.	
	
Claim 1. $\relement[2]$ acts on $V_{\blam}$ as zero if and only if 
the $d$-partition $\blam$ has at most one row per component.

Proof of the Claim 1. The case where 
$\blam$ has one node is easy, so assume that 
$\blam$ has at least two nodes. We use that $V_{\blam}$ has a $\K$-basis 
given by $v_{T}$ for $T$ a standard $d$-tableaux of shape $\blam$ 
on which $\jm$ acts by \autoref{Eq:EquivalenceJM}.
Using \autoref{Eq:EquivalenceJM}, we can calculate the 
action of $\relement[2]$ on the $\K$-basis 
given by the $v_{T}$. There are three cases 
depending on the positions of $1$ and $2$ in $T$ one needs to check:
\begin{gather*}
\ytableausetup{centertableaux,textmode,boxsize=0.4cm}
\text{(i)}\colon
\big(\dots,
\begin{ytableau}
1
\end{ytableau}
,\dots,
\begin{ytableau}
2
\end{ytableau}
,\dots\big)
\text{ or {\vive}}
,\quad
\text{(ii)}\colon
\big(
\dots,
\begin{ytableau}
1 & 2
\end{ytableau}
,\dots
\big)
,\quad
\text{(iii)}\colon
\Big(
\dots,
\begin{ytableau}
1 \\ 2
\end{ytableau}
,\dots
\Big)
.
\end{gather*}
All of these are annoying but 
straightforward calculations, and details are omitted.

Claim 2. $\akideal[2]=\akideal[>1]$. 

Proof of the Claim 2. To see this, we note that 
on $v_{T}\in V_{\blam}$ for $\blam$ a one row $d$-partition and $T$ a standard $d$-tableau of shape $\blam$ we have
\begin{gather*}
\epolyg{k}{i}\acts v_{T}=\epoly[k](\qpar^{2\alpha_1}\upar[1],\dots,\qpar^{2\alpha_d}\upar[d])v_{T},
\end{gather*}
where $\alpha_{r}=|\set{s<i|s\text{ is in the }r\text{th component of }T}|$.
Using this formula and \autoref{Eq:EquivalenceJM} one 
can recursively check that the above claim holds 
for $\relement[i]$ for $i\geq 2$, which implies that $\akideal[>1]\subset\akideal[2]$, 
and the proof of the claim is complete.

The first claim implies that $\akideal[2]=\idealak[2]$, and this together with the 
second claim and \autoref{P:EquivalenceGenBlob} proves the lemma.
\end{proof}

Recall $\aff[\web]$ from \autoref{L:AnnularWebsAffinization}. 
Mimicking the construction of $\awebstar[\levi]$ from $\awebstar$ 
as in \autoref{D:AnnularWebsLeviEvalCat}, we denote by 
$\aff[\webstar]/\text{LI}$ (LI stand for Levi ideal) the quotient of 
$\aff[\web]$ by essential circles.
The following will be compared with \cite[Conjectures 10.2 and 10.3]{CaKa-q-satake-sln} in \autoref{P:EquivalenceGenCKTwo} below.

\begin{Proposition}\label{P:EquivalenceGenCK}
We have $\Kstar$-algebra isomorphisms
\begin{gather*}
\akstar/\akideal[2]
=
\akstar/\akideal[>1]
\xrightarrow{\cong}
\End_{\awebstar[\levi]}(1^{\hcirc m})
\xrightarrow{\cong}
\End_{\aff[\webstar]/\text{LI}}(1^{\hcirc m})
,
\end{gather*}
with the first map being induced by \autoref{P:EquivalenceGenAK} 
and the second being induced by \autoref{L:AnnularWebsAffinization}.
\end{Proposition}

\begin{proof}
Combine \autoref{L:AnnularWebsAffinization} and
\autoref{P:EquivalenceGenBlobSecond}.
\end{proof}

\begin{Proposition}\label{P:EquivalenceGenCKTwo}
\autoref{P:EquivalenceGenCK} answers \cite[Conjectures 10.2 and 10.3]{CaKa-q-satake-sln} affirmatively (up to different ground rings; we address that in \autoref{R:EquivalenceGenBlobThird} below).
\end{Proposition}

\begin{proof}
The common object of interest when comparing 
\cite{CaKa-q-satake-sln} to this paper is the annular web category 
and its various flavors:
\begin{enumerate}[itemsep=0.15cm,label=\emph{\upshape(\roman*)}]

\item First, we have $\awebstar$ and 
$\awebstar[\levi]$ and these are the main objects of study in this paper.

\item We also have $\aff[\web]$
and $\aff[\webstar]/\text{LI}$, the versions defined by affinization. These are studied in \cite{CaKa-q-satake-sln} 
for $\sln$ instead of $\gln$.

\end{enumerate}
We first note that the difference between $\sln$ instead of $\gln$ plays no key role in the sense that the relevant hom-spaces are of the same dimension and all maps are defined {\ver}. On the representation theoretical side this is well-known, for the webs see {\eg} \cite[Remark 1.1]{TuVaWe-super-howe}.

In \cite[Section 10]{CaKa-q-satake-sln} Cautis--Kamnitzer define
a map from the affine Hecke algebra $\ak[\text{aff}]$ to 
the category $\aff[\web]$. In \cite[Conjecture 10.3]{CaKa-q-satake-sln} they conjecture what the kernel of this map is. In 
\autoref{P:EquivalenceGenCK} we identify the kernel of the 
respective map from the Ariki--Koike algebra $\akstar$ to $\awebstar[\levi]$.
Thus, we get the following comparison diagram:
\begin{gather*}
\begin{tikzcd}[ampersand replacement=\&,column sep=2.5cm,row sep=1cm]
\ak[\text{aff}]\ar[r,"\scalebox{0.6}{\text{map in \cite[Section 10]{CaKa-q-satake-sln}}}"]\ar[d,two heads]
\&
\aff\ar[r,"\cong","\text{\autoref{P:EquivalenceGenCK}}"',<->]\ar[d,two heads] \& \awebstar\ar[d,two heads]
\\
\akstar\ar[r,"\scalebox{0.6}{\text{map in \autoref{P:EquivalenceGenAK}}}","\text{via the right equivalence}"']
\&
\aff[\webstar]/\text{LI}\ar[r,"\text{\autoref{P:EquivalenceGenCK}}","\cong"',<->] \& \awebstar[\levi].
\end{tikzcd}
\end{gather*}
Comparison of definitions shows that this diagram commutes.
Similarly for \cite[Conjecture 10.2]{CaKa-q-satake-sln} which 
follows from the $\qpar=1$ version of the above.
\end{proof}

\begin{Remark}\label{R:EquivalenceGenBlobThird}
Recall from \autoref{N:PrelimField} that we work over a field
containing variables $\U=\set{\upar[1],\dots,\upar[d]}$ as well as
$\U^{-1}$. We use this crucially in \autoref{Eq:AnnularWebsLeviEval}
where we evaluate essential circles to elementary symmetric polynomials 
$\epoly$ in these variables and their inverses.

On the other hand, the ground used in \cite[Section 10]{CaKa-q-satake-sln} is
\begin{gather*}
E=\C(\qpar)[\efunction[1],\dots,\efunction[n]]
\end{gather*}
where $\efunction$ 
is the $k$th elementary symmetric function (not the polynomial), and these elementary symmetric function compare to our variables $\U$.

Note that \cite[Section 10]{CaKa-q-satake-sln} does not have $\U^{-1}$ and in this sense \autoref{P:EquivalenceGenCK} and \cite[Conjectures 10.2 and 10.3]{CaKa-q-satake-sln} are strictly speaking not comparable.
\end{Remark}


\subsection{Working integrally}\label{SS:EquivalenceIntegral}


Note that our main result \autoref{T:EquivalenceMain} 
is not stated or proven over $\Zv$, and 
we work over $\K$ in which case $\aweb[\levi]$ 
and $\rep[\levi]$ are semisimple. Working \emph{integrally}, that is, over
$\Zv$ or even $\Z[\vpar,\vpar^{-1},\U]$ needs some nontrivial extra steps:
\begin{enumerate}[label=\emph{\upshape(\roman*)}]

\item \autoref{T:EquivalenceGln} works over $\Zv$, see \cite[Theorem 2.58]{El-ladders-clasps}, which uses the light ladder strategy from that 
paper and \cite[Theorem 3.1]{AnStTu-cellular-tilting}.
Passing to an appropriate field 
({\eg} $\overline{\F_{p}}$ for a prime $p$) 
one gets an equivalence of pivotal categories between the additive idempotent completion of $\web$ and the 
category of tilting modules $\tilt$. 

\item The relation from diagram categories to tilting modules 
is a folk observation in the field, see {\eg} \cite[Theorem 2.58]{El-ladders-clasps}, \cite[Section 5A]{AnStTu-cellular-tilting}, 
\cite[Proposition 2.28]{TuWe-quiver-tilting} or \cite[Theorem 1.1]{Bo-c2-tilting} for some examples.

\item Thus, it is tempting to conjecture that 
integral versions of \autoref{T:EquivalenceMain} 
and its consequences involve $\tilt[\levi]$, {\eg} under appropriate 
assumptions on the underlying field 
the additive idempotent completion of $\awebv[\levi]$ should be 
equivalent to $\tilt[\levi]$. However, there is a nontrivial catch: 
the quantization does not behave very well, see {\eg} \autoref{SS:EquivalenceMonoidal} or \cite[Section 10]{CaKa-q-satake-sln}. 
As sketched in \autoref{R:EquivalenceMonoidalChoice}, this might 
indicate that quantum groups are not the correct objects to use 
in this setting.

\item Note that the blob algebra is not defined integrally, and 
it is also not clear from the definition how to work integrally.
Let us however point out that the description of $\blob$ from 
\cite[Theorem 2.15]{LaVa-schur-weyl-ariki-koike} works integrally 
and might play a role in the integral story.

\end{enumerate}
We decided not to pursue these points further in this work.


\section{Proof of the main theorem}\label{S:Proof}



\subsection{Well-definedness}\label{SS:ProofWellDefined}


Recall that the images of the coils are defined by explosion, that is, we mimic \autoref{Eq:AnnularWebsProperties} on the side of the representation theory. We define $\morstuff{u}\colon\vecrep^{\otimes k}\to\vecrep^{\otimes k}$ as the map sending $v_{i_{1},j_{1}}\otimes\cdots\otimes v_{i_{k},j_{k}}$ to $\upar[{i_{1}}]\cdots \upar[{i_{k}}]v_{i_{1},j_{1}}\otimes\cdots\otimes v_{i_{k},j_{k}}$. We use a box to denote this map in illustrations.

\begin{Lemma}\label{L:ProofLeviCoil}
We have
\begin{gather*}
\apfunctor(\winding[{(k,\obstuff{K})}]) = \qfac{k}^{-1}\cdot
\scalebox{0.75}{$\begin{tikzpicture}[anchorbase,scale=0.8]
\draw[spinach!50,fill=spinach!50] (-1.25,0.5) rectangle (0.25,1);
\draw[tomato!50,fill=tomato!50] (-1.25,1) rectangle (0.25,1.5);
\draw[usual] (0.5,-0.5)node[below]{$\obstuff{K}$} to (0.5,2.5);
\draw[usual] (-0.5,-0.5)node[below]{$k$} to (-0.5,0) to (-1,0.5);
\draw[usual] (-0.5,0) to (0,0.5);
\draw[usual] (-1,1.5) to (-0.5,2) to (-0.5,2.5);
\draw[usual] (0,1.5) to (-0.5,2);
\draw[usual] (0.5,2.5) to (-0.5,3.5)node[above]{$\phantom{\obstuff{K}}$};
\draw[usual,levi,crossline] (-0.5,2.5) to (0.5,3.5)node[above]{$\phantom{k}$};
\node at (-0.5,0.3) {$\dots$};
\node at (-0.5,1.25) {$\fulltwist$};
\node at (-0.5,0.75) {$\morstuff{u}$};
\node at (-0.5,1.7) {$\dots$};
\end{tikzpicture}$}
.
\end{gather*}
\end{Lemma}

\begin{proof}
Using the Levi crossings introduced in \autoref{SS:LeviCrossings}, we easily check that $\apfunctor$ sends the coil $\winding[{(k,\obstuff{K})}]$ to the claimed picture.
\end{proof}

\begin{Lemma}\label{L:ProofFullTwistFunnyFormula}
In $\fund[\levi]$ we have
\begin{gather*}
\qfac{k}^{-1}\qfac{l}^{-1}\cdot
\scalebox{0.75}{$\begin{tikzpicture}[anchorbase,scale=0.8]
\draw[tomato!50,fill=tomato!50] (-1,1) rectangle (1,1.5);
\draw[usual] (0,-0.25)node[below]{$k+l$} to (0,0.5) to (-0.5,1);
\draw[usual] (0,0.5) to (0.5,1);
\draw[usual] (-0.9,1.5) to (-0.5,2) to (-0.5,2.75)node[above]{$k$};
\draw[usual] (-0.1,1.5) to (-0.5,2);
\draw[usual] (0.1,1.5) to (0.5,2) to (0.5,2.75)node[above]{$l$};
\draw[usual] (0.9,1.5) to (0.5,2);
\node at (0,1.25) {$\fulltwist$};
\node at (0,0.8) {$\dots$};
\node at (-0.5,1.7) {\scalebox{0.8}{$\dots$}};
\node at (0.5,1.7) {\scalebox{0.8}{$\dots$}};
\end{tikzpicture}$}
=
\qfac{k+l}^{-1}\cdot
\scalebox{0.75}{$\begin{tikzpicture}[anchorbase,scale=0.8]
\draw[tomato!50,fill=tomato!50] (-1,1) rectangle (1,1.5);
\draw[usual] (0,0)node[below]{$k+l$} to (0,0.5) to (-0.5,1);
\draw[usual] (0,0.5) to (0.5,1);
\draw[usual] (-0.5,1.5) to (0,2) to (0,2.5) to (-0.5,3)node[above]{$k$};
\draw[usual] (0.5,1.5) to (0,2);
\draw[usual] (0,2.5) to (0.5,3)node[above]{$l$};
\node at (0,1.25) {$\fulltwist$};
\node at (0,1.7) {$\dots$};
\node at (0,0.8) {$\dots$};
\end{tikzpicture}$}
\,.
\end{gather*}
\end{Lemma}

\begin{proof}
We first suppose that the strands 
are oriented upward. The image of 
$\expl[k+l]$ is an $\Ul$-subrepresentation of 
$\vecrep^{\otimes{k{+}l}}$. The Levi full 
twist $\fulltwist$ is $\Ul$-invariant, so this subrepresentation
remains invariant. But this 
subrepresentation is isomorphic to $\ext{\qpar}{k{+}l}\vecrep$ 
and since the weight spaces of 
$\ext{\qpar}{k{+}l}\vecrep$ are of dimension one, 
we deduce that the image of a vector $v_{S}$ 
by $\fulltwist\circ\expl[k{+}l]$ is a multiple 
of $\expl[k{+}l](v_{S})$, say $\lscalar_{S}\expl[k{+}l](v_{S})$. 
Hence, $\splitmap{k,l}{k{+}l}\circ\iexpl[k+l]\circ
\fulltwist\circ\expl[k{+}l]$ sends $v_{S}$ to $\lscalar_{S}\splitmap{k,l}{k{+}l}\circ\iexpl[k{+}l]
\circ\expl[k{+}l](v_{S})=\lscalar_{S}\qfac{k{+}l}
\splitmap{k,l}{k{+}l}(v_{S})$ and $(\iexpl[k]\otimes\iexpl[l])\circ\fulltwist\circ\expl[k{+}l]$ 
sends $v_{S}$ to $\lscalar_{S}(\iexpl[k]\otimes\iexpl[l])
\circ\expl[k{+}l](v_{S})
=\lscalar_{S}\qfac{k}\qfac{l}\splitmap{k,l}{k{+}l}(v_{S})$.

The same arguments shows that the equality also holds with strands oriented downward.
\end{proof}

\begin{Lemma}\label{L:ProofInverse}
The relation \autoref{Eq:AnnularWebsRel1} is satisfied after applying the functor $\apfunctor$.
\end{Lemma}

\begin{proof}
Clear.
\end{proof}

\begin{Lemma}\label{L:ProofSplitCoil}
The relation \autoref{Eq:AnnularWebsRel2} is satisfied after applying the functor $\apfunctor$.
\end{Lemma}

\begin{proof}
Using the associativity of splits and merges and the fact that Levi crossings satisfy the Reidemeister III relation, it remains to prove
\begin{gather*}
\qfac{k{+}l}^{-1}\cdot
\scalebox{0.75}{$\begin{tikzpicture}[anchorbase,scale=0.8]
\draw[tomato!50,fill=tomato!50] (-1.25,1) rectangle (0.25,1.5);
\draw[usual] (0.5,0)node[below]{$\obstuff{K}$} to (0.5,2.5);
\draw[usual] (-0.5,0)node[below]{$k{+}l$} to (-0.5,0.5) to (-1,1);
\draw[usual] (-0.5,.5) to (0,1);
\draw[usual] (-1,1.5) to (-0.5,2) to (-0.5,2.5);
\draw[usual] (0,1.5) to (-0.5,2);
\draw[usual] (0.5,2.5) to (-0.5,3.5);
\draw[usual] (-0.5,3.5) to (-0.5,4)node[above]{$\obstuff{K}$};
\draw[usual,levi,crossline] (-0.5,2.5) to (0.5,3.5);
\draw[usual] (0.25,4)node[above]{$k$} to (0.5,3.5) to (0.75,4)node[above]{$l$};
\node at (-0.5,0.8) {$\dots$};
\node at (-0.5,1.25) {$\fulltwist$};
\node at (-0.5,1.7) {$\dots$};
\end{tikzpicture}$}
=
\qfac{k}^{-1}\qfac{l}^{-1}\cdot
\scalebox{0.75}{$\begin{tikzpicture}[anchorbase,scale=0.8]
\draw[tomato!50,fill=tomato!50] (-1,1) rectangle (1,1.5);
\draw[usual] (1.5,0)node[below]{$\obstuff{K}$} to (1.5,2.5) to (-0.5,3.5) to (-0.5,4)node[above]{$\obstuff{K}$};
\draw[usual] (0,0)node[below]{$k+l$} to (0,0.5) to (-0.5,1);
\draw[usual] (0,0.5) to (0.5,1);
\draw[usual] (-0.9,1.5) to (-0.5,2) to (-0.5,3);
\draw[usual] (-0.1,1.5) to (-0.5,2);
\draw[usual,levi,crossline] (-0.5,3) to (0.5,3.5);
\draw[usual] (0.5,3.5) to (0.5,4)node[above]{$k$};
\draw[usual] (0.1,1.5) to (0.5,2) to (0.5,2.5);
\draw[usual] (0.9,1.5) to (0.5,2);
\draw[usual,levi,crossline] (0.5,2.5) to (1.5,3);
\draw[usual] (1.5,3) to (1.5,4)node[above]{$l$};
\node at (0,1.25) {$\fulltwist$};
\node at (0,0.8) {$\dots$};
\node at (-0.5,1.7) {$\dots$};
\node at (0.5,1.7) {$\dots$};
\end{tikzpicture}$}
,
\end{gather*}
which follows immediately from \autoref{L:ProofFullTwistFunnyFormula}.
\end{proof}

We can show similarly that merges slide through coils.

\begin{Lemma}\label{L:ProofCapCoil}
The relation \autoref{Eq:AnnularWebsRel3} 
is satisfied after applying the functor $\apfunctor$.
\end{Lemma}

\begin{proof}
The image of the right-hand side of 
\autoref{Eq:AnnularWebsRel3} is given, 
up to the multiplication by elements of $\U$ 
that cancel out, by
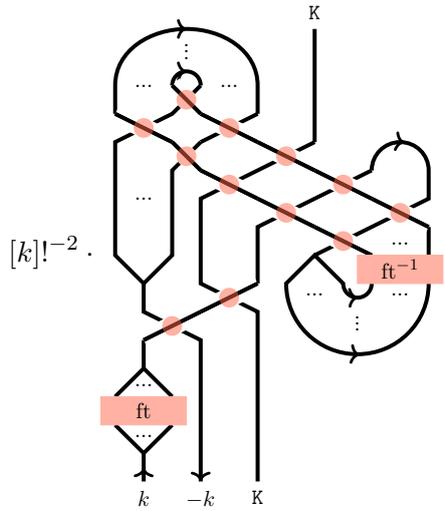
\begin{gather*}
\qfac{k}^{-2}\cdot
\scalebox{0.75}{$\begin{tikzpicture}[anchorbase]
\draw[usual,directed=0.2] (0,0)node[below] {$k$} to (0,0.5) to (-0.5,1);
\draw[usual] (0,0.5) to (0.5,1);
\draw[usual] (-0.5,1.5) to (0,2) to (0,2.5);
\draw[usual] (0.5,1.5) to (0,2);
\draw[usual] (2,3.5) to (2,4.5) to (4,5.5);
\draw[usual,directed=0.2] (4,5.5) to[out=90,in=180] (4.5,6) to[out=0,in=90] (5,5.5) to (5,5) to (3,4);
\draw[usual,rdirected=0.02] (1,0)node[below] {${-}k$} to (1,2.5) to (0,3) to (0,3.5) to (-0.5,4) to (-0.5,6) to (0.5,6.5) to (1,7);
\draw[usual] (0,3.5) to (0.5,4) to (0.5,5.5) to (1,6);
\draw[usual,directed=0.5] (0.5,7) to[out=90,in=180] (0.75,7.25) to[out=0,in=90] (1,7);
\draw[usual,directed=0.5] (-0.5,6.5) to (-0.5,7) to[out=90,in=180] (.75,8) to[out=0,in=90] (2,7) to (2,6.5);
\draw[usual] (1,6) to (2,6.5);
\draw[usual] (2,0)node[below] {$\obstuff{K}$} to (2,3) to (1,3.5) to (1,5) to (3,6);
\draw[usual] (3,6) to (3,8)node[above] {$\obstuff{K}$};
\draw[usual] (5,4.5) to (5,4);
\draw[usual,rdirected=0.341,directed=0.937] (5,3.5) to[out=270,in=0] (3.75,2.25) to[out=180,in=270] (2.5,3.5) to (3,4) to (3.5,3.5) to[out=270,in=180] (3.75,3.25) to[out=0,in=270] (4,3.5);
\draw[usual,levi,crossline] (0,2.5) to (1,3);
\draw[usual,levi,crossline] (1,3) to (2,3.5);
\draw[usual,levi,crossline] (4,4) to (3,4.5);
\draw[usual,levi,crossline] (3,4.5) to (2,5);
\draw[usual,levi,crossline] (5,4.5) to (4,5);
\draw[usual,levi,crossline] (4,5) to (3,5.5);
\draw[usual,levi,crossline] (2,5) to (1,5.5);
\draw[usual,levi,crossline] (3,5.5) to (2,6);
\draw[usual,levi,crossline] (2,6) to (1,6.5);
\draw[usual,levi,crossline] (1,6.5) to (0.5,7);
\draw[usual,levi,crossline] (1,5.5) to (0.5,6);
\draw[usual,levi,crossline] (0.5,6) to (-0.5,6.5);
\draw[tomato!50,fill=tomato!50] (3.75,3.5) rectangle (5.25,4);
\draw[tomato!50,fill=tomato!50] (-0.75,1) rectangle (0.75,1.5);
\node at (0,0.8) {$\dots$};
\node at (0,1.25) {$\fulltwist$};
\node at (0,1.7) {$\dots$};
\node at (0,5) {$\dots$};
\node at (4.5,3.75) {$\fulltwist^{-1}$};
\node at (4.5,3.3) {$\dots$};
\node at (3,3.3) {$\dots$};
\node at (3.75,2.8) {$\vdots$};
\node at (4.5,4.2) {$\dots$};
\node at (0,7) {$\dots$};
\node at (0.75,7.6) {$\vdots$};
\node at (1.5,7) {$\dots$};
\end{tikzpicture}$}
.
\end{gather*}

Using Reidemeister II relations and the fact that 
mates of merges are splits and {\vive}, 
it remains to prove the following equality:
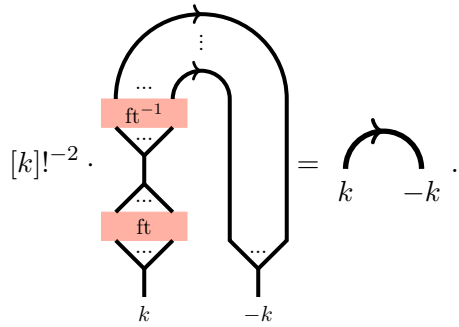
\begin{gather*}
\qfac{k}^{-2}\cdot
\scalebox{0.75}{$\begin{tikzpicture}[anchorbase]
\draw[tomato!50,fill=tomato!50] (-0.75,1) rectangle (0.75,1.5);
\draw[tomato!50,fill=tomato!50] (-0.75,3) rectangle (0.75,3.5);
\draw[usual] (0,0)node[below] {$k$} to (0,0.5) to (-0.5,1);
\draw[usual] (0,0.5) to (0.5,1);
\draw[usual] (-0.5,1.5) to (0,2) to (0.5,1.5);
\draw[usual] (0,2) to (0,2.5);
\draw[usual] (-0.5,3) to (0,2.5) to (0.5,3);
\draw[usual,directed=0.15] (0.5,3.5) to[out=90,in=180] (1,4) to[out=0,in=90] (1.5,3.5) to (1.5,1) to (2,0.5) to (2,0)node[below] {${-}k$};
\draw[usual,directed=0.3] (-0.5,3.5) to[out=90,in=180] (1,5) to[out=0,in=90] (2.5,3.5) to (2.5,1) to (2,0.5);
\node at (0,0.8) {$\dots$};
\node at (0,1.25) {$\fulltwist$};
\node at (0,1.7) {$\dots$};
\node at (0,2.8) {$\dots$};
\node at (0,3.25) {$\fulltwist^{-1}$};
\node at (0,3.7) {$\dots$};
\node at (1,4.5) {$\vdots$};
\node at (2,0.8) {$\dots$};
\end{tikzpicture}$}
=
\begin{tikzpicture}[anchorbase]
\draw[usual,directed=.5] (0,-1)node[below] {$k$} to[out=90,in=180] (0.5,-0.5) to[out=0,in=90] (1,-1)node[below] {${-}k$};
\end{tikzpicture}
.
\end{gather*}
In order to get rid of the explosions between the Levi full twists, we apply repeatedly \autoref{L:ProofFullTwistFunnyFormula} with $l=1$ and use the fact that the Levi full twist on $k$ strands can be obtain from the Levi full twist on $k-1$ strands and some extra crossings. We can conclude using explosion as in \autoref{L:WebsGlnProperties}.

The argument for the leftward oriented cap is similar and omitted.
\end{proof}

We can show similarly that cups slide through coils.

\begin{Lemma}\label{L:ProofEssentialCircles}
The relations in \autoref{Eq:AnnularWebsLeviEval} are in the kernel of $\apfunctor$.
\end{Lemma}

\begin{proof}
The image of the leftward oriented essential circle with a strand of thickness $k$ through $\apfunctor$ is given by
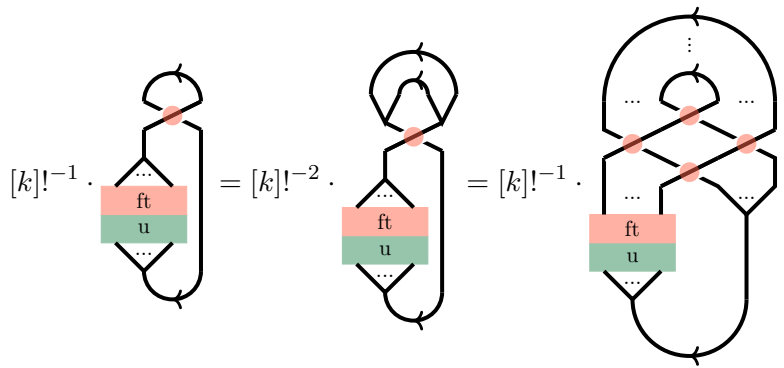
\begin{gather*}
\qfac{k}^{-1}\cdot
\scalebox{0.75}{$\begin{tikzpicture}[anchorbase]
\draw[spinach!50,fill=spinach!50] (-0.75,0.5) rectangle (0.75,1);
\draw[tomato!50,fill=tomato!50] (-0.75,1) rectangle (0.75,1.5);
\draw[usual] (0,3) to (1,2.5);
\draw[usual] (1,2.5) to (1,0);
\draw[usual,directed=0.5] (1,0) to[out=270,in=0] (0.5,-0.5) to[out=180,in=270] (0,0);
\draw[usual] (-0.5,0.5) to (0,0) to (0.5,0.5);
\draw[usual] (-0.5,1.5) to (0,2) to (0.5,1.5);
\draw[usual,levi,crossline] (0,2.5) to (1,3);
\draw[usual] (0,2) to (0,2.5);
\draw[usual,directed=0.5] (1,3) to[out=90,in=0] (0.5,3.5) to[out=180,in=90] (0,3);
\node at (0,0.3) {$\dots$};
\node at (0,0.75) {$\morstuff{u}$};
\node at (0,1.25) {$\fulltwist$};
\node at (0,1.7) {$\dots$};
\end{tikzpicture}$}
=
\qfac{k}^{-2}\cdot
\scalebox{0.75}{$\begin{tikzpicture}[anchorbase]
\draw[spinach!50,fill=spinach!50] (-0.75,0.5) rectangle (0.75,1);
\draw[tomato!50,fill=tomato!50] (-0.75,1) rectangle (0.75,1.5);
\draw[usual] (0,3) to (1,2.5);
\draw[usual] (1,2.5) to (1,0);
\draw[usual,directed=0.5] (1,0) to[out=270,in=0] (0.5,-0.5) to[out=180,in=270] (0,0);
\draw[usual] (-0.5,0.5) to (0,0) to (0.5,0.5);
\draw[usual] (-0.5,1.5) to (0,2) to (0.5,1.5);
\draw[usual,levi,crossline] (0,2.5) to (1,3);
\draw[usual] (0,2) to (0,2.5);
\draw[usual] (-0.25,3.5) to (0,3) to (0.25,3.5);
\draw[usual] (0.75,3.5) to (1,3) to (1.25,3.5);
\draw[usual,directed=0.5] (1.25,3.5) to[out=90,in=0] (0.5,4.25) to[out=180,in=90] (-0.25,3.5);
\draw[usual,directed=0.5] (0.75,3.5) to[out=90,in=0] (0.5,3.75) to[out=180,in=90] (0.25,3.5);
\node at (0,0.3) {$\dots$};
\node at (0,0.75) {$\morstuff{u}$};
\node at (0,1.25) {$\fulltwist$};
\node at (0,1.7) {$\dots$};
\end{tikzpicture}$}
=
\qfac{k}^{-1}\cdot
\scalebox{0.75}{$\begin{tikzpicture}[anchorbase]
\draw[spinach!50,fill=spinach!50] (-0.75,0.5) rectangle (0.75,1);
\draw[tomato!50,fill=tomato!50] (-0.75,1) rectangle (0.75,1.5);
\draw[usual] (2,0) to (2,1.5) to (1.5,2) to (-0.5,3) to (-0.5,3.5);
\draw[usual] (2,1.5) to (2.5,2) to (2.5,2.5) to (0.5,3.5);
\draw[usual,directed=0.5] (2,0) to[out=270,in=0] (1,-1) to[out=180,in=270] (0,0);
\draw[usual] (-0.5,0.5) to (0,0) to (0.5,0.5);
\draw[usual] (-0.5,1.5) to (-0.5,2.5);
\draw[usual] (0.5,1.5) to (0.5,2);
\draw[usual,levi,crossline] (0.5,2) to (1.5,2.5);
\draw[usual,levi,crossline] (1.5,2.5) to (2.5,3);
\draw[usual,levi,crossline] (-0.5,2.5) to (0.5,3);
\draw[usual,levi,crossline] (0.5,3) to (1.5,3.5);
\draw[usual] (2.5,3) to (2.5,3.5);
\draw[usual,directed=0.5] (2.5,3.5) to[out=90,in=0] (1,5) to[out=180,in=90] (-0.5,3.5);
\draw[usual,directed=0.5] (1.5,3.5) to[out=90,in=0] (1,4) to[out=180,in=90] (0.5,3.5);
\node at (0,0.3) {$\dots$};
\node at (0,0.75) {$\morstuff{u}$};
\node at (0,1.25) {$\fulltwist$};
\node at (0,1.8) {$\dots$};
\node at (2,1.8) {$\dots$};
\node at (1,4.5) {$\vdots$};
\node at (0,3.5) {$\dots$};
\node at (2,3.5) {$\dots$};
\end{tikzpicture}$}
.
\end{gather*}
In this calculation the first equality is obtained from explosion of strands, and the last equality is obtained from \autoref{L:ProofFullTwistFunnyFormula}.

Now, since
\begin{gather*}
\scalebox{0.75}{$\begin{tikzpicture}[anchorbase]
\draw[usual] (-0.5,1.5) to (-0.5,2.5);
\draw[usual] (2.5,1.5) to (2.5,2.5);
\draw[usual] (1.5,1.5) to (1.5,2) to (-0.5,3) to (-0.5,3.5);
\draw[usual] (2.5,1.5) to (2.5,2.5) to (0.5,3.5);
\draw[usual] (0.5,1.5) to (0.5,2);
\draw[usual,levi,crossline] (0.5,2) to (1.5,2.5);
\draw[usual,levi,crossline] (1.5,2.5) to (2.5,3);
\draw[usual,levi,crossline] (-0.5,2.5) to (0.5,3);
\draw[usual,levi,crossline] (0.5,3) to (1.5,3.5);
\draw[usual] (2.5,3) to (2.5,3.5);
\draw[usual,directed=0.5] (2.5,3.5) to[out=90,in=0] (1,5) to[out=180,in=90] (-0.5,3.5);
\draw[usual,directed=0.5] (1.5,3.5) to[out=90,in=0] (1,4) to[out=180,in=90] (0.5,3.5);
\node at (0,1.8) {$\dots$};
\node at (2,1.8) {$\dots$};
\node at (1,4.5) {$\vdots$};
\node at (0,3.5) {$\dots$};
\node at (2,3.5) {$\dots$};
\end{tikzpicture}$}
=
\scalebox{0.75}{$\begin{tikzpicture}[anchorbase]
\draw[tomato!50,fill=tomato!50] (-0.25,-0.5) rectangle (2.25,0);
\draw[usual] (0,-0.5) to (0,-1);
\draw[usual] (2,-0.5) to (2,-1);
\draw[usual] (3,1) to (3,-1);
\draw[usual] (4,1) to (4,-1);
\draw[usual] (0,0) to (0,0.5);
\draw[usual] (-0.5,0.5) to (0,1);
\draw[usual] (-0.5,1) to[out=90,in=0] (-0.75,1.25) to[out=180,in=90] (-1,1) to (-1,0.5) to[out=270,in=180] (-0.75,0.25) to[out=0,in=270] (-0.5,0.5);
\draw[usual,levi,crossline] (0,0.5) to (-0.5,1);
\draw[usual] (2,0) to (2,0.5);
\draw[usual] (1.5,0.5) to (2,1);
\draw[usual] (1.5,1) to[out=90,in=0] (1.25,1.25) to[out=180,in=90] (1,1) to (1,0.5) to[out=270,in=180] (1.25,0.25) to[out=0,in=270] (1.5,0.5);
\draw[usual,levi,crossline] (2,0.5) to (1.5,1);
\draw[usual,directed=.5] (0,1) to[out=90,in=180] (2,3) to[out=0,in=90] (4,1);
\draw[usual,directed=.5] (2,1) to[out=90,in=180] (2.5,1.5) to[out=0,in=90] (3,1);
\node at (0.5,0.75) {$\dots$};
\node at (1,-0.25) {$\fulltwist^{-1}$};
\node at (1,-0.75) {$\dots$};
\node at (3.5,-0.75) {$\dots$};
\node at (2.5,2.25) {$\vdots$};
\end{tikzpicture}$}
,
\end{gather*}
it remains to compute the following scalar:
\begin{gather*}
\qfac{k}^{-1}\cdot
\scalebox{0.75}{$\begin{tikzpicture}[anchorbase]
\draw[spinach!50,fill=spinach!50] (-1.25,0.5) rectangle (1.25,1);
\draw[usual] (2,0) to (2,3);
\draw[usual,directed=0.5] (2,0) to[out=270,in=0] (1,-1) to[out=180,in=270] (0,0);
\draw[usual] (-1,0.5) to (0,0) to (1,0.5);
\draw[usual] (-1,1) to (-1,1.5);
\draw[usual] (-1.5,2) to[out=90,in=0] (-1.75,2.25) to[out=180,in=90] (-2,2) to (-2,1.5) to[out=270,in=180] (-1.75,1.25) to[out=0,in=270] (-1.5,1.5) to (-1,2);
\draw[usual,levi,crossline] (-1,1.5) to (-1.5,2);
\draw[usual] (1,1) to (1,1.5);
\draw[usual] (0.5,2) to[out=90,in=0] (0.25,2.25) to[out=180,in=90] (0,2) to (0,1.5) to[out=270,in=180] (0.25,1.25) to[out=0,in=270] (0.5,1.5) to (1,2);
\draw[usual,levi,crossline] (1,1.5) to (0.5,2);
\draw[usual] (-1,2) to (-1,2.5) to (0,3) to (1,2.5) to (1,2);
\draw[usual] (0,3) to[out=90,in=180] (1,4) to[out=0,in=90] (2,3);
\node at (0,0.3) {$\dots$};
\node at (0,0.75) {$\morstuff{u}$};
\node at (-0.5,1.75) {$\dots$};
\end{tikzpicture}$}
.
\end{gather*}
We easily check that the twist sends $v_{i,j}$ to $\qpar^{-n+2(l_{1}+\cdots +l_{d})}v_{i,j}$ and therefore the previous scalar is equal to
\begin{align*}
&\sum_{k_{1}+\cdots+k_d=k}\sum_{\substack{1\leq j_{1,1}<\cdots<j_{1,k_{1}}\leq l_{1}\\\cdots\\1\leq j_{d,1}<\cdots<j_{d,k_{d}}\leq l_{d}}}\upar[1]^{k_{1}}\cdots \upar[d]^{k_d}\prod_{i=1}^d \qpar^{-n+2(l_{1}+\cdots +l_d)k_{i}+\sum_{r=1}^{k_{i}}(n+1-2(l_{1}+\cdots+l_{i-1}+j_{i,r}))}\\
=&\sum_{k_{1}+\cdots+k_d=k}\sum_{\substack{1\leq j_{1,1}<\cdots<j_{1,k_{1}}\leq l_{1}\\\cdots\\1\leq j_{d,1}<\cdots<j_{d,k_{d}}\leq l_{d}}}\upar[1]^{k_{1}}\cdots\upar[d]^{k_d}\prod_{i=1}^d \qpar^{k_{i}-\sum_{r=1}^{k_{i}}2j_{i,r}}\\
=&\epoly[k](\qpar^{-1}\upar[1],\qpar^{-3}\upar[1],\dots,\qpar^{-2l_{1}+1}\upar[1],\dots,\qpar^{-1}\upar[d],\qpar^{-3}\upar[d],\dots,\qpar^{-2l_{d}+1}\upar[d]),
\end{align*}
which is what we needed to show.
\end{proof}

\begin{Lemma}\label{L:ProofWelldefined}
The functor $\apfunctor$ is well-defined 
and descends to the functor $\apfunctorl{\qpar}{\levi}$. Moreover, 
for $\qpar\neq 1$ the Levi ideal $\lideal$
contains the two-sided $\vcirc$-ideal 
and right $\wprod$-ideal generated by the Levi evaluations, while for $\qpar=1$ the Levi ideal
$\lideal$ contains the two-sided $\vcirc$-$\wprod$-ideal generated by the Levi evaluations.
\end{Lemma}

\begin{proof}
We need to check that the relations in \autoref{D:AnnularWebsDefinition} 
are satisfied and that the Levi evaluations from \autoref{Eq:AnnularWebsLeviEval} are in the kernel of $\apfunctor$. This follows as a combination of \autoref{L:ProofCapCoil}
and \autoref{L:ProofEssentialCircles}. By \autoref{L:EquivalenceMonoidal}, the 
statement about the Levi ideal for $\qpar=1$ follows from that.
To verify the claim for $\qpar\neq 1$, we recall from \autoref{R:AnnularWebsCrossings} that the coils pass in front of the annulus.
Now we observe that, for example, (the $\leftrightsquigarrow$ refers to 
\autoref{R:AnnularWebsCrossings})
\begin{gather}\label{Eq:ProofLeftRight}
\begin{aligned}
\text{right $\wprod$-ideal:}&\quad
\begin{tikzpicture}[anchorbase,scale=1]
\draw[usual] (-0.5,0.5) to (0,0.5) to (0.5,0.5);
\draw[affine] (-0.5,0) to (-0.5,1);
\draw[affine] (0.5,0) to (0.5,1);
\end{tikzpicture}
\wprod
\begin{tikzpicture}[anchorbase,scale=1]
\draw[usual] (0,0) to (0,1);
\draw[affine] (-0.5,0) to (-0.5,1);
\draw[affine] (0.5,0) to (0.5,1);
\end{tikzpicture}
=
\begin{tikzpicture}[anchorbase,scale=1]
\draw[usual] (0,0) to (0,1);
\draw[usual,crossline] (-0.5,0.5) to (0,0.5) to (0.5,0.5);
\draw[affine] (-0.5,0) to (-0.5,1);
\draw[affine] (0.5,0) to (0.5,1);
\end{tikzpicture}
\leftrightsquigarrow
\begin{tikzpicture}[anchorbase,scale=1]
\draw[white] (-0.5,0.5) to[out=180,in=180] (-0.5,1.25);
\draw[usual] (0,0) to (0,1);
\draw[usual,crossline] (-0.5,0.5) to (0,0.5) to (0.5,0.5);
\draw[affine] (-0.5,0) to (-0.5,1);
\draw[affine] (0.5,0) to (0.5,1);
\draw[usual,spinach] (-0.5,0.5) to[out=180,in=180] (-0.5,-0.25) to (0.5,-0.25)to[out=0,in=0] (0.5,0.5);
\end{tikzpicture}
\;,
\\
\text{left $\wprod$-ideal:}&\quad
\begin{tikzpicture}[anchorbase,scale=1]
\draw[usual] (0,0) to (0,1);
\draw[affine] (-0.5,0) to (-0.5,1);
\draw[affine] (0.5,0) to (0.5,1);
\end{tikzpicture}
\wprod
\begin{tikzpicture}[anchorbase,scale=1]
\draw[usual] (-0.5,0.5) to (0,0.5) to (0.5,0.5);
\draw[affine] (-0.5,0) to (-0.5,1);
\draw[affine] (0.5,0) to (0.5,1);
\end{tikzpicture}
=
\begin{tikzpicture}[anchorbase,scale=1]
\draw[usual] (-0.5,0.5) to (0,0.5) to (0.5,0.5);
\draw[usual,crossline] (0,0) to (0,1);
\draw[affine] (-0.5,0) to (-0.5,1);
\draw[affine] (0.5,0) to (0.5,1);
\end{tikzpicture}
\;\leftrightsquigarrow
\begin{tikzpicture}[anchorbase,scale=1]
\draw[white] (-0.5,0.5) to[out=180,in=180] (-0.5,1.25);
\draw[usual] (-0.5,0.5) to (0,0.5) to (0.5,0.5);
\draw[usual,crossline] (0,0) to (0,1);
\draw[affine] (-0.5,0) to (-0.5,1);
\draw[affine] (0.5,0) to (0.5,1);
\draw[usual,spinach] (-0.5,0.5) to[out=180,in=180] (-0.5,-0.25) to (0.5,-0.25)to[out=0,in=0] (0.5,0.5);
\end{tikzpicture}
\;.
\end{aligned}
\end{gather}
In the top picture the essential circle is not pierced 
by the identity morphism, while in the bottom it is.
More generally, multiplying essential 
circles from the right (but not from the left)
by any morphism $\morstuff{f}$ produces a picture with 
$\morstuff{f}$ not interfering with the essential circles and the essential circles are on the outside. Thus, by our convention 
from \autoref{R:AnnularWebsCrossings}, the essentially circle evaluations will not change by right multiplication.
\end{proof}

For the remainder of the paper, we assume \autoref{L:ProofWelldefined}.


\subsection{Proof for \texorpdfstring{$\qpar=1$}{q equal 1}}\label{SS:ProofOne}


Recall that a colored permutation on 
$\set{1,\dots,d}^{k}$ is a permutation 
on $\set{1,\dots,k}$ that preserves the labels. In terms of 
classical permutation diagrams these are 
permutations with colored strands 
such that crossings preserve the colors, {\eg}:
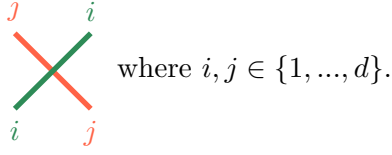
\begin{gather*}
\begin{tikzpicture}[anchorbase,scale=1]
\draw[usual,tomato] (1,0)node[below]{$j$} to (0,1)node[above]{$j$};
\draw[usual,spinach] (0,0)node[below]{$i$} to (1,1)node[above]{$i$};
\end{tikzpicture}
\text{ where }i,j\in\set{1,\dots,d}.
\end{gather*}
Recall also the web block projectors $\wproj$
from \autoref{D:AnnularWebsProjectors} respectively 
the Levi block projectors $\lproj$ from \autoref{D:LeviProjectors}.
We use these to define 
a certain basis in 
the following definition, where $\sym[k]$ is the symmetric group on 
$\set{1,\dots,k}$:

\begin{Definition}\label{D:ProofBasis}
Fix $k\in\N$. Let $s=(s_{1},\dots,s_{k}),t\in\set{1,\dots,d}^{k}$ and let $\sigma$ be a 
colored permutation from $s$ to $t$ such that the longest
word of $\sym[{l_{i}+1}]\subset\sym[k]$ does not appear in the permutation of color $i$. 
(This condition is vacuous for $l_{i}\geq k$ or if the color $i$ 
does not appear strictly more than $l_{i}$ times.)
Define
\begin{gather*}
\basisone=
\sigma\vcirc(\wproj[s_{1}]\wprod\dots\wprod\wproj[s_{d}])
,\quad
\lbasisone=
\sigma\vcirc(\lproj[s_{1}]\lprod\dots\lprod\lproj[s_{d}])
,
\end{gather*}
where we view $\sigma$ as an element of 
$\awebone[\levi]$ and $\fundone[\levi]$, respectively.
The respective sets 
(collecting these elements for 
all $k\in\N$) of these are denoted by $\basissetone$ 
and $\lbasissetone$.
\end{Definition}

\begin{Lemma}\label{L:ProofLinInd}
We have $\apfunctorl{1}{\levi}\big(\basissetone\big)=\lbasissetone$ 
and this set is a $\Kone$-linear 
independent set in $\coprod_{k}\End_{\fundone[\levi]}(\vecrep[1]^{\otimes k})$.
\end{Lemma}

\begin{proof}
We get $\apfunctorl{1}{\levi}\big(\basissetone\big)=\lbasissetone$
directly from \autoref{L:EquivalenceMonoidal} and the definition.

To prove faithfulness, let $\lproj[s]=\wproj[s_{1}]\wprod\dots\wprod\wproj[s_{d}]$.
By \autoref{L:LeviProjectors} and construction we have
$\lbasisone\vcirc\lproj[s]=\lproj[t]\vcirc\lbasisone=\lbasisone$ 
and $\lbasisone\vcirc\lproj[q]=\lproj[r]\vcirc\lbasisone=0$ for 
$q\neq s$ and $r\neq t$. Thus, it suffices to show that 
$(\lbasisone)_{\sigma}$ is $\Kone$-linear 
independent (that is, we can fix $s$ and $t$). 
After sorting the colors of $s$ and $t$, it remains to verify 
faithfulness within a block, {\ie} for $\gln[l_{i}]$. 
Thus, classical theory applies: for $\gln[l_{i}]$ and $\sym[k]$ 
it is known that Schur--Weyl duality gives the generalized Temperley--Lieb
algebra as the endomorphism algebra. This algebra 
admits a description in terms of 
a quotient of $\sym[k]$ by the longest word of $\sym[{l_{i}+1}]$, see {\eg} 
\cite[Section 3]{Ha-murphy-generalized-tl}, and the elements of 
$\lbasissetone$ describe the associated standard-type basis within 
each block.
\end{proof}

\begin{Lemma}\label{L:ProofSpanning}
The set $\basissetone$ $\Kone$-linearly
spans $\coprod_{k}\End_{\awebone[\levi]}(1^{\hcirc k})$. 
Moreover, the set $\lbasissetone$ is a $\Kone$-linear 
spanning set of $\coprod_{k}\End_{\fundone[\levi]}(\vecrep[1]^{\otimes k})$.
\end{Lemma}

\begin{proof}
We first recall that the full antisymmetrizer is zero in $\fundone[\levi]$, 
see {\eg} \cite[Lemma 11]{MaSt-complex-reflection-groups}. 
Thus, the same holds in $\awebone[\levi]$ by 
\autoref{L:ProofLinInd}. That is, in formulas we have
\begin{gather*}
\left(\sum_{\sigma\in\sym[{l_{i}+1}]}
(-1)^{\mathrm{l}(\sigma)}
\sigma(\lproj[i])^{\lprod(l_{i}+1)}=0\right)
\Rightarrow
\left(\sum_{\sigma\in\sym[{l_{i}+1}]}
(-1)^{\mathrm{l}(\sigma)}
\sigma(\wproj[i])^{\wprod(l_{i}+1)}=0\right),
\end{gather*}
and both hold.

To address the first statement of the lemma, observe that 
$\coprod_{k}\End_{\awebone[\levi]}(1^{\hcirc k})$ is generated
as a $\Kone$-algebra by crossings and coils since we can remove essential circles 
by \autoref{L:ProofWelldefined}. Moreover, the sliding relations 
\autoref{Eq:AnnularWebsRel2} and \autoref{Eq:AnnularWebsRel3} imply that 
we have the usual $\Kone$-linear 
spanning set given by first coils and then crossings. Observe 
next that the web block projectors $\Kone$-linear 
span the subalgebra generated by coils, so it remains to see that the 
symmetric group part is $\Kone$-linear spanned by $\sigma$ such that the longest
word of $\sym[l_{i}{+}1]$ does not appear in the permutation of color $i$. This 
however is a consequence of the vanishing of the antisymmetrizer. 

For the second statement of the lemma 
we use \autoref{L:LeviModules}, 
\autoref{L:LeviProjectors} and Schur's lemma. 
(As explained in the proof of 
\autoref{L:LeviReidemeisterOne}, Schur's lemma still holds in this setting 
although $\Kone$ is not necessary algebraically closed.)
\end{proof}

\begin{Lemma}\label{L:ProofBasis}
The set $\basissetone$ is a $\Kone$-basis of 
$\coprod_{k}\End_{\awebone[\levi]}(1^{\hcirc k})$. Furthermore,
the set $\lbasissetone$ is a $\Kone$-basis 
of $\coprod_{k}\End_{\fundone[\levi]}(\vecrep[1]^{\otimes k})$. 
\end{Lemma}

\begin{proof}
Combine \autoref{L:ProofLinInd} and \autoref{L:ProofSpanning}. 
(Note that \autoref{L:ProofLinInd} also proves that $\basissetone$ 
is a $\Kone$-linear independent set.)
\end{proof}

\begin{Proposition}\label{P:ProofFullyFaithful}
The functor $\apfunctorl{1}{\levi}$ is fully faithful and 
the Levi ideal $\lideal$ is the two-sided $\vcirc$-$\wprod$-ideal generated by the Levi evaluations.
\end{Proposition}

\begin{proof}
\autoref{L:WebsExplosion} implies that
we need to show 
that $\coprod_{k}\End_{\awebone[\levi]}(1^{\hcirc k})$ and 
$\coprod_{k}\End_{\fundone[\levi]}(\vecrep[1]^{\otimes k})$ are
isomorphic $\Kone$-vector spaces with $\apfunctorl{1}{\levi}$ inducing 
an isomorphism. This follows from \autoref{L:ProofLinInd} and 
\autoref{L:ProofBasis}. The proof is complete.
\end{proof}

\begin{proof}[Proof of \autoref{T:EquivalenceMain}.(b)]
\autoref{L:ProofWelldefined} shows that the functors $\apfunctor[1]$
and $\apfunctorl{1}{\levi}$
are well-defined and \autoref{L:EquivalenceMonoidal} shows the 
statements involving the pivotal structure.
Fully faithfulness follows from \autoref{P:ProofFullyFaithful}, and
\autoref{L:LeviRepSemi} ensures that we have that
$\completion{\fundone[\levi]}$ is equivalent to 
$\repone[\levi]$. These statements taken together complete the proof 
using the usual properties of the additive idempotent 
completion.
\end{proof}


\subsection{Proof for \texorpdfstring{$\qpar\neq 1$}{q not 1}}\label{SS:Proofq}


We start with the following lemma.

\begin{Lemma}\label{L:ProofSurjective}
The functor $\apfunctor$ is full.
\end{Lemma}

\begin{proof}
It is clear that the image of the crossings and the coils span. 
(Note that there is no issue with essential circles in $\fund[\levi]$.)
\end{proof}

For $\qpar\neq 1$, we note that we can mimic \autoref{D:ProofBasis} 
on the Levi side (the only difference is that we use a positive lift, in Levi crossings, of $\sigma$ 
instead of $\sigma$ itself) 
to define $\lbasis$ as well as $\lbasisset$.
Now we use that and \autoref{L:ProofSurjective} to define 
$\basis$ as well as $\basisset$ by pulling back the elements 
from $\lbasis$ by choosing a preimage.

\begin{Lemma}\label{L:ProofLinIndTwo}
The set $\basisset\subset\coprod_{k}\End_{\aweb[\levi]}(1^{\hcirc k})$ is $\K$-linearly 
independent. 
Moreover, the set $\lbasisset$ is $\K$-linearly 
independent in $\coprod_{k}\End_{\fund[\levi]}(\vecrep^{\otimes k})$.
\end{Lemma}

\begin{proof}
The claim on the Levi side can be proven {\ver} as in 
\autoref{L:ProofLinInd}, 
so our focus is on the web side. 
However, by construction, the set $\basisset$ is then 
sent to  
$\lbasisset$, so $\basisset$ is $\K$-linearly 
independent because $\lbasisset$ is.
\end{proof}

\begin{Lemma}\label{L:ProofSpanningTwo}
The set $\basisset$ $\K$-linearly
spans $\coprod_{k}\End_{\aweb[\levi]}(1^{\hcirc k})$. 
Moreover, the set $\lbasisset$ is a $\K$-linear 
spanning set of $\coprod_{k}\End_{\fund[\levi]}(\vecrep^{\otimes k})$.
\end{Lemma}

\begin{proof}
We can remove essential circles in front of the cylinder by 
definition of the monoidal product, see also \autoref{Eq:ProofLeftRight},
and the Levi ideal.
The essential circles in the back of the cylinder labeled by $1$ 
can then be inductively removed by using
\begin{gather*}
\begin{tikzpicture}[anchorbase,scale=1]
\draw[usual,directed=0.99] (1,0)node[below]{$1$} to (0,1)node[above]{$1$};
\draw[usual,directed=0.99,crossline] (0,0)node[below]{$1$} to (1,1)node[above]{$1$};
\end{tikzpicture}
-
\begin{tikzpicture}[anchorbase,scale=1]
\draw[usual,directed=0.99] (0,0)node[below]{$1$} to (1,1)node[above]{$1$};
\draw[usual,directed=0.99,crossline] (1,0)node[below]{$1$} to (0,1)node[above]{$1$};
\end{tikzpicture}
=
(\qpar-\qpar^{-1})\cdot
\begin{tikzpicture}[anchorbase,scale=1]
\draw[usual,directed=0.99] (1,0)node[below]{$1$} to (1,1)node[above]{$1$};
\draw[usual,directed=0.99] (0,0)node[below]{$1$} to (0,1)node[above]{$1$};
\end{tikzpicture}
,
\end{gather*}
which is the classical skein relation (that holds in our setting by using \autoref{Eq:PrelimThinCrossing}).
Using explosion, the rest of the argument is the same as in the proof of \autoref{L:ProofSpanning} 
by using that the $\qpar\neq 1$ basis agrees with $\basissetone$ on the associated graded 
by filtration by number of crossings (using the skein relations).
\end{proof}

\begin{Lemma}\label{L:ProofBasisTwo}
The set $\basisset$ is a $\K$-basis of 
$\coprod_{k}\End_{\aweb[\levi]}(1^{\hcirc k})$. Furthermore,
the set $\lbasisset$ is a $\K$-basis 
of $\coprod_{k}\End_{\fund[\levi]}(\vecrep^{\otimes k})$. 
\end{Lemma}

\begin{proof}
By \autoref{L:ProofLinIndTwo} and \autoref{L:ProofSpanningTwo}.
\end{proof}

\begin{Remark}\label{R:ProofBasisTwo}
We do not have or need any explicit description of 
the elements of $\basisset$ in terms of webs.
\end{Remark}

\begin{Proposition}\label{P:ProofFullyFaithfulTwo}
The functor $\apfunctorl{\qpar}{\levi}$ is fully faithful
and the Levi ideal $\lideal$
is the two-sided $\vcirc$-$\lprod$-ideal generated by the Levi evaluations.
Here $\lprod$ is the pullback of the monoidal structure from $\fund[\levi]$ to 
$\aweb$.
\end{Proposition}

\begin{proof}
As in the proof of 
\autoref{P:ProofFullyFaithful}.
\end{proof}

\begin{proof}[Proof of \autoref{T:EquivalenceMain}.(a)]
Using the above statements, this can be proven 
{\ver} as \autoref{T:EquivalenceMain}.(b).
\end{proof}

\end{document}